\documentclass[a4paper,11pt,oneside]{article}
\usepackage[dvipdfmx]{graphicx,color}
\usepackage{tikz}
\input{mymacros.sty}
\usepackage{pst-all}

 \newtheorem*{ack}{Acknowledgements}

\DeclareFontFamily{U}{mathx}{\hyphenchar\font45}
\DeclareFontShape{U}{mathx}{m}{n}{
      <5> <6> <7> <8> <9> <10>
      <10.95> <12> <14.4> <17.28> <20.74> <24.88>
      mathx10
      }{}
\DeclareSymbolFont{mathx}{U}{mathx}{m}{n}
\DeclareFontSubstitution{U}{mathx}{m}{n}
\DeclareMathAccent{\widecheck}{0}{mathx}{"71}

\newcommand{\Nk}{{\mathsf N}}
\def\cal{\mathop{\mathrm{AL}}\nolimits^{\vee}}
\def\al{\mathop{\mathrm{AL}}\nolimits}
\def\td{\mathop{\mathrm{Td}}\nolimits}

\def\fund{\mathop{\mathrm{fund}}\nolimits}
\def\adj{\mathop{\mathrm{adj}}\nolimits}
\def\euk{\mathop{\mathrm{Eu}}\nolimits_{K}}
\def\eukt{\mathop{\mathrm{Eu}}\nolimits_{K}^{t}}
\def\intk{\int^{[K]}}

\excludeversion{NB}

\begin{document}
\title{$K$-theoretic wall-crossing formulas and 
multiple basic
hypergeometric series}
\author{R. Ohkawa and J. Shiraishi}
\maketitle

\begin{abstract}
We study $K$-theoretic integrals over famed quiver moduli 
via wall-crossing phenomena.
We study the chainsaw quiver varieties, and 
consider generating functions defined by two types of 
$K$-theoretic classes.
In particular, we focus on integrals over the handsaw quiver 
varieties of type $A_{1}$, and get functional equations for each of them. 
We also give explicit formula for these partition functions. 
In particular, we obtain geometric interpretation of 
transformation formulas for multiple basic hypergeometric series
including the Kajihara transformation formula, and the
one studied by Langer-Schlosser-Warnaar 
and Halln\"as-Langman-Noumi-Rosengren. 
\end{abstract}

\section{Introduction}

\subsection{Framed quiver moduli}

Framed quiver moduli spaces
$M_{Q}^{\zeta}(\alpha)$ 
{\color{black}introduced by
King \cite{Ki} and Reineke \cite{R}
}
are moduli spaces of framed quiver representations.
We can identify them with moduli spaces of 
representations 
of a quiver {\color{black} 
$(Q_{0}, Q_{1})$} with 
{\color{black} a
set $Q_{2}$ of relations and} 
a distinguished vertex denoted by $\infty \in Q_{0}$, 
{\color{black} where $Q_{0}$ is the
set of vertices, and $Q_{1}$ is 
the set of arrows
{\color{black}(see \S \ref{subsec:setting} for 
more detailed explanation)}.
We call such a collection of a quiver
$(Q_{0}, Q_{1})$ with $Q_{2}, \infty$ a {\it framed quiver} $Q$.
We set $I = Q_{0} \setminus \lbrace
\infty \rbrace$.}
For a generic stability 
parameter $\zeta=(\zeta_{i})_{i \in I} \in
\R^{I}$ sitting
inside the complement of union of hyperplanes called 
{\it walls}, we consider corresponding moduli spaces
$M_{Q}^{\zeta}(\alpha)$ of 
$\zeta$-stable framed quiver representations where 
$\alpha=(\alpha_{v})_{v \in Q_{0}} \in 
\Z^{Q_{0}}$ is a dimension vector of a fixed underlying 
graded vector space of $Q$-representations.
{\color{black} For each vertex $v$ of 
	$Q_{0}$, we have a {\it tautological bundle} $\mc V_{v}$
over $M_{Q}^{\zeta}( \alpha )$, and 
a natural homomorphism of vector bundles $\phi_{h} \colon \mc V_{\out(h)} \to \mc V_{\inn(h)}$ corresponding to each arrow 
$h \in Q_{1}$.  
More presicely, this means that the restriction
$\lbrace \phi_{h} |_{x} \colon \mc V_{\out(h)}|_{x} 
\to \mc V_{\inn(h)}|_{x} \rbrace_{h \in Q_{1}}$ to each point $x \in M_{Q}^{\zeta}(\alpha)$
represents an isomorphism class $x$ of $Q$-representations}.

We call arrows departing from $\infty$ {\it framings}, and
arrows arriving at $\infty$ {\it co-framings}.
The natural diagonal torus $\mb T$ acts on moduli spaces 
$M_{Q}^{\zeta}(\alpha)$ 
associated with
framings and co-framings.
Consequently, we can consider
$\mb T$-equivariant $K$-theory 
over $M_{Q}^{\zeta}(\alpha)$.   
To define $\mb T$-equivariant integrals over 
$M_{Q}^{\zeta}(\alpha)$, 
we add other types of 
torus to $\mb T$ if necessary so that 
the affine algebro-geometric quotients 
$M^{0}_{Q}(\alpha)$
have the fixed points sets 
$M^{0}_{Q}(\alpha)^{\mb T} = \pt$, and
we have $\mb T$-equivariant morphisms 
from $M_{Q}^{\zeta}(\alpha)$ to $M^{0}_{Q}(\alpha)$.
For a $\mb T$-equivariant $K$-theoretic class $\Lambda$ 
on $M_{Q}^{\zeta}(\alpha)$, we define the $K$-theoretic
Euler class $\eukt(\Lambda)=\ch (\wedge_{-t} \Lambda^{\vee})$,
where $\Lambda^{\vee}$ is the dual of $\Lambda$,
$\wedge_{-t} \Lambda^{\vee}$ is the 
$t$-wedge product $\sum_{i \ge 0} (-t)^{i} \cdot 
\wedge^{i} \Lambda^{\vee}$, and $\ch$ denotes the Chern character
(see \S \ref{subsec:k-th} for the precise definition). 
When $t=1$, set $\euk(\Lambda)=\eukt(\Lambda)|_{t=1}$.
We define $K$-theoretic integrals 
$\intk_{[M^{\zeta}(\alpha)]^{vir}} \eukt( \Lambda)$ in 
\S \ref{subsec:k-th}.
{\color{black} We mainly study the case where $\Lambda$ is 
a linear combinations of $\mc V_{v}$ and 
$\mc V_{v}^{\vee}$ ($v \in Q_{0}$)
with the coefficients of $\mb T$-characters.
Then we write such a class by $\Lambda(\mc V)$.}
{\color{black} 
We introduce $t$-integer $[m]_{t} =(1-t^{m} ) / (1-t)$ for a non-negative integer $m$,
and we set $\lvert \mk I \rvert_{t} = 
[\lvert \mk I \rvert ]_{t}$ 
for a finite set $\mk I$.
We also set $[m]_{t} ! = [1]_{t} \cdot [2]_{t} \cdots [m]_{t}$.}

For two generic stability 
parameters $\zeta^{+}, \zeta^{-}$
inside adjacent chambers along a wall, 
we consider corresponding moduli spaces
$M_{Q}^{\zeta^{+}}(\alpha), M_{Q}^{\zeta^{-}}(\alpha)$ of 
$\zeta^{\pm}$-stable 
framed quiver representations.
A {\it wall-crossing formula} describes
the difference 
$\intk_{[M_{Q}^{\zeta^{+}}(\alpha)]^{vir}} \eukt (\Lambda)
-\intk_{[M_{Q}^{\zeta^{-}}(\alpha)]^{vir}} \eukt (\Lambda)$ {\color{black} for $\Lambda=\Lambda(\mc V)$ on 
$M_{Q}^{\zeta^{\pm}}(\alpha)$}. 
In the first part of the paper (until \S 3), 
we give a framework 
using enhanced master spaces
via Mochizuki's method \cite{M}, 
and deduce the wall-crossing formula
in Theorem \ref{thm:wcm}.
{\color{black} There are many articles that use Mochizuki-style wall-crossing 
(cf. \cite{KLT} for $K$-theoretic
Donaldson-Thomas theory, and references therein).}

{\color{black}
To state our results, we need a choice
of a vector $\beta=(\beta_{i})_{i \in I}
\in \Z^{I}$, 
a vertex $\ast \in I$, and 
a generic parameter $\bar{\zeta}$ 
such that $0 \le \beta_{i} \le \alpha_{i}$, 
$\beta_{\ast} \neq 0$, and 
$\bar{\zeta}$ belongs to the 
hyperplane $\beta^{\perp}$.
By the genericity assumption of $\bar{\zeta}$, 
we have only two chambers whose boundaries contain $\bar{\zeta}$
(see the beginning of 
\S \ref{subsec:enhancement} for
the presice choice of $\bar{\zeta}$).
Let $\mc C^{-}$ denotes one of them 
characterized by the condition $\sum_{i \in I } \zeta^{-}_{i} \beta_{i} < 0$ for any $\zeta^{-}=
(\zeta^{-}_{i})_{i \in I} \in \mc C^{-}$, and $\mc C^{+}$ the other one.

Let $\Dec{}_{\beta_{\ast}, j}^{\alpha_{\ast}}$ denote the set 
of collections $\mbi{\mk I}=
(\mk I^{(1)}, \ldots, \mk I^{(j)})$ 
(decomposition data)
such that
\begin{enumerate}
\item[$\bullet$] $\mk I^{(1)}, \ldots, \mk I^{(j)}$
are disjoint non-empty subsets of $[\alpha_{\ast}]=
\lbrace 1, 2, \ldots, \alpha_{\ast} \rbrace$ 
\item[$\bullet$] $|\mk I^{(i)}| = d^{(i)} 
\beta_{\ast}$ with $d^{(i)} \in \Z_{>0}$ for $i=1, \ldots,j$, and
\item[$\bullet$] $\min(\mk I^{(1)}) > \cdots > \min(\mk I^{(j)})$.
\end{enumerate}
We set $\mk I^{\infty} =
\lbrace 1, \ldots, \alpha_{\ast} \rbrace
\setminus \mk I^{(1)} \sqcup \cdots \sqcup 
\mk I^{(j)}$, 
$\mbi d_{\mbi{\mk I}}
=(d^{(1)}, \ldots, d^{(j)})$ and 
$\lvert \mbi d_{\mbi{\mk I}} \rvert 
= d^{(1)} + \cdots + d^{(j)}$.

We take stability parameters 
$\zeta^{\pm}=(\zeta^{\pm}_{i})_{i \in I}$ from the chambers $\mc C^{\pm}$ introduced
as above.
}

\begin{thm}
\label{thm:introwcm}
{\color{black} For $\Lambda=\Lambda(\mc V)$
on $M_{Q}^{\zeta^{\pm}}(\alpha)$,} 
we have
\begin{align}
&
\intk_{[M_{Q}^{\zeta^{+}}(\alpha)]^{vir}} 
\eukt (\Lambda(\mc V))  - 
\intk_{[M_{Q}^{\zeta^{-}}(\alpha )]^{vir}} 
\eukt (\Lambda(\mc V)) 
\notag
\\
&=
\nonumber
\sum_{j=1}^{\lfloor \alpha_{\ast}/\beta_{\ast} \rfloor } 
\sum_{\mbi{\mk I} \in \Dec{}_{\beta_{\ast}, j}^{\alpha_{\ast}}} 
{ |\mk I^{\infty}|_{t}! \over [\alpha_{\ast}]_{t}! }
\Res_{u_{1}=0, \infty} \cdots 
\Res_{u_{j}=0, \infty}
\int_{[\wt M^{0} (\alpha - 
	\lvert \mbi d_{\mbi{\mk I}} \rvert 
\beta )]^{vir}} 
\wt C_{\mbi{\mk I}} (\mc V)
\end{align}
{\color{black}
where $\lfloor \alpha_{\ast} / \beta_{\ast} 
\rfloor$ is a round down of the rational number
$\alpha_{\ast} / \beta_{\ast}$, and 
\nonumber
	$\wt C_{\mbi{\mk I}} (\mc V)$
is a localized 
$\mb T \times \prod_{k=1}^{j} 
\C^{\ast}_{\hbar_{k}}$-equivariant 
cohomology class defined in 
\eqref{itcoho}.
The integral $\int_{[\wt M^{0} (\alpha - 
	\lvert \mbi d_{\mbi{\mk I}} \rvert 
\beta )]^{vir}} 
\wt C_{\mbi{\mk I}} (\mc V)
$ is a rational expression of
$u_{1} = e^{\hbar_{1}}, \ldots, 
u_{j}=e^{\hbar_{j}}$ 
with the coefficients of the $\mb T$-equivariant 
parameters.
The symbol $\displaystyle \Res_{u = 0,\infty}$ denotes 
taking a sum of residues at $u=0$ and $u=\infty$.
}
\end{thm}

{\color{black}
We take $\mb T$-equivariant 
$K$-theory classes 
$\Lambda_{\adj}=\Lambda_{\adj}(\mc V)$
defined by a 
linear combination 
of $\mc V_{v}$ and $\mc V_{v}^{\vee}$ 
($v \in Q_{0}$)
corresponding to 
$T^{\ast} M^{\zeta^{\pm}}_{Q}(\alpha)$ on
$M^{\zeta^{\pm}}_{Q}(\alpha)$.}
Then we get a more concrete formula in Theorem \ref{thm:introadj} 
(Theorem \ref{thm:adjoint} in the main
body of the text).
{\color{black}
We introduce the Euler form
\begin{align}
\label{eulerform}
\chi(\alpha', \beta')
=
\sum_{h \in Q_{1}} 
\alpha'_{\out(h)} \beta'_{\inn(h)}
- 
\sum_{ \gamma \in Q_{2}} \alpha'_{\out(\gamma)} 
\beta'_{\inn(\gamma)}
-
\sum_{i \in I} 
\alpha'_{i} \beta'_{i}
\end{align}
for dimension vectors 
$\alpha'=(\alpha'_{v})_{v \in Q_{0}},  
\beta'=(\beta'_{v})_{v \in Q_{0}} \in \Z^{Q_{0}}$.}
\begin{thm}
\label{thm:introadj}
{\color{black}
For $\Lambda_{\adj}=\Lambda_{\adj}(\mc V)$ on 
$M_{Q}^{\zeta^{\pm}}(\alpha)$,} we have
\begin{align*}
&
\intk_{[M_{Q}^{\zeta^{+}}(\alpha)]^{vir}} 
\euk^{t}(\Lambda_{\adj})
- \intk_{[M_{Q}^{\zeta^{-}}(\alpha)]^{vir}} 
\euk^{t}(\Lambda_{\adj})\\
&=
\sum_{k=1}^{\lfloor \alpha_{\ast}/\beta_{\ast} \rfloor}
{\color{black}
\sum_{j=1}^{k} 
\sum_{\mbi{\mk I} \in \Dec^{\alpha_{\ast}}_{\beta_{\ast}, j}
\atop
|\mbi d_{\mbi{\mk I}}|=k}
} 
{ |\mk I^{\infty}|_{t} !  
\over 
[\alpha_{\ast}]_{t} ! }
\prod_{i=1}^{j} 
{[d^{(i)} \beta_{\ast} -1]_{t}!
\over
t-1} 
\gamma_{d^{(i)}} (
{\color{black}t}) 
\\
& \cdot  
( t^{s( \mk I^{(i)}, \mbi{\mk I}^{>i})+ {\color{black}
\chi \left( d^{(i)} \beta,
\alpha- d^{\le i} \beta \right)
} 
} 
- t^{s( \mbi{\mk I}^{>i},  \mk I^{(i)}) + {\color{black} 
\chi \left( 
\alpha- d^{\le i} \beta,
d^{(i)} \beta \right)
}
})  
\int_{ [M^{\zeta^{-}}(\alpha - k \beta )]^{vir}} 
\euk^{t}(\Lambda_{\adj}),
\end{align*}
{\color{black} where
$d^{\le i} = d^{(1)} + \cdots + d^{(i)}$ for 
$\mbi d_{\mbi{\mk I}} =(d^{(1)}, \ldots,  d^{(j)})$, 
$\Dec (\alpha_{\ast})= \bigsqcup_{j=1}
^{\lfloor \alpha_{\ast}/\beta_{\ast} \rfloor} 
\Dec{}_{\beta_{\ast}, j}^{\alpha_{\ast}}$,
$\gamma_{d^{(i)}} ({\color{black}t})$, and
{\color{black}
$s(\mk I, \mk I')$
}
are defined in 
\eqref{gamma},
and \eqref{sii}.
} 
\end{thm}

\subsection{Chainsaw quiver variety}
  
In the latter part of the paper, we study the chainsaw 
quiver 
\begin{center}
\includegraphics[scale=1]{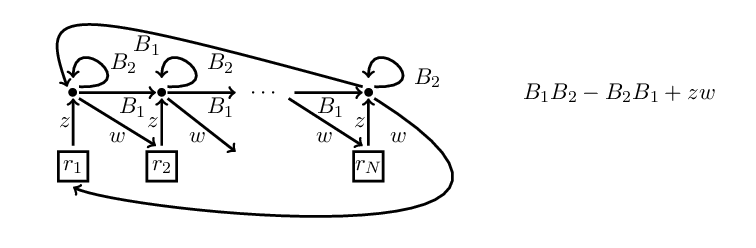}
\end{center}
where we index the vertices by $I=\lbrace 1, \ldots, N
\rbrace=\Z/N\Z$, and
take $r_{i}$ framings and $r_{i+1}$ co-framings 
for each vertex $i \in I$. 
{\color{black} The relation $B_{1}B_{2}-B_{2}B_{1} +zw=0$ is 
imposed at each vertex $i \in I$}.
The associated framed quiver moduli spaces
$M^{\zeta}(\mbi{r}, \mbi{v})$ are called 
the {\it chainsaw quiver variety}
of type $A^{(1)}_{N-1}$ where
$\mbi{r} = (r_{i})_{i \in I}$ is the framing, 
and $\mbi{v}=(v_{i})_{i \in I}$ is the dimension vector 
of a fixed underlying graded vector space.

We take two types of $K$-theory classes $\Lambda_{\adj} =
T^{\ast}M^{\zeta}(\mbi{r}, \mbi{v})$ and  
$\Lambda_{\fund}$ (see \eqref{ladj} and \eqref{lfund}
as to the definitions in terms of virtual vector bundles
using tautological bundles).
We take two stability parameters 
$\zeta^{\al}=(-1, \ldots, -1), 
\zeta^{\cal}=(1, \ldots, 1) \in \mb R^{I}$, and 
call $M^{\al} (\mbi{r}, \mbi{v})=M^{\zeta^{\al}}(\mbi{r}, \mbi{v})$
the {\it affine Laumon space} of type $A^{(1)}_{N-1}$, and
$M^{\cal} (\mbi{r}, \mbi{v})=M^{\zeta^{\cal}}(\mbi{r}, \mbi{v})$
the {\it dual} of $M^{\al} (\mbi{r}, \mbi{v})$.
Notice that there may be many walls separating $\zeta^{\al}$
and $\zeta^{\cal}$ in general.

We call the generating series 
\begin{align*}
Z_{\adj}^{\mbi{r}}
&=
Z_{\adj}^{\mbi{r}}( \mbi e, t| \mbi p | q,\kappa)
= \sum_{\mbi{v} \in \Z_{\ge 0}^{I}} 
\mbi p^{\mbi{v}} \intk_{[M^{\al}(\mbi{r}, \mbi{v})]^{vir}} 
\eukt (\Lambda_{\adj}), 
\\
\widecheck{Z}_{\adj}^{\mbi{r}}
&=
\widecheck{Z}_{\adj}^{\mbi{r}}( \mbi e, t| \mbi p | q,\kappa)
= \sum_{\mbi{v} \in \Z_{\ge 0}^{I}} 
\mbi p^{\mbi{v}} 
\intk_{[M^{\cal}(\mbi{r}, \mbi{v})]^{vir}} 
\eukt (\Lambda_{\adj}), 
\\
Z_{\fund}^{\mbi{r}}
&=
Z_{\fund}^{\mbi{r}}(\mbi e, \mbi \mu, 
\mbi \nu | \mbi p | q, \kappa)
= \sum_{\mbi{v} \in \Z_{\ge 0}^{I}} 
\mbi p^{\mbi{v}} \intk_{[M^{\al}(\mbi{r}, \mbi{v})]^{vir}} 
\euk (\Lambda_{\fund}),
\\
\widecheck{Z}_{\fund}^{\mbi{r}}
&=
\widecheck{Z}_{\fund}^{\mbi{r}}(\mbi e, \mbi \mu, 
\mbi \nu | \mbi p | q, \kappa)
= \sum_{\mbi{v} \in \Z_{\ge 0}^{I}} 
\mbi p^{\mbi{v}} 
\intk_{[M^{\cal}(\mbi{r}, \mbi{v})]^{vir}} 
\euk (\Lambda_{\fund})
\end{align*}
{\it partition functions}, where $q, \kappa, t, \mbi e=
\left( 
(e_{(i,\alpha)})_{\alpha=1}^{r_{i}} 
\right)_{i \in I}$,
$\mbi \mu=
\left(
(\mu_{(i, \alpha)})_{\alpha=1}^{r_{i}} 
\right)_{i \in I}$, and 
$\mbi \nu=
\left(
(\nu_{(i, \alpha)})_{\alpha=1}^{r_{i}} 
\right)_{i \in I}$
are
$\mb T$-equivariant parameters.
Remark that we have substituted $t=1$ 
in the latter series $
Z_{\fund}^{\mbi{r}}$ and 
$\widecheck{Z}_{\fund}^{\mbi{r}}$.

Originally, Nekrasov \cite{Nek} introduced
similar partition functions for $N=1$ and conjectured
that these partition functions give 
deformations of the Seiberg-Witten prepotentials for 
$N=2$ SUSY Yang-Mills theory.
This conjecture is proven in Braverman-Etingof \cite{BE}, Nekrasov-Okounkov \cite{NO} 
and Nakajima-Yoshioka \cite{NY1} independently.

{\color{black}
Combinatorial descriptions of 
torus fixed points of framed moduli spaces may have 
different appearance according to a choice of stability parameters. 
Neverthless they sometimes have the same 
torus equivariant integrals
as in Conjecture 
\ref{conj:adjintro} (2) below.
This turns out to be the
case for integrals of the adjoint matter classes 
over the Nakajima quiver varieties
(cf. \cite{AKL} ). 
They are studied in \cite{AO}, \cite{Ok}, \cite{OS}, and wall-crossing phenomena are
applied in different context from 
the present paper.
On the other hand, 
we study the $qq$-Painlv\'e
VI equation and its relationship with
the torus equivariant integrals over the
affine Laumon spaces in \cite{AHKOSSY1},
\cite{AHKOSSY2}.
It may be worth investigating this relation and possible representation
theoretic backgroud
in future.
}

We give some conjectures concerning relationships among
the dual partition functions for each type of
$K$-theory classes.
{\color{black} We define several infinite products
(see \cite{GR})
\begin{align*}
(x;q)_{\infty}
&=
\prod_{l=0}^{\infty} (1- xq^{l})
\\
(x;q, t)_{\infty}
&=
\prod_{l, m=0}^{\infty} (1- xq^{l} t^{m})
\\
(x;q, \kappa, t)_{\infty}
&=
\prod_{l,m,n=0}^{\infty} (1- xq^{l} \kappa^{m}
t^{n}).
\end{align*} 
}

\begin{conj}
\label{conj:adjintro}
(1)
For $\mbi \e_{l} = (\delta_{il})_{i \in I} \in \Z^{I}$, 
we have
\begin{align*}
\widecheck{Z}^{\mbi \e_{l}}_{\adj}(\mbi e, t|\mbi p|q, \kappa)
&=
\prod_{k=1}^{N-1}
{(qp_{l} \cdots p_{l+k-1}; q, t p_{1} \cdots p_{N})_{\infty}
\over
(tp_{l} \cdots p_{l+k-1}; q, t p_{1} \cdots p_{N})_{\infty}}
\\
& \cdot
{(qtp_{1} \cdots p_{N}; q, \kappa^{N}, t p_{1} \cdots p_{N})_{\infty}
(\kappa^{N} tp_1 \cdots p_{N}; q, \kappa^{N}, t p_{1} \cdots p_{N})_{\infty}
\over 
(t^{2}p_1 \cdots p_{N}; q, \kappa^{N}, t p_{1} \cdots p_{N})_{\infty}
(q \kappa^{N}p_1 \cdots p_{N}; q, \kappa^{N}, t p_{1} \cdots p_{N})_{\infty}},
\\
Z^{\mbi \e_{l}}_{\adj}(\mbi e, t|\mbi p|q, \kappa)
&=
\prod_{k=1}^{N-1} 
{(qp_{l-1} \cdots p_{l-k}; q, t p_{1} \cdots p_{N})_{\infty}
\over
(tp_{l-1} \cdots p_{l-k}; q, t p_{1} \cdots p_{N})_{\infty}}
\\
& \cdot
{(qtp_1 \cdots p_{N}; q, \kappa^{N}, t p_{1} \cdots p_{N})_{\infty}
(\kappa^{N} tp_1 \cdots p_{N}; q, \kappa^{N}, t p_{1} \cdots p_{N})_{\infty}
\over 
(t^{2}p_1 \cdots p_{N}; q, \kappa^{N}, t p_{1} \cdots p_{N})_{\infty}
(q \kappa^{N}p_1 \cdots p_{N}; q, \kappa^{N}, t p_{1} \cdots p_{N})_{\infty}}.
\end{align*}
(2)
When $r_{1} = \cdots = r_{N}=r$, we have 
\begin{align*}
Z^{\mbi r}_{\adj}(\mbi e, t|\mbi p|q, \kappa)
=
\widecheck{Z}^{\mbi{r}}_{\adj}(\mbi e, t|\mbi p|q, \kappa).
\end{align*}
\end{conj}
We describe behaviour under taking limits of 
parameters $\mbi e$ in Proposition \ref{prop:str},
from which we deduce functional equations
for arbitrary $\mbi r \in \Z^{I}$ from the Conjecture \ref{conj:adjintro}.
\begin{thm}
\label{thm:adj}
We assume Conjecture \ref{conj:adjintro}.
Then for arbitrary $\mbi r \in \Z^{I}$, we have
\begin{align}
\nonumber
{\widecheck{Z}^{\mbi{r}}_{\adj} (\mbi e, t|\mbi p|q, \kappa)
\over 
Z^{\mbi{r}}_{\adj} (\mbi e, t|\mbi p|q, \kappa)
}
=&
\prod_{k=1}^{N-1}
\prod_{l \in I}
{(qt^{r_{l+1} + \cdots + r_{l+k}}p_{l} \cdots p_{l+k-1}; 
q, t^{\lvert \mbi r \rvert} p_{1} \cdots p_{N}, t)_{\infty}
\over
(t^{r_{l+1} + \cdots + r_{l+k}+1}p_{l} \cdots p_{l+k-1}; 
q, t^{\lvert \mbi r \rvert} p_{1} \cdots p_{N}, t)_{\infty}}
\\
\label{adjrel}
\cdot &
{
(t^{r_{l} + \cdots + r_{l+k-1}+1}p_{l} \cdots p_{l+k-1}; 
q, t^{\lvert \mbi r \rvert} p_{1} \cdots p_{N}, t)_{\infty}
\over
(qt^{r_{l} + \cdots + r_{l+k-1}}p_{l} \cdots p_{l+k-1}; 
q, t^{\lvert \mbi r \rvert} p_{1} \cdots p_{N}, t)_{\infty}
}
\end{align} 
where $\lvert \mbi r \rvert =r_{1} + \cdots + r_{N}$.
\end{thm}
{\color{black}
In recent preprint \cite{AKL}, Conjecture 
\ref{conj:adjintro} is applied to study Vafa-Witten
invariants introduced by \cite{TT1} and \cite{TT2}.
}

For the fundamental matter class, we give functional 
equations including ``$q$ to the minus Laplacian'' 
$q^{-\Delta}$ defined in \eqref{qlaplacian}.
This gives factorized form of $\widecheck{Z}^{\mbi r}_{\fund}
/ Z^{\mbi r}_{\fund}$ for $N=1,2$. 
\begin{conj}
\label{conj:fundintro}
We have
\begin{align}
\nonumber
&
\left( {e_1e_2 \cdots e_N  p_1 p_2 \cdots p_N
\over 
\nu_1\nu_2 \cdots \nu_N};\kappa^N \right)_\infty \cdot
\sum_{l_{1} \geq 0} \cdots \sum_{l_{N} \geq 0}
{q^{(l_{1} l_{2} + l_{2} l_{3} + \cdots + l_{N} l_{1})/2}
\over (q;q)_{l_{1}} (q;q)_{l_{2}} \cdots (q;q)_{l_{N}}}
\\
\nonumber
& \cdot
\left(-q^{1/2} e_2 p_1\over \nu_2\right)^{l_{1}} 
\left(-q^{1/2} e_3 p_2\over \nu_3\right)^{l_{2}}
\cdots
\left(-q^{1/2} e_1 p_{N} \over \nu_1\right)^{l_{N}} 
\\
\nonumber 
& \cdot
q^{\sum_{i \in I} (l_{i+1}- l_{i-1})\vartheta_i/2}
q^{-{1\over 2}\Delta}
\widecheck{Z}^{\mbi{r}}_{\fund}
(\mbi e, \mbi \mu, \mbi \nu|\mbi p| q, \kappa)\\
\nonumber
=&
\left({\mu_1\mu_2 \cdots \mu_N  p_1 p_2 \cdots p_N 
\over 
e_1 e_2 \cdots e_N};\kappa^N \right)_\infty \cdot 
\sum_{l_{1} \geq 0} \cdots \sum_{l_{N} \geq 0}
{q^{(l_{1} l_{2} + l_{2} l_{3} + \cdots + l_{N} l_{1} )/2}
\over 
(q;q)_{l_{1}} (q;q)_{l_{2}} \cdots (q;q)_{l_{N}}}
\\
\nonumber
&\cdot
\left(-q^{1/2} \mu_1 p_1\over e_1\right)^{l_{1}} 
\left(-q^{1/2} \mu_2 p_2\over e_2\right)^{l_{2}} 
\cdots
\left(-q^{1/2} \mu_N p_N \over e_N\right)^{l_{N}} 
\\
\nonumber
& \cdot
q^{\sum_{i \in I} (-l_{i+1} + l_{i-1}) \vartheta_i/2}
q^{-{1\over 2}\Delta}
Z^{\mbi{r}}_{\fund}
(\mbi e, \mbi \mu, \mbi \nu|\mbi p| q, \kappa),
\end{align}
where we use symbols $e_{i}=
\prod_{\alpha=1}^{r_{i}} e_{(i, \alpha)}, \mu_{i}=
\prod_{\alpha=1}^{r_{i}} \mu_{(i, \alpha)}, 
\nu_{i}=\prod_{\alpha=1}^{r_{i}} \nu_{(i, \alpha)}$
for $i=1, \ldots, N$, and 
{\color{black}$\vartheta_{i}= p_{i} \partial/ \partial p_{i}$
for $i \in I=\mathbb{Z} /N \mathbb{Z}$.}
\end{conj}

\subsection{Handsaw quiver variety of type $A_{1}$}
We consider the case where $N=2$ and $\mbi{v}=(0,n)$
to get the {\it Laumon space} $M^{-} (\mbi{r}, n)=
M^{\al}(\mbi{r}, (0,n))$ of type $A_{1}$, and the dual
$M^{+} (\mbi{r}, n)=M^{\cal}(\mbi{r}, (0,n))$.
\begin{center}
\includegraphics[scale=1]{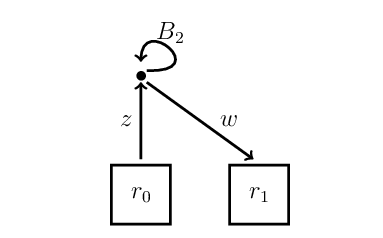}
\end{center}
Here we 
use $r_{0}=r_{2}$ and $p=p_{2}$, and substitute 
$p_{1}=0$ and $\kappa=1$. 
Set $e_{\alpha} = e_{(0, \alpha)}, \mu_{\alpha} =
q/\mu_{(0,\alpha)}$ for $\alpha=1, \ldots, r_{0}$, 
$e_{r_{0}+\alpha} = e_{(1, \alpha)}$, 
$\mu_{r_{0} +\alpha} = 1/\nu_{(1, \alpha)}$ for 
$\alpha=1, \ldots, r_{1}$, and
\begin{align*}
Z_{\pm \adj}^{\mbi{r}}
&=
Z_{\pm \adj}^{\mbi{r}}(q,t, \mbi e; p)
= \sum_{n=0}^{\infty} 
p^{n} \intk_{[M^{\pm}(\mbi{r}, n)]^{vir}} 
\eukt (\Lambda_{\adj}), 
\\
Z_{\pm \fund}^{\mbi{r}}
&=
Z_{\pm \fund}^{\mbi{r}}(q, \mbi e, \mbi \mu; p)
= \sum_{n=0}^{\infty} 
p^{n} \intk_{[M^{\pm}(\mbi{r}, n)]^{vir}} 
\euk (\Lambda_{\fund}).
\end{align*}

\begin{thm}
\label{thm:intromain1}
(a) We study
combinatorial description of torus fixed points, and 
derive the explicit formula (in Theorem \ref{explicit1}) 
for $Z^{\mbi{r}}_{\pm \adj}$
\begin{align}
\label{introadj+}
Z_{+\adj}^{\mbi{r}}(q,t, \mbi e; p)
&=
\sum_{ \substack{\mbi k \in \mathbb{Z}_{\ge 0} ^{r_{1}} } }
\left( t^{r_{0}} p  \right)^{\lvert \mbi k \rvert}
\prod_{ \substack{\alpha \in J_{1}\\ \beta \in J_{1}}}
\frac{(tq^{-k_{\beta}} e_{\alpha} / e_{\beta} ;q)_{k_{\alpha}}}
{(q^{-k_{\beta}} e_{\alpha} / e_{\beta} ;q)_{k_{\alpha}}} 
\prod_{\substack{ \beta \in J_{1} \\ \alpha \in J_{0}}} 
\frac{(q e_{\beta}/t e_{\alpha}  ;q)_{k_{\beta}}}
{( qe_{\beta}/e_{\alpha};q)_{k_{\beta}}},
\\
\label{introadj-}
Z_{-\adj}^{\mbi{r}}(q,t, \mbi e; p)
&=
\sum_{ \substack{\mbi k \in \mathbb{Z}_{\ge 0}^{r_{0}} } }
\left( t^{r_{1}} p \right)^{\lvert \mbi k \rvert}
\prod_{ \substack{\alpha \in J_{0}\\ \beta \in J_{0}}}
\frac{(tq^{-k_{\beta}} e_{\beta} / e_{\alpha} ;q)_{k_{\alpha}}}
{(q^{-k_{\beta}} e_{\beta} / e_{\alpha} ;q)_{k_{\alpha}}} 
\prod_{\substack{ \beta \in J_{0} \\ \alpha \in J_{1}}} 
\frac{( qe_{\alpha}/t e_{\beta} ;q)_{k_{\beta}}}
{( qe_{\alpha}/ e_{\beta} ;q)_{k_{\beta}}
},
\end{align}
where $J_{0}=\lbrace 1, \ldots, r_{0} \rbrace$ and 
$J_{1}=\lbrace r_{0}+1 ,\ldots, r_{0}+r_{1} \rbrace$.\\
(b)
By the wall-crossing formula Theorem \ref{thm:adjoint}, 
we get functional equation
(in Theorem \ref{thm:main1})
\begin{align}
\label{intromain1}
Z_{+ \adj}^{\mbi{r}} (q, t, \mbi e; p)
=
\frac{
(q t^{r_{0}} p ;q,t )_{\infty} (t^{r_{1}+1} p ;q,t )_{\infty} 
}{
(t^{r_{0}+1 } p ;q, t)_{\infty}(qt^{r_{1} } p ;q, t)_{\infty}
}
\cdot Z_{- \adj}^{\mbi{r}} (q, t, \mbi e; p),
\end{align}
where $(x;q, t)_{\infty}=
\prod_{m,n=0}^{\infty} (1- xq^{n} t^{m})$ 
is the double infinite product.
\end{thm}
Some remarks are in order concerning our proof of the wall-crossing
formula \eqref{intromain1} in Theorem 
\ref{thm:intromain1}. 
Recall that our argument uses enhanced master spaces via Mochizuki 
method summarized in Theorem \ref{thm:adjoint}.
One has an involved combinatorial structure
in Theorem \ref{thm:adjoint}
which prevents us from computing the coefficients 
explicitly in general.
However when $r_{0}=r_{1}=m$, we can show the vanishing of the
coefficients (Lemma \ref{lem:vanish}).
Then from Theorem \ref{thm:adjoint}, 
Theorem \ref{explicit1}, and Lemma \ref{lem:vanish} we have that
\begin{align}
\label{adjm}
Z_{+ \adj}^{(m,m)} (q, t, \mbi e; p)
=
Z_{- \adj}^{(m,m)} (q, t, \mbi e; p).
\end{align}
This transformation formula \eqref{adjm} has been 
already studied in the study of hypergeometric series
from various contexts.

\begin{thm}
[Langer-Schlosser-Warnaar \cite{LSW}, Halln\"as-Langmann-Noumi-Rosengren \cite{HLNR2}]
\label{thm:LSW}
By taking trigonometric limit of 
\cite[Cor. 4.3]{LSW}, or \cite[(6.7)]{HLNR2}, we have
\begin{align}
&
\sum_{ \substack{\mbi k \in \mathbb{Z}_{\ge 0} ^{m} } }
\left( q^{m} p / \tau^{m}  \right)^{\lvert \mbi k \rvert}
\prod_{\alpha, \beta=1}^{m}  
\frac{ \displaystyle
(q^{- k_{\beta} +1 } x_{\alpha} / \tau x_{\beta}; q)_{k_{\alpha}}}
{ \displaystyle
(q^{-k_{\beta}} x_{\alpha} /  x_{\beta}; q)_{k_{\alpha}}}
\frac{(x_{\beta} y_{\alpha}; q)_{k_{\beta}}}
{(qx_{\beta} y_{\alpha} / \tau ; q)_{k_{\beta} }}
\notag
\\
&=
\sum_{ \substack{\mbi k \in \mathbb{Z}_{\ge 0}^{m} } }
\left( q^{m} p / \tau^{m}  \right)^{\lvert \mbi k \rvert}
\prod_{\alpha, \beta=1}^{m}  
\frac{ \displaystyle
(q^{-k_{\beta}+1} y_{\alpha} / \tau y_{\beta}; q)_{k_{\alpha}}}
{ \displaystyle
(q^{-k_{\beta}} y_{\alpha} / y_{\beta}; q)_{ k_{\alpha}}}
\frac{(y_{\beta} x_{\alpha}; q)_{k_{\beta}}}
{(qy_{\beta} x_{\alpha} / \tau; q)_{k_{\beta} }}.
\label{noumi1}
\end{align}
\end{thm}
\begin{rem} 
One finds that the formula \eqref{noumi1} 
is identical to 
the formula \eqref{adjm}
after substituting $\tau=q/t$, 
$x_{\alpha} = e_{\alpha}$ for $\alpha \in J_{1}$, and 
$y_{\alpha} = q/te_{\alpha}$ for $\alpha \in J_{0}$
in \eqref{introadj+} and \eqref{introadj-}.
Therefore our geometric proof obtained in this article 
gives another interpretation
of \eqref{noumi1}. 
\end{rem}

A proof of \eqref{intromain1} in the general case where
$r_{0} \neq r_{1}$ is given by the 
following result due to M. Noumi (private communication).
\begin{thm}[Noumi]
\label{thm:noumi}
From Theorem \ref{thm:LSW}, when $r_{0} \le r_{1}$ 
it follows that
\begin{align}
&
\sum_{ \substack{\mbi k \in \mathbb{Z}_{\ge 0} ^{r_{1}} } }
\left( q^{r_{0}} p / \tau^{r_{0}}  \right)^{\lvert \mbi k \rvert}
\prod_{\alpha, \beta=1}^{r_{1}}  
\frac{ \displaystyle
(q^{- k_{\beta} +1 } x_{\alpha} / \tau x_{\beta}; q)_{k_{\alpha}}}
{ \displaystyle
(q^{-k_{\beta}} x_{\alpha} /  x_{\beta}; q)_{k_{\alpha}}}
\prod_{\beta=1}^{r_{1}}  
\prod_{\alpha=1}^{r_{0}}  
\frac{(x_{\beta} y_{\alpha}; q)_{k_{\beta}}}
{(qx_{\beta} y_{\alpha} / \tau ; q)_{k_{\beta} }}
\notag
\\
&=
\prod_{s=1}^{r_{1}-r_{0}} 
\frac{(q^{r_{0}+s}p/\tau^{r_{0}+s-1} ;q)_{\infty}}{(q^{r_{0}+s } p /
\tau^{r_{0}+s};q)_{\infty}}
\notag\\
&\cdot
\sum_{ \substack{\mbi k \in \mathbb{Z}_{\ge 0}^{r_{0}} } }
\left( q^{r_{1}} p / \tau^{r_{1}}  \right)^{\lvert \mbi k \rvert}
\prod_{\alpha, \beta=1}^{r_{0}}  
\frac{ \displaystyle
(q^{-k_{\beta}+1} y_{\alpha} / \tau y_{\beta}; q)_{k_{\alpha}}}
{ \displaystyle
(q^{-k_{\beta}} y_{\alpha} / y_{\beta}; q)_{ k_{\alpha}}}
\prod_{\beta=1}^{r_{0}}  
\prod_{\alpha=1}^{r_{1}}  
\frac{(y_{\beta} x_{\alpha}; q)_{k_{\beta}}}
{(qy_{\beta} x_{\alpha} / \tau; q)_{k_{\beta} }}.
\label{noumi2}
\end{align}
\end{thm}
\proof
We reproduce Noumi's proof.
Set
\[
F_{r_{1}, r_{0}}(x, y ; p)=
\sum_{ \substack{\mbi k \in \Z_{\ge 0} ^{r_{1}} } }
p^{|\mbi k| }
\prod_{\alpha, \beta=1}^{r_{1}}  
\frac{ \displaystyle
(q^{- k_{\beta} +1 } x_{\alpha} / \tau x_{\beta};q)_{k_{\alpha}}}
{ \displaystyle
(q^{-k_{\beta}} x_{\alpha} /  x_{\beta};q)_{k_{\alpha}}}
\prod_{\beta=1}^{r_{1}}  
\prod_{\alpha=1}^{r_{0}}  
\frac{(x_{\beta} y_{\alpha};q)_{k_{\beta}}}
{(qx_{\beta} y_{\alpha} / \tau;q)_{k_{\beta} }}.
\]
We need to compare 
$F_{r_{1}, r_{0}}(x,y;p)$ with $F_{r_{0}, r_{1}}(y,x;p)$.
From Theorem \ref{thm:LSW}, it follows that 
$F_{m,m}(x,y;p)=F_{m,m}(y,x;p)$.

For general $r_{0}, r_{1}$, we set $y'=(y_{1}, \ldots, y_{r_{0}-1})$ 
and take the limit $y_{r_{0}} \to 0$.
We have
\[
\displaystyle \lim_{y_{r_{0}} \to 0} 
F_{r_{1}, r_{0}}(x, y', y_{r_{0}}; p) = 
F_{r_{1}, r_{0}-1} (x, y'; p).
\]
On the other hand, we can arrange the factors as
\begin{align*}
F_{r_{0}, r_{1}}(y', y_{r_{0}}, x ; p)
&=
\sum_{ \substack{\mbi k \in \Z_{\ge 0} ^{r_{0}} } }
p^{|\mbi k'| }
\prod_{\alpha, \beta=1}^{r_{0} -1}  
\frac{ \displaystyle
(q^{- k_{\beta} +1 } y_{\alpha} / \tau y_{\beta};q)_{k_{\alpha}}}
{ \displaystyle
(q^{-k_{\beta}} y_{\alpha} /  y_{\beta};q)_{k_{\alpha}}}
\prod_{\beta=1}^{r_{0}-1}  
\prod_{\alpha=1}^{r_{1}}  
\frac{(y_{\beta} x_{\alpha};q)_{k_{\beta}}}
{(q y_{\beta} x_{\alpha} / \tau;q)_{k_{\beta} }}
\cdot
p^{k_{r_{0}}} \frac{(q^{-k_{r_{0}}+1}/ \tau;q)_{k_{r_{0}}}}
{(q^{-k_{r_{0}}};q)_{k_{r_{0}}}} 
\\
&\cdot \prod_{\alpha=1}^{r_{0} -1}
\frac{ \displaystyle
(q^{- k_{r_{0}} +1 } y_{\alpha} / \tau y_{r_{0}};q)_{k_{\alpha}}}
{ \displaystyle
(q^{-k_{r_{0}}} y_{\alpha} /  y_{r_{0}};q)_{k_{\alpha}}}
 \cdot
\prod_{\beta=1}^{r_{0}-1}
\frac{ \displaystyle
(q^{- k_{\beta} +1 } y_{r_{0}} / \tau y_{\beta};q)_{k_{r_{0}}}}
{ \displaystyle
(q^{-k_{\beta}} y_{r_{0}} /  y_{\beta};q)_{k_{r_{0}}}}
\cdot 
\prod_{\alpha=1}^{r_{1}}  
\frac{(y_{r_{0}} x_{\alpha};q)_{k_{r_{0}}}}
{(q y_{r_{0}} x_{\alpha} / \tau;q )_{k_{r_{0}} }},
\end{align*}
where $\mbi k'=(k_{1}, \ldots, k_{r_{0}-1})$ 
for $\mbi k=(k_{1}, \ldots, k_{r_{0}}) \in \Z_{\ge 0}^{r_{0}}$.
Hence we have 
\begin{align*}
\displaystyle \lim_{y_{r_{0}} \to 0}F_{r_{0}, r_{1}}
(y', y_{r_{0}}, x ; p)
&=
F_{r_{0}-1, r_{1}}(y', x ; qp/\tau) \cdot
\sum_{k_{r_{0}}=0}^{\infty} \frac{(\tau;q)_{k_{r_{0}}}}
{(q;q)_{k_{r_{0}}}} \left( \frac{qp}{\tau} \right)^{k_{r_{0}}} 
\\
&=
F_{r_{0}-1, r_{1}}(y', x ; qp/\tau ) \cdot
\frac{(qp;q)_{\infty}}{(qp/\tau; q)_{\infty}}. 
\end{align*}

For any $r_{0} \le r_{1}$,
we start from $F_{r_{1}, r_{1}}(x, y;p)=F_{r_{1}, r_{1}}(y, x;p)$
and get
\begin{align*}
F_{r_{1}, r_{0}}(x, y;p) 
&=
\lim_{y_{r_{0}+1} \to 0} \cdots \lim_{y_{r_{1}} \to 0} 
F_{r_{1}, r_{1} }(y, y_{r_{0}+1} , \ldots, y_{r_{1}}, x; p)
\\
&=
\prod_{s=1}^{r_{1} - r_{0}} \frac{(q^{s}p/ \tau^{s-1};q)_{\infty}}{(q^{s}p/ \tau^{s};q)_{\infty}} 
\cdot
F_{r_{0}, r_{1}}(y, x; q^{r_{1} - r_{0}}p/ \tau^{r_{1} - r_{0}}).
\end{align*}
\endproof

To write the explicit formula for 
$Z^{\pm}_{\fund}$,  
let $\mc B$ be the $q$-Borel transformation acting on $p$-series 
defined by
$\mc B \cdot p^{n} = q^{ n(n+1)/2} p^{n}$ for $n \in \Z_{\ge 0}$.
Then in Theorem \ref{subsec:explicit2}, we have 
\begin{align}
\label{introfund+}
\mc B^{-1} Z_{+ \fund}^{\mbi{r}}
&= 
\sum_{ \substack{\mbi k \in (\Z_{\ge 0})^{J_{1}} \\ |\mbi k|=n} }  
\left( -p \prod_{\alpha \in J_{0}} 
\frac{q}{ e_{\alpha} \mu_{\alpha} } \right)^{n}
\frac{\Delta_{J_{1}}(q^{\mbi k} \mbi e )}{\Delta_{J_{1}}(\mbi e)}
\prod_{ \substack{\alpha \in J_{1}\\ \beta \in J_{1}}}
\frac{(e_{\beta} \mu_{\alpha} ;q)_{k_{\beta}}}
{(qe_{\beta}/e_{\alpha} ;q)_{k_{\beta}}} 
\prod_{\substack{ \beta \in J_{1} \\ \alpha \in J_{0}}} 
\frac{( e_{\beta} \mu_{\alpha} ;q)_{k_{\beta}}}
{(qe_{\beta}/e_{\alpha} ;q)_{k_{\beta}}},
\\
\label{introfund-}
\mc B^{-1} Z_{- \fund}^{\mbi{r}}
&= 
\sum_{ \substack{\mbi k \in (\Z_{\ge 0})^{J_{0}} \\ |\mbi k|=n} }  
\left(p \prod_{\alpha \in J_{1}} e_{\alpha} \mu_{\alpha} 
\right)^{n}
\frac{\Delta_{J_{0}}(q^{\mbi k} \mbi e^{-1})}{\Delta_{J_{0}}(\mbi e^{-1})}
\prod_{ \substack{\alpha \in J_{0} \\ \beta \in J_{0}}}
\frac{
(q/ e_{\beta} \mu_{\alpha} ;q)_{k_{\beta}}
}
{
(qe_{\alpha}/e_{\beta} ;q)_{k_{\beta}} 
}
\prod_{ \substack{\beta \in J_{0} \\ \alpha \in J_{1}}} 
\frac{
(q/ e_{\beta} \mu_{\alpha} ;q )_{k_{\beta}}
}
{
(qe_{\alpha}/e_{\beta} ;q)_{k_{\beta}}
}.
\end{align}
\begin{thm}
\label{thm:intromain2}
Applying the wall-crossing formula Theorem \ref{thm:wcm}, 
we get functional equations
(in Theorem \ref{thm:main2}),
\begin{align}
\label{intromain2}
\mc B^{-1} Z_{+ \fund}^{\mbi{r}} (q, \mbi e, \mbi \mu; p)=
\frac{
( -e_{r_{0}+1} \mu_{r_{0}+1} \cdots e_{r} \mu_{r} p ;q)_{\infty} }
{
(- q^{r_{0}} p/e_{1}\mu_{1} \cdots e_{r_{0}}\mu_{r_{0}} ;q)_{\infty}
}
\cdot 
\mc B^{-1} Z_{- \fund}^{\mbi{r}} (q, \mbi e, \mbi \mu; p).
\end{align}
\end{thm}

We recall the Kajihara transformation \cite[(4.8)]{K} 
for multiple basic
hypergeometric series.
 following the formulation in \cite{HLNR1}.
Given four collections of variables
$(a_{1}, \ldots, a_{m}), 
\mbi x=(x_{1}, \ldots, x_{m}),  
 (b_{1}, \ldots, b_{n})$, and 
$(c_{1}, \ldots, c_{n})$,
Kajihara and Noumi's \cite{KN} multiple basic hypergeometric 
series are
\begin{align}
&
\phi^{m, n} \left( 
\begin{array}{c|c}
a_{1}, \ldots, a_{m} & b_{1}, \ldots, b_{n} \\ 
x_{1}, \ldots, x_{m} & c_{1}, \ldots, c_{n} 
\end{array} 
; p \right) 
\notag\\
&
=
\sum_{ \mbi k \in (\Z_{\ge 0})^{m} }  
p^{|\mbi k|} \frac{\Delta(q^{\mbi k} \mbi x)}{\Delta(\mbi x)}
\prod_{\alpha, \beta=1}^{m} 
\frac{\displaystyle
(a_{\alpha} x_{\beta} / x_{\alpha}; q)_{k_{\beta}}
}
{
(q x_{\beta} / x_{\alpha}; q)_{k_{\beta}}
}
\prod_{\beta =1}^{m} \prod_{\alpha =1}^{n}
\frac{
(x_{\beta} b_{\alpha} ; q)_{k_{\beta}}
}{
(x_{\beta} c_{\alpha}  ; q)_{k_{\beta}}
}.
\label{multhyp}
\end{align}
By \cite[Theorem 1.1]{K}, 
the following {\it Euler transformation formula} holds
for general parameters $a_{i}, b_{k}$, and $c$ : 
\begin{align}
\label{kajihara}
&
\phi^{m, n} \left( 
\begin{array}{c|c}
a_{1}, \ldots, a_{m} & b_{1}y_{1} , \ldots, b_{n} y_{n} \\ 
x_{1}, \ldots, x_{m} & c y_{1}, \ldots, c y_{n} 
\end{array} 
; u \right) 
\notag\\
& 
=
\frac{( abu/c^{n}; q)_{\infty}}{(u :q)_{\infty}} 
\cdot
\phi^{n, m} \left( 
\begin{array}{c|c}
c/b_{1}, \ldots, c/b_{n} & cx_{1}/a_{1}, \ldots, cx_{m}/a_{m} \\ 
y_{1}, \ldots, y_{n} & cx_{1}, \ldots, cx_{m} 
\end{array} 
; abu/c^{n} \right), 
\end{align}
where $a=a_{1} \cdots a_{m}$ and $b=b_{1} \cdots b_{n}$.
We substitute $m=r_{1}, n=r_{0}$, $u=-p/b$, $c=1$, and
$a_{\alpha}= e_{\alpha}\mu_{\alpha}, x_{\alpha}= qe_{\alpha}$ for
$\alpha \in J_{1}$, and 
$b_{\alpha}= e_{\alpha} \mu_{\alpha}/q, 
y_{\alpha}= 1/e_{\alpha}$ for $\alpha \in J_{0}$.
Then we see that \eqref{intromain2} and \eqref{kajihara} 
are identical.
Our proof using the wall-crossing formula is 
independent from the one in \cite{K}.
Recently, Yoshida \cite{Y} also deduce these transformation
formula from physics arguments.
 
This paper is organized as follows.
First we recall Mochizuki method to prove our main results.
In section 2, we introduce framed quiver moduli, 
the enhanced master space, and describe decomposition 
into
connected components of the fixed points set of the 
enhanced master space.
In section 3, we present wall-crossing formula to compute 
equivariant integrals over framed moduli recursively. 
In section 4, we introduce the chainsaw quiver variety,
give combinatorial descriptions of the fixed points
sets, and write explicit formula, and state 
Conjecture \ref{conj:adj} and Conjecture \ref{conj:fund} 
concerning relationships with the dual partition 
functions.
Then in section 5, we introduce the handsaw quiver variety of type $A_{1}$, 
and give an explicit formula.
In section 6, we evaluate residues in the wall-crossing 
formula 
in this setting, and give a proof of our main results
Theorem \ref{thm:main1} and Theorem \ref{thm:main2}.

\section{Framed quiver moduli and enhanced master space}
\label{sec:framed}

\subsection{Setting}
\label{subsec:setting}

\indent 
Let {\color{black} $(Q_{0}, Q_{1})$} be a quiver, i.e. a directed graph, where
$Q_{0}$ denotes the set of the vertices, 
$Q_{1}$ the set of the arrows.
For each arrow $h \in Q_{1}$, denote by 
$\text{out}(h) \slash \text{in}(h) \in Q_{0}$ 
the beginning$\slash$ending of $h$.
A path $p$ is a composable sequence $h_{1} h_{2} \cdots h_{k}$ 
of arrows satisfying $\out(h_{i})=\inn(h_{i+1})$ 
for $i=1, \ldots, k-1$.
For the path $p=h_{1} h_{2} \cdots h_{k}$, we set $\out(p)=\out(h_{k})$ and $\inn (p)=\inn(h_{1})$. 
Let ${\color{black} A (Q_{0}, Q_{1})}$ be the $\C$-vector space consisting of formal finite $\C$-linear combinations of paths.
An associative algebra structure on 
${\color{black} A (Q_{0}, Q_{1})}$ is induced by compositions of paths.
We call ${\color{black} A (Q_{0}, Q_{1})}$ 
a path algebra.

We consider a factored algebra 
${\color{black} A (Q_{0}, Q_{1})} / J$, 
where $J$ is 
a two-sided ideal generated by
some relations.
To be more precise, a relation is defined to be an 
element $\gamma=c_{1} p_{1} + \cdots + c_{m} p_{m} \in 
{\color{black} A (Q_{0}, Q_{1})}$
satisfying
$\out(p_{1}) = \cdots = \out(p_{m})$ and $\inn(p_{1}) = \cdots = \inn(p_{m})$.
For such $\gamma$, we set $\out(\gamma)=\out(p_{1})$ 
and $\inn(\gamma)=\inn(p_{1})$.
Let $Q_{2}=\lbrace \gamma_{1}, \ldots, \gamma_{s} 
\rbrace$ be a set of the relations.
Then set $J = \langle Q_{2} \rangle$ 
{\color{black} as a two-sided ideal}.
A collection $({\color{black} Q_{0}, Q_{1}}, Q_{2})$ is called a {\it quiver with 
relations}.

For geometric constructions, we take a collection 
$({\color{black} Q_{0}, Q_{1}}, Q_{2}, \infty)$ 
consisting of a quiver 
${\color{black} (Q_{0}, Q_{1})}$, 
relations $Q_{2}$, 
and a distinguished vertex $\infty \in Q_{0}$ 
called {\it a framing vertex}. 
We set $I=Q_{0} \setminus \lbrace \infty \rbrace$ and 
assume that $\out(\gamma), \inn(\gamma) \in I$ for any 
$\gamma \in Q_{2}$.
Such a collection ${\color{black} 
	Q=(Q_{0}, Q_{1}, Q_{2}, \infty)}$ is called a 
{\it framed quiver}.

Let $V=\bigoplus_{v \in Q_{0}} V_{v}$ be a 
finite dimensional $Q_{0}$-graded vector space
with $\dim V_{\infty} = 1$, and set
\begin{align}
\label{dimvect}
\alpha = ( \dim V_{v})_{v \in Q_{0}} \in \Z^{Q_{0}}.
\end{align}
Let
$\rep_{Q}(V)
= 
\lbrace
\rho \colon A (Q_{0}, Q_{1})/ \langle Q_{2} \rangle \to \End_{\C}(V) \mid
\text{ algebra homomorphisms }
\rbrace
$
be the set of $Q$-representations on $V$.
This is a zero set of a 
{\it \color{black} map} $\mu$ defined as follows.
We set
\begin{align}
\label{mbM}
\mb M_{Q}= \mb M_{Q}(V)=\prod_{h \in Q_{1}} \Hom_{\C} (V_{\out (h)}, V_{\inn (h)}), \quad 
\mb L_{Q}=\mb L_{Q}(V)=\prod_{\gamma \in R} \Hom_{\C} (V_{\out (\gamma)}, V_{\inn (\gamma)}), 
\end{align}
and define $\mu=\mu_{Q} \colon \mb M_{Q}(V) \to \mb L_{Q}(V)$ by sending
$B=(b_{h})_{h \in Q_{1}} \in \mb M_{Q}(V)$ to 
$\mu(B)=(b_{\gamma})_{\gamma \in Q_{2}}$ where 
$b_{\gamma} \in \Hom_{\C} (V_{\out (\gamma)}, 
V_{\inn (\gamma)})$ is defined from $B$ first by 
taking a composition along each path, and then making a 
linear combination of the contributions from the paths, 
according to $\gamma$.
We identify $\rep_{Q}(V)$ with $\mu^{-1}(0)$.
 
Recall that $I=Q_{0} \setminus \lbrace \infty \rbrace$.
For describing stability conditions, we need $(\zeta_{i} )_{i \in I} \in \mb R^{I}$ and 
\begin{align}
\label{zetainfty}
\zeta_{\infty} =
- \sum_{i \in I} \zeta_{i} \dim V_{i}. 
\end{align}
For a $Q_{0}$-graded subspace $S=\bigoplus_{v \in Q_{0}} S_{v}$ of $V$
with $S_{v} \subset V_{v}$, set
$\zeta (S)=\sum_{v \in Q_{0}} \zeta_{v} \cdot \dim S_{v}$.
Note that
$\zeta(V) =0$
always holds when $\dim V_{\infty}=1$ by \eqref{zetainfty}. 
In what follows, we only consider the case where a $Q_{0}$-graded vector space $V=\bigoplus_{v \in Q_{0}} V_{v}$ satisfies
the condition $\dim V_{\infty} =1$.
Given $Q$-representation 
{\color{black} $B=(b_{h})_{h \in Q_{1}}$} 
on $V$, when $b_{h} (S_{\out(h)}) \subset 
S_{\inn(h)}$ for any arrow $h$ of $Q$, we say that $S$ is a sub-representation of $B$. 
\begin{defn}
\label{def:stab}
A $Q$-representation $B$ on $V$ is said to be $\zeta$-semistable 
if for any sub-representation $S$ of $B$, we have $\zeta(S) \le 0$.
Furthermore $B$ is said to be $\zeta$-stable if $\zeta(S)<0$ for any non-trivial
proper sub-representation $S$ of $B$.
\end{defn}
For any $Q_{0}$-graded sub-space $S$ of $V$, we have $\zeta(S)+\zeta(V/S)=\zeta(V)=0$.
Hence we can replace $\zeta(S) \le 0$ ( resp. $\zeta(S) < 0$ ) with $\zeta(V/S) \ge 0$ ( resp. $\zeta(V/S) >0$ ) 
in the above definition.

Set $\Delta(\alpha)=\lbrace \beta=(\beta_{i})_{i \in I} \in \Z^{I} \mid 0 \le \beta_{i} \le \alpha_{i}=\dim V_{i} 
\ \text{ for any } i \in I, \text{ and } \beta \neq (0)_{i \in I} \rbrace$.
Then each element $\beta \in \Delta (\alpha)$ defines a hyperplane 
$\beta^{\perp}=\lbrace \zeta \in \mb R^{I} \mid \sum_{i \in I} \zeta_{i} \beta_{i} = 0 \rbrace$
called a {\it wall}.
A connected component of $\mb R^{I} \setminus \bigcup_{\beta \in \Delta (\alpha)} \beta^{\perp}$ is called
a {\it chamber} where stability and semi-stability in Definition \ref{def:stab} coincide.
When $\zeta$ belongs to a single chamber, the sign of $\zeta(S)$ does not change.
Hence $\zeta$-stability conditions for $\zeta$ in a single chamber $\mc C$ are all equivalent.
We sometimes say that a $Q$-representation $\rho$ is $\mc C$-semistable (resp. $\mc C$-stable) if 
$\rho$ is $\zeta$-semistable (resp. $\zeta$-stable) for $\zeta \in \mc C$.

Set $\mu^{-1}(0)^{\zeta}=\lbrace B \in \mu^{-1}(0) \colon
\zeta\text{-semistable} \rbrace$.
Then the group $\prod_{v \in Q_{0}} \GL(V_{v})$ acts 
on $\mu^{-1}(0)^{\zeta}$ naturally by conjugation.
Since we are considering the case where $V_{\infty} = \C$, 
the gauge group $G= \prod_{i \in I} \GL(V_{i})$ can be identified 
with $\prod_{v \in Q_{0}} \GL(V_{v})/\C^{\ast} \id_{V}$ by 
normalizing the component in $\GL(V_{\infty})$ to $\id_{V_{\infty}}$, 
noting that $\C^{\ast} \id_{V}$ trivially acts on $\rep_{Q}(V)$.
We set $M_{Q}^{\zeta}(V) = \mu^{-1}(0)^{\zeta} / G$ in the sense of 
GIT quotient.
Denoting the dimension vector by $\alpha = (\dim V_{v})_{v \in Q_{0}} \in (\Z_{\ge 0})^{Q_{0}}$, we also write
$M^{\zeta}(\alpha)=M_{Q}^{\zeta}(V)$ by abuse of notation.
We also set $M^{\mc C}(\alpha)=M^{\zeta}(\alpha)$ for $\zeta \in \mc C$.
Since elements in $G$ are regarded 
{\color{black}as }isomorphisms of framed quiver representations, which are 
isomorphisms of quiver representations with the component on 
framing vertex $\infty$ equal to the identity,  
$M_{Q}^{\zeta}(V)$ naturally regarded {\color{black}as }moduli of framed quiver representations.
This trick is used in \cite{C} to construct the Nakajima quiver variety
as a moduli space of quiver representations.

For each vertex $i \in I$, we consider a {\it tautological bundle} 
\[
\mc V_{i} = \left( \lbrace \rho \in \rep_{Q}(V) \mid 
\zeta\text{-semistable} \rbrace \times V_{i} \right) / 
G
\]
 over $M_{Q}^{\zeta}(V)$.
We have a natural homomorphism of vector bundles $\phi_{h} \colon \mc V_{\out(h)} \to \mc V_{\inn(h)}$ for each 
arrow $h \in Q_{1}$ such that 
$\phi_{h} |_{[B]} \colon \mc V_{\out(h)}|_{[B]} \to \mc V_{\inn(h)}|_{[B]}$ is 
identified with $b_{h} \colon V_{\out(h)} \to V_{\inn(h)}$.

\subsection{Enhancement of quiver}
\label{subsec:enhancement}
A dimension vector $\beta =(\beta_{i} )_{i \in I} \in (\mb Z_{\ge 0 })^{I}$ 
is said to be primitive if the greatest common divisor of $(\beta_{i})_{i \in I}$ is equal to 1.
We fix a primitive dimension vector $\beta
{\color{black}
\in \Delta(\alpha)}$.
We choose a generic parameter $\bar{\zeta} \in \beta^{\perp}$
such that $\bar{\zeta}$ does not lie on any hyperplane $(\beta')^{\perp}$ other than $\beta^{\perp}$ defined by
$\beta' \in \Delta(\alpha)
{\color{black}
\setminus \lbrace \beta \rbrace}$.
By this {\color{black}genericity  assumption for $\bar{\zeta}$}, we see that
$\bar{\zeta}(S)=0$ with $S_{\infty}=0$ implies that
$(\dim S_{i})_{i \in I} = d \beta \in (\Z_{\ge 0})^{I}$
for some integer $d \ge 0$.
{\color{black}
Furthermore we have only two chambers whose boundaries contain $\bar{\zeta}$.
Let $\mc C^{-}$ denotes one of them such that $(\zeta^{-}, \beta)<0$ for any 
$\zeta^{-} \in \mc C^{-}$, and $\mc C^{+}$ the other one.
}
\begin{lem}
\label{lem:cstab}
A $Q$-representation $(V, B)$ is $\mc C^{-}$-stable (resp. $\mc C^{+}$-stable) if and only if
the following two conditions hold :
\\
(a) $(V, B)$ is $\bar{\zeta}$-semistable.
\\
(b) We have no sub-representation $S$ of $(V,B)$ such that 
$\bar{\zeta}(S) =0$, $S_{\infty}=\C$ and $S \neq V$ (resp. $S_{\infty}=0$ and $S=0$). 
\end{lem}
\proof
We show that the negations of both sides are equivalent
by classifying non-trivial proper $Q_{0}$-graded subspace $S$ of $V$ such that $\zeta(S) > 0$.
When $\pm \bar{\zeta}(S)>0$, we have $\pm \zeta(S)>0$
since a stability parameter $\zeta$ inside $\mc C^{-}$ or $\mc C^{+}$ can converge to $\bar{\zeta}$ lying in the boundaries of $\mc C^{-}$ and $\mc C^{+}$.
When $\bar{\zeta}(S)=0$ and $S_{\infty}=\C$, we have $\zeta (S) = -\zeta(V/S) > 0 $ for $\zeta \in \mc C^{-}$ 
(resp. $\zeta (S) = -\zeta(V/S) < 0 $ for $\zeta \in \mc C^{+}$) 
since the dimension vector of $V/S$ is proportional to $\beta$.
Similarly when $\bar{\zeta}(S)=0$ and $S_{\infty}=0$, we have $\zeta (S) < 0 $ for $\zeta \in \mc C^{-}$ 
(resp. $\zeta (S) >0$ for $\zeta \in \mc C^{+}$).
Hence existence of a sub-representation $S$ of $(V, B)$ with $\zeta(S)>0$, 
meaning that $(V, B)$ is not $\zeta$-stable, is equivalent to saying that
(a) or (b) does not hold. 
\endproof

To describe wall-crossing from $\mc C^{-}$ to $\mc C^{+}$, we choose a vertex $i \in I$ such that $\beta_{i} \neq 0$.
Without loss of generality, such a fixed vertex is written by $\ast \in I$.
Then we define a new quiver 
{\color{black}$(\wt{Q}_{0}, \wt{Q}_{1})$} 
from {\color{black}$(Q_{0}, Q_{1})$} as follows.
The set of vertices in $\wt{Q}$ is 
\[
\wt{Q}_{0} = Q_{0} \sqcup \lbrace  \ast (k) \mid k= 1, 2, \ldots, L \rbrace,
\]
where $L \ge \alpha_{\ast} = \dim V_{\ast}$ is a positive integer.
The set $\wt{Q}_{1}$ of arrows in $\wt{Q}$ is a disjoint union of $Q_{1}$ and
\[ 
\lbrace \tilde{h}_{k} \colon \ast (k) \to \ast (k+1) \rbrace_{ k = 1}^{L - 1} \sqcup
\lbrace \tilde{h}_{L} \colon \ast (L) \to \ast \rbrace.
\]
We call such a quiver 
{\color{black} 
	$(\wt{Q}_{0}, \wt{Q}_{1})$} 
an {\it enhancement} of 
{\color{black} $(Q_{0},Q_{1})$}.

We also write by $\wt{Q}$ a framed quiver $({\color{black} 
\wt{Q}_{0}, \wt{Q}_{1}}, Q_{2}, \infty)$.
We consider a $\wt{Q}_{0}$-graded vector space $\wt{V}= V \oplus \bigoplus_{k=1}^{L} 
\wt{V}_{\ast (k)}$ with $\dim \wt{V}_{\ast (1)} \le 1$, $\dim \wt{V}_{\ast (L)} = \alpha_{\ast}=\dim V_{\ast}$, and
$0 \le \dim \wt{V}_{\ast (k)} 
-\dim \wt{V}_{\ast (k-1)} \le 1$
for $k=2, \ldots, L - 1$.
We set 
\begin{align}
\label{indexv}
\mk I=\lbrace k \in [L] \mid \dim \wt{V}_{\ast (k)} - \dim \wt{V}_{\ast (k-1)}= 1\rbrace,
\end{align}
where
$[L]=\lbrace 1,2,\ldots, L \rbrace$.

Let us consider $\wt{Q}$-representations 
$\wt{B}=(\wt{b}_{h})_{h \in \wt{Q}_{1}} \colon A (\wt{Q}_{0}, \wt{Q}_{1}) 
/ \langle Q_{2} \rangle \to \End_{\C}(\wt{V})$ on $\wt{V}$ satisfying the condition
that
linear maps 
$\wt{b}_{\tilde{h}_{k} } \colon \wt{V}_{\ast (k)} \to \wt{V}_{\ast (k+1)}$ for $k=1, \ldots, L-1$
and $\wt{b}_{\tilde{h}_{L}} \colon \wt{V}_{\ast (L)} \to V_{\ast}$ corresponding to the all added arrows are injective.
Such $\wt{Q}$-representations are called
{\it enhanced $Q$-representations}, which are regarded, up to $\prod_{k=1}^{L} \GL(V_{\ast(k)})$-action, 
as pairs $(B, F_{\bullet})$ of $Q$-representations $B$
and full flags $F_{\bullet}$ of $V_{\ast}$ satisfying 
\begin{align}
\label{index}
\lbrace i \in \Z \mid F_{i}/F_{i-1} \neq 0 \rbrace = \mk I
\end{align}
via taking image of compositions $V_{\ast(k)} \to V_{\ast(k+1)} \to \cdots \to V_{\ast(L)} \to V_{\ast}$ as $F_{k}$.
These elements $(B, F_{\bullet})$ forms a full flag bundle of $\mc V_{\ast}$ over $M^{\zeta}(\alpha)$. 

The stability parameters of $\wt{Q}$-representations 
are taken from the augmented space as
{\color{black} 
$\wt{\zeta}=(\zeta, \eta) \in \R^{I} \times \R^{L}$, 
where $\zeta= (\zeta_{i} )_{i \in I}
\in \R^{I}$, and 
$\eta=(\eta_{\ast(k)})_{k=1}^{L}$ belongs to 
the component $\R^{L}$ corresponding 
to the added vertices
$\ast(1), \ldots, \ast(L) \in \wt{Q}_{0}$.
We identify $\R^{I}$ with the subset $\R^{I} \times \lbrace 0 \rbrace \subset \R^{I} \times \R^{L}$.
We deform $\zeta^{-} = (\zeta^{-}, 0) \in \R^{I} \times \R^{L}$ to get a parameter $\wt\zeta=(\zeta, \eta)$ satisfying
conditions 
\cite[(9), (10)]{O5}, and find
distinguished chambers $\wt{\mc C}_{\ell} \subset \R^{I} \times \R^{L}$ containing
$\wt{\zeta}$ for $\ell=0, 1, \ldots, L$, 
which enhance the chamber $\mc C^{-}, \mc C^{+} \subset \R^{I}$ in the following sense.
}
Below, the $\wt{\mc C}_{\ell}$-stability conditions are interpreted in terms of enhanced $Q$ representations 
in Lemma \ref{lem:ellstab}, and
we see that these conditions connect $\mc C^{-}$-stability ($\ell=0$) and 
$\mc C^{+}$-stability ($\ell=L$)
in Lemma \ref{lem:flag}.
\begin{lem}
\label{lem:ellstab}
For $\ell \ge 0$, an enhanced $Q$-representation $\wt{B}=(B, F_{\bullet})$ is $\wt{\mc C}_{\ell}$-stable if $B$ is 
$\bar{\zeta}$-semistable and any sub-representation $S$ of $B$ with $\bar{\zeta}(S)=0$ satisfies the following two conditions:
\begin{enumerate}
\item[\textup{(a)}] If $S_{\infty}=\C$ and $S \neq V$, we have $F_{\ell} \not\subset S_{0}$. 
\item[\textup{(b)}] If $S_{\infty}=0$ and $S \neq 0$, we have $S_{0} \cap F_{\ell} =0$.
\end{enumerate}
\end{lem}
{\color{black}
	\proof 
See \cite[Definition 3.2, Proposition 3.3]{O5}.
\endproof
}
For a dimension vector $\alpha \in (\Z_{\ge 0})^{Q_{0}}$ and a finite subset
$\mk I \subset [L]=\lbrace 1,2, \ldots, L \rbrace$,
we write by $\wt{M}^{\ell}(\alpha, \mk I)$ moduli of $\wt{\mc C}_{\ell}$-stable enhanced
$Q$-representations $(B, F_{\bullet})$
on $V$ satisfying \eqref{dimvect} with full flags $F_{\bullet}$ of
$V_{\ast}$ satisfying \eqref{index}.
They are constructed as moduli $M^{\wt{\mc C}_{\ell}}_{\wt{Q}}(\wt{V})$ of 
$\wt{\mc C}_{\ell}$-stable $\wt{Q}$-representations 
on $\wt{V}=\bigoplus_{v \in Q_{0}} V_{v} \oplus \bigoplus_{k=1}^{L} V_{\ast(k)}$ with
\eqref{indexv}. 
\begin{lem}
\label{lem:flag}
When $\ell=0$ (resp. $\ell= L$), a pair $\wt{B}=(B, F_{\bullet})$ is $\wt{\mc C}_{\ell}$-stable if and only 
if $B$ is $\zeta$-stable for $\zeta \in \mc C^{-}$ (resp. $\zeta \in \mc C^{+}$).
\end{lem}
\proof
First we remark that $\zeta$-stability for $\zeta \in \mc C^{-}$, or $\zeta \in \mc C^{+}$
implies $\bar{\zeta}$-semi-stability since inside $\mc C^{-}$, or $\mc C^{+}$, we can take limit converging to $\bar{\zeta}$.
Under this assumption, $\zeta$-stability for $\zeta \in \mc C^{-}$ (resp. $\zeta \in \mc C^{+}$) 
is equivalent to prohibition of
sub-representation $S$ with $\bar{\zeta}(S)=0$ satisfying the assumption in (1) (resp. (2))
in Lemma \ref{lem:ellstab}. 
On the other hand, by the impossible conclusion $F_{0}=\lbrace 0 \rbrace \not \subset S_{0}$ when $\ell=0$ 
(resp. $S_{0} \cap F_{L} =S_{0} \cap V =0$ when $\ell=L$), the prohibition is also equivalent 
to the $\wt{\mc C}_{0}$-stability (resp. $\wt{\mc C}_{L}$-stability) condition
under the assumption of $\bar{\zeta}$-semi-stability.
\endproof
We need notation.
In general, for a vector bundle $\mc E$ over $M$ and a finite set $\mk I \subset \Z$ with $\rk \E$ elements, 
we write by $Fl_{M}(\mc E, \mk I)$ the full flag bundle over $M$ whose fiber over $x \in M$ consists of {\color{black} full flags}
$F_{\bullet}=(F_{i})_{i \in \Z}$ of $\mc E_{x}$ satisfying \eqref{index}.
These are all isomorphic to $Fl_{M}(\mc E, [\rk \mc E])$, but we use products of these flag bundles
and combinatorial description in terms of the indices later.
We always set $F_{0}=0$.
When $\ell =0$ (resp. $\ell=L$), we have $\wt{M}^{\ell}(\alpha, \mk I) \cong Fl_{M^{\zeta}(\alpha)}(\mc V_{\ast}, \mk I)$  for $\zeta \in \mc C^{-}$ 
(resp. $\zeta \in \mc C^{+}$) by this lemma.
Taking cap product of the Euler class of the relative tangent bundle of 
$Fl_{M^{\zeta}(\alpha)}(\mc V_{\ast}, \mk I) \to M^{\zeta}(\alpha)$ with the fundamental cycle
of $Fl_{M^{\zeta}(\alpha)}(\mc V_{\ast}, \mk I)$, 
integrals over $M^{\zeta}(\alpha)$
are reduced to ones over $Fl_{M^{\zeta}(\alpha)}(\mc V_{\ast}, \mk I)$.

\subsection{Enhanced master space}
\label{subsec:master}

Recall that we fix $0 \le \ell \le L$ and 
$\emptyset \neq \mk I \subset [L]=
\lbrace 1, \ldots, L \rbrace$.
We want to recursively compute integrals over 
$\wt M^{\ell}(\alpha, \mk I)$ in terms of ones over 
$\wt M^{0}(\alpha, \mk I)$
and $\wt M^{\ell^{\sharp}}(\alpha - d^{\sharp} \beta, 
\mk I \setminus \mk I^{\sharp})$ with 
$\ell^{\sharp} < \ell$ and 
non-empty subset $\mk I^{\sharp} \subset \mk I$
where $d^{\sharp}=|\mk I^{\sharp}|/\beta_{\ast}$
together with some coefficients defined by integral of cohomology classes (see \eqref{itcoho} below). 
To this end, we introduce the {\it enhanced master space} $\mc M=\mc M(\ell, \alpha, \mk I, \theta^{-}, \theta^{+})$ containing these 
moduli spaces (Definition \ref{defn:ems}).

We fix a $Q_{0}$-graded vector space $V=\bigoplus_{v \in Q_{0}} V_{v}$ satisfying
\eqref{dimvect},
and a $\wt{Q}_{0}$-graded vector space $\wt{V}=V \oplus \bigoplus_{k=1}^{L} V_{\ast(k)}$ satisfying \eqref{indexv}. 
We set $\mb M =\mb M_{Q} (V)$ and $\wt{\mb M} =\mb M_{\wt{Q}}(\wt{V})$.
The {\color{black} map} $\tilde{\mu} = \mu_{\wt{Q}} \colon \wt{\mb M} \to \mb L_{\wt{Q}}(\wt{V})=\mb L_{Q}(V)$
is obtained by composing the natural projection $\wt{\mb M} \to \mb M$ and $\mu=\mu_{Q} \colon \mb M \to \mb L_{Q}(V)$.
We also set $G=\prod_{i \in I} \GL(V_{i})$ and $\wt{G} = G \times \prod_{k=1}^{L} \GL(V_{k})$.
A moduli $\wt{M}^{\ell}(\alpha, \mk I)$ of $\wt{\mc C}_{\ell}$-stable enhanced $\wt{Q}$-representations on $\wt{V}$ is
constructed as a quotient stack $\wt{M}^{\ell}(\alpha, \mk I)=\tilde{\mu}^{-1}(0)^{\ell}/\wt G$,
where $\tilde{\mu}^{-1}(0)^{\ell}$ is $\wt{\mc C}_{\ell}$-stable locus in $\tilde{\mu}^{-1}(0)$.

Let us take suitable stability parameters $\theta^{+} \in \wt{\mc C}_{\ell} \cap \Z^{\wt{Q}_{0}}$ and
$\theta^{-} \in \wt{\mc C}_{0} \cap \Z^{\wt{Q}_{0}}$.
From $\theta^{\pm} =(\theta^{\pm}_{i})_{i \in \wt{Q}_{0}}\in \Z^{\wt{Q}_{0}}$, we define characters $\chi_{\pm} \colon \wt{G} \to \C^{\ast}$
by sending $g = ((g_{i})_{i \in I}, (g_{k})_{k=1}^{L} ) \in \wt{G}$ to
$\chi_{\pm}(g)=\prod_{i \in I} (\det g_{i} )^{\theta^{\pm}_{i}} 
\cdot \prod_{k=1}^{L} (\det g_{\ast(k)} )^{\theta^{\pm}_{\ast(k)}}$.
We define weight spaces $\C_{\chi_{\pm}}$ by $g \cdot 1=\chi_{\pm}(g)$ for $g \in \wt{G}$, 
and $\wt{G}$-equivariant line bundles $L^{\pm}$ over $\wt{\mb M}$ by $L^{\pm} = \wt{\mb M} \times \C_{\chi_{\pm}}$. 
Then we write by $\wh{\mb M}$ the projective bundle $\PP_{\wt{\mb M}}( L^{-} \oplus L^{+})$ over $\wt{\mb M}$, and 
consider the hyperplane bundle $\mo(1)$ (dual of the tautological bundle) over $\wh{\mb M}$.

We introduce the homogeneous coordinate $[x_{-}, x_{+}]$ of the fiber of 
the projective bundle $\wh{\mb M}=\PP_{\wt{\mb M}}(L^{-} \oplus L^{+})$.
Recall that in Definition \ref{def:stab} is given a numerical criterion of GIT-stability for $Q$-representation.
On $\hat{\mb M}$, GIT-stability is defined as follows. 
\begin{defn}
A point $\wh{B}=(\wt{B}, [x_{-}, x_{+}]) \in \wh{\mb M}$ with $\wt{B} \in \wt{\mb M}$ is said to be semi-stable if 
there exists a section $\sigma \in \Gamma(\wh{\mb M}, \mo(n))$ such that $\sigma(\wh{B}) \neq 0$ for a positive integer $n$,
where $\mo(n) = \mo(1)^{\otimes n}$. 
\end{defn}

We define a semi-stable locus $\wh{\mb M}^{ss}=\lbrace \wh{B} \in \wh{\mb M} \colon \text{ semi-stable } \rbrace$.
We consider the composition $\hat{\mu} \colon \wh{\mb M} \to \mb L_{Q}(V)$ of the projection $\wh{\mb M} \to \wt{\mb M}$ and $\tilde{\mu}
\colon \wt{\mb M} \to \mb L_{Q}(V)$.
Set $\hat{\mu}^{-1}(0)^{ss}=\hat{\mu}^{-1}(0) \cap \wh{\mb M}^{ss}$.
\begin{defn}
\label{defn:ems}
The quotient stack $\mc M =\mc M(\ell, \alpha, \mk I, \theta^{-}, \theta^{+})= \hat{\mu}^{-1}(0)^{ss} / \wt{G}$ is called the
enhanced master space. 
\end{defn}
{\color{black}
The enhanced master space $\mc M$ has non-trivial 
stabilizer groups as we remark
below Corollary \ref{cor:setbij}.
But we can sometimes take  
stability parmeters $\theta^{\pm}$ such
that $\mc M$ is a variety with the 
only trivial stabilizer 
group as in \cite{O3} where 
moduli space is empty set on one side of two chambers
along a wall.} 
We consider
the family $\theta^{t}=t \theta^{+} + (1-t) \theta^{-}$ of stability parameters for $0 \le t \le 1$. 
We have the following description of points in $\hat{\mu}^{-1}(0)^{ss}$.
\begin{lem}
\label{lem:ss}
For $\wh{B}=( \wt{B}, [x_{-}, x_{+}] ) \in \hat{\mu}^{-1}(0)$, we have the following.\\
(a) When $x_{-}=0,$ then $\wh{B}$ is semi-stable if and only if $\wt{B}$ is $\theta^{+}$-semistable. \\
(b) When $x_{+}=0$, then $\wh{B}$ is semi-stable if and only if $\wt{B}$ is $\theta^{-}$-semistable. \\
(c) When $x_{-} \cdot x_{+} \neq 0$, then $\wh{B}$ is semi-stable if and only if  $\wt{B}$ is  $\theta^{t}$-semistable for some $0 \le t \le 1$. 
\end{lem}

By (a) and (b) in Lemma \ref{lem:ss}, we have the inclusion of 
$\wt M^{\ell}(\alpha, \mk I)$ and $\wt M^{0}(\alpha, \mk I)$ into $\mc M$
as the hyperplane locus $\mc M_{\pm} = \lbrace x_{\mp}=0 \rbrace$. 
These spaces $\mc M_{+} \cong \wt M^{\ell}(\alpha, \mk I)$ and 
$\mc M_{-} \cong \wt M^{0}(\alpha, \mk I)$ are connected components of 
fixed points set $\mc M^{\C^{\ast}_{\hbar}}$ for 
the fiber-wise $\C^{\ast}_{\hbar}$-action on $\mc M$ induced by 
\begin{align}
\label{fiberwise}
e^{\hbar} \cdot [x_{-}, x_{+}] = [e^{\hbar} x_{-}, x_{+}]
\end{align}
for $e^{\hbar} \in \C^{\ast}_{\hbar}$.

We consider $\mc M_{exc} = \mc M^{\C^{\ast}_{\hbar}} \setminus (\mc M_{+} \sqcup \mc M_{-})$.
A point of the quotient stack $\mc M$ can be regarded as the orbit 
$\wt{G} \cdot \wh{B}$, where the $\wt{G}$-action is defined by
$g \cdot (\wt{B}, [x_{-}, x_{+}] )= (g \cdot \wt{B}, [\chi_{-}(g) x_{-}, \chi_{+}(g) x_{+}])$ for 
$g \in \wt{G}$.
This orbit is fixed by the fiberwise $\C^{\ast}_{\hbar}$-action if and only if $(\wt{B}, [e^{\hbar} x_{-}, x_{+}] )$ is contained in
$\wt{G} \cdot (\wt{B}, [x_{-}, x_{+}] )$.
When $x_{-} x_{+} \neq 0$, this implies existence of $g \in \tilde{G}$ 
such that $g \cdot \wt{B}=\wt{B}$ and $g \neq \id_{\tilde{V}}$, 
meaning that $\wt{B}$ has non-trivial stabilizer group in $\wt{G}$.

Depending on a choice of $\theta^{+} \in \wt{\mc C}_{\ell}$ and $\theta^{-} \in \wt{\mc C}_{0}$, 
the line segment
$\lbrace \theta^{t} \mid  0 \le t \le 1 \rbrace$ hits different places in walls between $\wt{C}_{\ell}$ and $\wt{C}_{0}$, and
determines a semi-stable locus $\hat{\mu}^{-1}(0)^{ss}$ by Lemma \ref{lem:ss}. 
\begin{lem}
\label{lem:flatsharp}
There exists $\theta^{+} \in \wt{\mc C}_{\ell}$ and $\theta^{-} \in \wt{\mc C}_{0}$ such that 
for any $\wh{B}=(\wt{B}, [x_{0}, x_{1}]) \in \hat{\mu}^{-1} (0)^{ss}$,  the following holds.
The {\color{black} representations 
	$\wt{B}$ on $\wt{V}$} has a non-trivial stabilizer group if and only if
we have decompositions 
{$\wt{V}=\wt{V}^{\sharp} 
	\oplus \wt{V}^{\flat}$
and 	
	\color{black} $\wt{B} =\wt{B}^{\sharp}
\oplus \wt{B}^{\flat}$ 
of representations $\wt{B}^{\sharp}$
on $\wt{V}^{\sharp}$ and $\wt{B}^{\flat}$ 
on $\wt{V}^{\flat}$}
such that the following conditions hold.
\\
(a) For $\mk I^{\sharp} = \lbrace k \in \mk I \mid \dim V^{\sharp}_{\ast(k)} \neq \dim V^{\sharp}_{\ast(k-1)}  \rbrace$,
we have $\min(\mk I^{\sharp}) \le \ell$.\\
(b) We have $V^{\sharp}_{\infty}=0$ and $\bar{\zeta}(V^{\sharp})=0$.\\
(c) $(\wt{V}^{\sharp}, \wt{B}^{\sharp})$ is $(\bar{\zeta}, 0)$-semistable, and we have no sub-representation 
$S$ of $(\wt{V}^{\sharp}, B^{\sharp})$ such that $\bar{\zeta}(S)=0$, $S_{\ast(\min(\mk I^{\sharp}))} \neq 0$ 
and $S \neq \wt V^{\sharp}$. \\
(d) $V^{\flat}_{\infty}=V_{\infty}$ and $(\wt{V}^{\flat}, \wt{B}^{\flat})$ is $\wt{\mc C}_{\min(\mk I^{\sharp}) -1}$-stable. 
\end{lem}
Objects $(\wt{V}^{\flat}, \wt{B}^{\flat})$ appearing in (d) in Lemma \ref{lem:flatsharp} are parametrized by
$\wt{M}^{\min(\mk I^{\sharp})-1}(\alpha - d^{\sharp} \beta, \mk I \setminus \mk I^{\sharp})$.
In the next subsection, we construct a moduli space parametrizing $(\wt{V}^{\sharp}, \wt{B}^{\sharp})$
satisfying (a), (b), (c) in Lemma \ref{lem:flatsharp}.
The connected components of this moduli space are indexed by an element in the following set
\begin{align}
\label{decompdata}
\mc D_{\beta_{\ast}}^{\ell}(\mk I) = \lbrace \mk I^{\sharp} \subset \mk I \mid |\mk I^{\sharp}| = d^{\sharp} \beta_{\ast} \text{ for } 
d^{\sharp} =1, \ldots, \lfloor \alpha_{\ast} / \beta_{\ast} \rfloor \text{ and } 
\min(\mk I^{\sharp}) \le \ell \rbrace.
\end{align}

\subsection{Moduli stack of destabilizing 
	{\color{black} representations $\wt B^{\sharp}$ 
on $\wt V^{\sharp}$}}
\label{subsec:destab}
From the condition (b) in Lemma \ref{lem:flatsharp}, we need a new framing vertex other than $\infty \in Q_{0}$ to regard  
$(\wt V^{\sharp}, \wt B^{\sharp} )$ as a framed quiver representation.
Comparing the condition (c) in Lemma \ref{lem:flatsharp} with the condition (b) in Lemma \ref{lem:cstab}, we take $\ast(\min(\mk I^{\sharp})) \in \wt{Q}_{0}$ as 
a framing vertex instead of $\infty \in Q_{0}$.
This leads to a new framed quiver $Q^{\sharp}=(Q^{\sharp}, R, \infty^{\sharp})$ as follows.
We define $Q^{\sharp}=(Q^{\sharp}_{0}, Q^{\sharp}_{1})$ by $Q^{\sharp}_{0}=Q_{0} \sqcup \lbrace \infty^{\sharp} \rbrace,
Q^{\sharp}_{1}= Q_{1} \sqcup \lbrace a^{\sharp} \colon \infty^{\sharp} \to \ast \rbrace$ 
adding a new arrow $\infty^{\sharp} \stackrel{a^{\sharp}}{\to} \ast$.
We only consider $Q^{\sharp}_{0}$-graded vector space 
$V^{\sharp} \oplus V^{\sharp}_{\infty^{\sharp}}$ where $V^{\sharp} =\bigoplus_{i \in I}V^{\sharp}_{i}$
such that $(\dim V^{\sharp}_{i})_{i \in I} = d^{\sharp} \beta \in \Z^{I}$ and 
$
V^{\sharp}_{\infty^{\sharp}}=\C.
$ 

For any parameter $\zeta \in \mb R^{I}$, we define 
stability parameter $\zeta^{\sharp}=(\zeta^{\sharp}_{v})_{v \in Q^{\sharp}_{0}} \in \mb R^{Q^{\sharp}_{0}}$ 
by
\begin{align}
\label{zetasharp}
\begin{cases}
\zeta^{\sharp}_{i} = \zeta_{i} \text{ for }i \in I \\
\zeta^{\sharp}_{\infty}=0\\
\zeta^{\sharp}_{\infty^{\sharp}} = - \sum_{i \in I} \zeta_{i} \dim V^{\sharp}_{i}=- \zeta(d^{\sharp} \beta).
\end{cases}
\end{align}
We set $H_{Q}(d^{\sharp} \beta)=M^{\zeta^{\sharp}}_{Q^{\sharp}} (V^{\sharp} \oplus V^{\sharp}_{\infty^{\sharp}})$
for $\zeta \in {\color{black}
\mc C^{-}}$.
\begin{rem}
From our choice \eqref{zetasharp} and $\zeta \in 
{\color{black} \mc C^{-}}$, the moduli space
$H_{Q}(d^{\sharp} \beta)$ parametrizes 
$Q^{\sharp}$-representations satisfying conditions (a) and (b) in Lemma \ref{lem:cstab}
replacing $\infty$ with $\infty^{\sharp}$.
When $Q$ is {\color{black} the double of 
the Jordan quiver together with the 
preprojective relation}
giving constructions of the Nakajima quiver variety of $A^{(1)}_{0}$-type,
$H_{Q}(d^{\sharp} \beta)$ is isomorphic to the Hilbert scheme of points on $\C^{2}$.
\end{rem}
We consider a tautological homomorphism $\mc V^{\sharp}_{\infty^{\sharp}} \to \mc V^{\sharp}_{\ast}$ 
on $H_{Q}(d^{\sharp} \beta)$ corresponding to the arrow from $\infty^{\sharp}$ to $\ast$.
It is a non-zero homomorphism from the line bundle $\mc V^{\sharp}_{\infty^{\sharp}}$ by the above remark.
In fact, otherwise we can take a sub-representation $V^{\sharp}_{\infty^{\sharp}} \subset V^{\sharp} \oplus V^{\sharp}_{\infty^{\sharp}}$ 
violating the condition (b) in Lemma \ref{lem:cstab}. 
We set $\bar{\mc V}^{\sharp}_{\ast} = \mc V^{\sharp}_{\ast} / \mc V^{\sharp}_{\infty^{\sharp}}$,
and consider the full flag bundle $Fl_{H_{Q}(d^{\sharp} \beta)}
(\bar{\mc V}^{\sharp}_{\ast}, \bar{\mk I}^{\sharp} )$ over $H_{Q}(d^{\sharp} \beta)$ where
$\bar{\mk I}^{\sharp} =\mk I^{\sharp} \setminus \lbrace \min(\mk I^{\sharp}) \rbrace$.
An element in $Fl_{H_{Q}(d^{\sharp} \beta)}(\bar{\mc V}^{\sharp}_{\ast}, \bar{\mk I}^{\sharp} )$ 
is represented by a collection $(B^{\sharp}, a^{\sharp}, \bar{F}^{\sharp}_{\bullet})$ 
of $\zeta^{\sharp}$-stable $Q$-representation $B^{\sharp}$ on $V^{\sharp}$, 
a injection $a^{\sharp} \colon V^{\sharp}_{\infty} =\C \to V^{\sharp}_{\ast}$, 
and a full flag $\bar{F}^{\sharp}_{\bullet}$ of $\bar{V}^{\sharp}=V^{\sharp}_{\ast} / a^{\sharp} (V^{\sharp}_{\infty})$.
It is identified with an object $(\wt V^{\sharp}, \wt B^{\sharp})$ which appeared in
Lemma \ref{lem:flatsharp} as follows.

We recall that $\wt{Q}$-representations $(\wt V^{\sharp}, \wt B^{\sharp})$ are regarded as
enhanced $Q$-representations $(B^{\sharp}, F^{\sharp}_{\bullet})$ as in \S \ref{subsec:enhancement}.
On the other hand, from the above collection $(B^{\sharp}, a^{\sharp}, \bar{F}^{\sharp}_{\bullet})$,
we set $F^{\sharp}_{\bullet}$ as follows.
For $i \ge \min(\mk I^{\sharp})$, we write 
by $F^{\sharp}_{i}$ the inverse images of $\bar{F}^{\sharp}_{i}$ by 
the quotient map $V^{\sharp}_{\ast} \to \bar{V}^{\sharp}_{\ast}$.
For $i < \min(\mk I^{\sharp})$, we set $F^{\sharp}_{i}=0$.
This also gives enhanced $Q$-representations $(B^{\sharp}, F^{\sharp} )$, and
we have $F^{\sharp}_{\bullet} \in Fl(V^{\sharp}_{\ast}, \mk I^{\sharp})$ since 
$\bar{F}^{\sharp}_{\bullet} \in Fl(\bar{V}^{\sharp}_{\ast}, \bar{\mk I}^{\sharp})$.

\begin{lem}
\label{lem:sharp}
For a collection $(B^{\sharp}, a^{\sharp}, \bar{F}^{\sharp}_{\bullet})$ of a $Q$-representation $B^{\sharp}$ on $V^{\sharp}$, 
an injection $a^{\sharp} \colon V^{\sharp}_{\infty^{\sharp}} \to V^{\sharp}_{\ast}$, 
and a full flag $\bar{F}^{\sharp}_{\bullet} \in Fl(\bar{V}^{\sharp}_{\ast}, \bar{\mk I}^{\sharp})$, 
the corresponding data 
$(B^{\sharp}, F^{\sharp}_{\bullet})$ satisfies condition (c) in Lemma \ref{lem:flatsharp}
if and only if $(B^{\sharp}, a^{\sharp})$
represents an element in $Fl_{H_{Q}(d^{\sharp} \beta)}( \bar{\mc V}^{\sharp}_{\ast}, \bar{\mk I}^{\sharp} )$ 
\end{lem}
\proof 
Since $F^{\sharp}_{\min(\mk I^{\sharp})} = \im a^{\sharp}$, it follows from comparison between 
the condition (c) in Lemma \ref{lem:flatsharp} and
the condition (b) in Lemma \ref{lem:cstab}.
\endproof

\begin{cor}
	\label{cor:setbij}
There exists a set theoretical bijection between $\mc M_{exc}$ and 
\[
\bigsqcup_{\mk I^{\sharp} \in \mc D_{\beta_{\ast}}^{\ell}(\mk I)} 
Fl_{H_{Q}(d^{\sharp} \beta)}(\bar{\mc V}^{\sharp}_{\ast}, \bar{\mk I}^{\sharp})
\times
\wt M^{\min(\mk I^{\sharp})-1}(\alpha - d^{\sharp} \beta, \mk I^{\flat}), 
\]
where we set $d^{\sharp} = \lvert \mk I^{\sharp} \rvert / \beta_{\ast}$ and $\mk I^{\flat}=\mk I \setminus \mk I^{\sharp}$
for $\mk I^{\sharp}$ in the index set $\mc D_{\beta_{\ast}}^{\ell}(\mk I)$ defined in \eqref{decompdata}.
\end{cor}
\proof
For each $\mk I^{\sharp} \in \mc D_{\beta_{\ast}}^{\ell}(\mk I)$, we take 
data $((B^{\sharp}, a^{\sharp}, \bar{F}^{\sharp}_{\bullet}), (B^{\flat}, F^{\flat}))$ representing elements in 
\[
Fl_{H_{Q}(d^{\sharp} \beta)}(\bar{\mc V}^{\sharp}_{\ast}, \bar{\mk I}^{\sharp})
\times
\wt M^{\min(\mk I^{\sharp})-1}(\alpha - d^{\sharp} \beta, \mk I^{\flat}).
\]
From these data we get representatives 
$((B^{\sharp}, F^{\sharp}_{\bullet}), (B^{\flat}, F^{\flat} ), [1,1] )$ of elements in $\mc M_{exc}$. 
By Lemma \ref{lem:flatsharp} and Lemma \ref{lem:sharp}, this gives a well-defined bijection.
\endproof
This corollary suggest that $\mc M_{exc}$ have connected components $\mc M_{\mk I^{\sharp}}$ indexed by
$\mk I^{\sharp} \in \mc D_{\beta_{\ast}}^{\ell}(\mk I)$ such that points of $\mc M_{\mk I^{\sharp}}$ are 
parametrized by $Fl_{H_{Q}(d^{\sharp} \beta)}(\bar{\mc V}^{\sharp}_{\ast}, \bar{\mk I}^{\sharp}) \times
\wt M^{\min(\mk I^{\sharp})-1}(\alpha - d^{\sharp} \beta, \mk I^{\flat})$.
But we have a non-trivial stabilizer group of $((B^{\sharp}, F^{\sharp}_{\bullet}), (B^{\flat}, F^{\flat} ), [1,1] )$
consisting of $\zeta \id_{V^{\sharp}}$ for the primitive $m^{\sharp}$-th root
$\zeta=e^{2 \pi \sqrt{-1}/m^{\sharp}}$ of unity for
$m^{\sharp} =d^{\sharp}\theta^{+} (\beta) - d^{\sharp} \theta^{-}(\beta)$.
Hence we need to take \'etale cover $\mc M'_{\mk I^{\sharp}}$ of $\mc M_{\mk I^{\sharp}}$
of degree $\frac{1}{m^{\sharp}}$.   
Precise statements are given in Theorem \ref{thm:decomp}
in the next subsection.

We consider a modified action 
$\C^{\ast}_{e^{\hbar/m^{\sharp} }} \times \mc M \to \mc M$ induced by 
\begin{align}
\label{act}
\left( \rho, F_{\bullet}, [x_{-},x_{+}] \right) \mapsto 
\left( e^{\hbar/ m^{\sharp}} \id_{\tilde{V}^{\sharp} \oplus } \id_{\tilde{V}^{\flat}} \right) 
\left( \rho, F_{\bullet}, [e^{\hbar}x_{-},x_{+}] \right). 
\end{align}
This action is equal to the original $\C^{\ast}_{\hbar}$-action \eqref{fiberwise}, since the difference is absorbed in $G$-action.

\subsection{Decomposition of $\mc M^{\C^{\ast}_{\hbar}}$}
For tautological bundles $\mc V^{\sharp}, \mc V^{\sharp}_{\infty^{\sharp}}$ over $H_{Q}(d^{\sharp} \beta)$
and $\mc V^{\flat}$ over $\wt M^{\min(\mk I^{\sharp})-1}(\alpha - d^{\sharp} \beta, \mk I^{\flat})$, we also
write by the same letter their pull-backs to 
$Fl_{H_{Q}(d^{\sharp} \beta)}(\bar{\mc V}^{\sharp}_{\ast}, \bar{\mk I}^{\sharp}) \times 
\wt M^{\min(\mk I^{\sharp})-1}(\alpha - d^{\sharp} \beta, \mk I^{\flat})$
by the projections.
We take universal full flags $\bar{\mc F}_{\bullet}^{\sharp}$ on 
$Fl_{H_{Q}(d^{\sharp} \beta)}( \bar{\mc V}^{\sharp}_{\ast}, \bar{\mk I}^{\sharp})$.
For $i \ge \min(\mk I^{\sharp})$, we write 
by $\mc F^{\sharp}_{i}$ the inverse images of $\bar{\mc F}^{\sharp}_{i}$ by 
the tautological homomorphism $\mc V^{\sharp}_{\infty^{\sharp}} \to \mc V^{\sharp}_{\ast}$, and
set $\mc F_{i}=0$ for $i < \min(\mk I^{\sharp})$.

\begin{thm}
\label{thm:decomp}
We have a decomposition into connected components
\begin{align}
\label{decomp1}
\mc M^{\C^{\ast}_{\hbar}} = \mc M_{+} \sqcup \mc M_{-} \sqcup 
\bigsqcup_{\mk I^{\sharp} \in \mc \mc D_{\beta_{\ast}}^{\ell}(\mk I) }\mc M_{\mk I^{\sharp}}
\end{align}
such that the followings hold.\\
(i) We have $\mc M_{+} \cong \wt M^{\ell}(\alpha, \mk I)$ and $\mc M_{-}$ is isomorphic to the full flag bundle 
$Fl_{M^{\zeta^{-}}(\alpha)}(\mc V_{\ast}, \mk I)$ of the tautological bundle $\mc V_{\ast}$ on $M^{\zeta^{-}}(\alpha)$
with $\lbrace k \in [L] \mid \mc F_{k} / \mc F_{k-1} \rbrace = \mk I$.
\\
(ii) For each $\mk I^{\sharp} \in \mc \mc D_{\beta_{\ast}}^{\ell}(\mk I)$, we have $\C^{\ast}_{\hbar/m}$-equivariant
finite \'etale morphisms 
$
\Phi \colon \mc M'_{\mk I^{\sharp}} \to \mc M_{\mk I^{\sharp}}$ of degree $1/ m^{\sharp}$, and 
$\Psi \colon \mc M_{\mk I^{\sharp}}' \to  
Fl_{H_{Q}(d^{\sharp} \beta)}( \bar{\mc V}^{\sharp}_{\ast}, \bar{\mk I}^{\sharp}) \times 
\wt M^{\min(\mk I^{\sharp})-1}(\alpha  - d^{\sharp} \beta, \mk I^{\flat})
$ 
of degree $1/ (m^{\sharp})^{2}$, where  
$\mk I^{\flat}=\mk I \setminus \mk I^{\sharp}$, $m^{\sharp}=d^{\sharp} \theta^{+} (\beta) - d^{\sharp} \theta^{-}(\beta)$, and
$d^{\sharp}=\lvert \mk I^{\sharp} \rvert / \beta_{\ast}$. \\
(iii) As $\C^{\ast}_{\hbar/m^{\sharp}}$-equivariant 
$Q_{0}$-graded vector bundles on $\mc M_{\mk I^{\sharp}}'$, 
we have a decomposition
\begin{align*}
\Phi^{\ast} \mc V|_{\mc M_{\mk I^{\sharp}}} 
& \cong
\left( \Psi^{\ast} \mc V^{\sharp} \otimes L 
\otimes \C_{e^{\hbar/m^{\sharp}}} \right)
\oplus
\Psi^{\ast} \mc V^{\flat}
\\
\Phi^{\ast} \mc F_{\bullet}|_{\mc M_{\mk I^{\sharp}}} 
& \cong 
\left( \Psi^{\ast} \mc F^{\sharp}_{\bullet} \otimes L \otimes \C_{e^{\hbar/m^{\sharp}}} \right) \oplus \Psi^{\ast} \mc F^{\flat}_{\bullet},
\end{align*} 
where $L$ is a line bundle on $\mc M_{\mk I^{\sharp}}'$. 
\end{thm}
\proof
See \cite[Theorem 4.3]{O5}.
\endproof
{\color{black}
The component $\mc M_{\mk I^{\sharp}}$
is constructed as a quotient stack
$U_{\mk I^{\sharp}}/
G_{\mk I^{\sharp}}$, where 
$U_{\mk I^{\sharp}}/
G_{\mk I^{\sharp}}$ is  
a Zariski open subset of an Affine space, and $G_{\mk I^{\sharp}}$ is 
a product of general linear groups
acting on $U_{\mk I^{\sharp}}$
with the finite stabilizer groups.
Then the \'etale cover $\mc M_{\mk I^{\sharp}}'$ in Theorem \ref{thm:decomp} is obtained by 
replacing a torus component 
$\C^{\ast}_{t}$ in $G_{\mk I^{\sharp}}$ with 
$\C^{\ast}_{t^{1/m^{\sharp}}}$
(see proof of \cite[Theorem 4.3 (iii)]{O5}).
This gives non-integer degree 
$1/m^{\sharp}$ by 
\cite[Definition 1.15]{V}.
}

We write by $N_{+}, N_{-}$, and $N_{\mk I^{\sharp}}$ normal bundles of $\mc M_{+}, \mc M_{-}$ and $\mc M_{\mk I^{\sharp}}$ in $\mc M$ respectively.
\begin{lem}[\protect{\cite[Proposition 5.9.2]{M}}]
\label{ob1}
We have 
$N_{\pm} = e^{\pm \hbar} \otimes \bigotimes_{i \in \wt{Q}_{0}\setminus \lbrace \infty \rbrace}
\det \mc V_{i}^{\otimes \mp (\theta^{+}_{i} - \theta^{-}_{i})}$.
\end{lem}

We take universal full flags $\bar{\mc F}_{\bullet}^{\sharp}$ on 
$Fl_{H_{Q}(d^{\sharp} \beta)}( \bar{\mc V}^{\sharp}_{\ast}, \bar{\mk I}^{\sharp})$, 
and $\mc F_{\bullet}^{\flat}$ on $M_{\wt{Q}}^{\bar{\zeta}, \min(\mk I_{\sharp}) -1} (\wt{V}_{\flat})$.
On $Fl_{H_{Q}(d^{\sharp} \beta)}( \bar{\mc V}^{\sharp}_{\ast}, \bar{\mk I}^{\sharp})$,
we define a full flag $\mc F_{\bullet}^{\sharp}$ of $\mc V^{\sharp}_{\ast}$ 
by the pull-back of $\bar{\mc F}_{\bullet}^{\sharp}$ to $\mc V^{\sharp}_{\ast}$ for $i \neq \min(\mk I^{\sharp})$
and $\mc F_{\min(\mk I^{\sharp})}^{\sharp} = \im (\mc V^{\sharp}_{\infty^{\sharp}} \to \mc V^{\sharp}_{\ast})$. 
We also write their pull-backs to the product $\mc M_{\sharp} \times \mc M_{\flat}$
by the same letter $\mc F^{\sharp}_{\bullet}, \mc F^{\flat}_{\bullet}$.

For vector bundles $\mc E, \mc F$ on a stack $\mc Z$, we write by 
$\mc Hom(\mc E, \mc F)$ the vector bundle $\mc E^{\vee} \otimes \mc F$ on $\mc Z$.
For two flags $\mc F_{\bullet}, \mc F_{\bullet}'$ of sheaves on the same Deligne-Mumford stack,  
we set 
\begin{align*}
\Theta (\mc F_{\bullet}, \mc F_{\bullet}') 
&=
\sum_{l < l' }  
\mc Hom\left(\mc F_{l} / \mc F_{l-1}, \mc F_{l'}' /\mc F_{l'-1}' \right),
\\
\mk H (\mc F_{\bullet}, \mc F_{\bullet}')
&=
\Theta (\mc F_{\bullet}, \mc F_{\bullet}') 
+
\Theta (\mc F_{\bullet}', \mc F_{\bullet}).
\end{align*}
When $\mc F_{\bullet} = \mc F_{\bullet}'$, we set $\Theta(\mc F_{\bullet}) = \Theta(\mc F_{\bullet}, \mc F_{\bullet})$.

\begin{lem} 
\label{nmki}
We have $\Phi^{\ast} N_{\mk I^{\sharp} } \cong 
\Psi^{\ast} \left[ \mathfrak{N}( \mc V^{\sharp} \otimes
L \otimes e^{\hbar/m^{\sharp}}, \mc V^{\flat})  + 
\mk H (\mc F^{\sharp}_{\bullet} \otimes
L \otimes e^{\hbar/m^{\sharp}}, \mc F^{\flat}_{\bullet}) \right]
$ where
\begin{align}
\mathfrak{N}( \mc V^{\sharp}, \mc V^{\flat})
&=
\sum_{(\heartsuit, \spadesuit)} 
\left(
\sum_{\substack{h \in Q_{1} }} 
\mc H om( \mc V^{\heartsuit}_{\out(h)}, 
\mc V^{\spadesuit}_{\inn(h)})
- 
\sum_{\gamma \in Q_{2}} \mc H om( 
\mc V^{\heartsuit}_{\out(\gamma)}, 
\mc V^{\spadesuit}_{\inn(\gamma)}) 
-
\sum_{i \in I} 
\mc Hom ( \mc V^{\heartsuit}_{i}, 
\mc V^{\spadesuit}_{i}) \right),
\label{normal}
\end{align}
where the sum is taken for only cross terms 
$(\heartsuit, \spadesuit)= (\sharp, \flat), 
(\flat, \sharp)$.
\end{lem}
Here $m^{\sharp}$ is the parameter appeared in 
Theorem \ref{thm:decomp} (ii).


\subsection{Relative tangent bundles for flags}
 \label{subsec:rela2}
We consider the pull-back $\Theta^{rel}$ to $\mc M$ of 
the relative tangent bundle of $[\wt{\mu}^{-1}(0)/G]$ 
over $[\mu^{-1}(0)/G]$. 
After restricting to $\mc M_{\mk I^{\sharp}}$ and 
pulling back to $\mc M_{\mk I^{\sharp}}'$, we have a 
decomposition 
\begin{align}
\label{rel}
\Phi^{\ast} \Theta^{rel} = \Psi^{\ast} 
( \Theta(\mc F_{\bullet}^{\sharp}) \oplus 
\Theta(\mc F_{\bullet}^{\flat}) \oplus 
\mk H (\mc F_{\bullet}^{\sharp} \otimes
L \otimes e^{\hbar/m^{\sharp}}, 
\mc F_{\bullet}^{\flat})  )
\end{align} 
We also consider the relative tangent bundle $\Theta'$ of 
the projection 
$Fl_{H_{Q} (d^{\sharp} \beta)}( \bar{\mc V}^{\sharp}_{\ast}, \bar{\mk I}^{\sharp}) \to
H_{Q} (d^{\sharp} \beta)$.
Then we have an exact sequence
\begin{align}
\label{ob4}
0 \to \Theta' \to \Theta(\mc F_{\bullet}^{\sharp}) \to \bar{\mc V}^{\sharp}_{\ast} \to 0,
\end{align}
where $\bar{\mc V}^{\sharp}_{\ast}$ is the quotient 
by the tautological homomorphism $\mc V^{\sharp}_{\infty^{\sharp}} \to \mc V^{\sharp}_{\ast}$, 
and it is identified with $\mc Hom (\mc F^{\sharp}_{\min(\mk I^{\sharp})}, 
\mc V^{\sharp}_{\ast}/ \mc F^{\sharp}_{\min(\mk I^{\sharp})})$.

\begin{NB}
For any full flag bundle $p \colon Y =Fl(\E, \underbar{n}) \to X$ 
for a rank $n$ vector bundle $\E$ on $X$, 
we can check $p_{\ast} \left( e(T_{Y/X}) \cap p^{\ast} \beta \right) = 
n! \beta$ for $\beta \in a_{\bullet}(X)$ as follows.
Since full flag bundles are obtained by iterations of projective bundles, 
it is enough to show that for $\pi \colon Z =\PP(\E, \underbar{n}) \to X$, 
we have $\pi_{\ast} \left( e(T_{Z/X}) \cap \pi^{\ast} \beta \right)  = n \beta$ for $\beta \in a_{\bullet}(X)$.
This follows from the description of $\pi_{\ast} \colon a_{\bullet}(Z) \to a_{\bullet}(X)$ in \cite{F}.
\end{NB}

\section{Wall-crossing formula}
\label{sec:wall}

We set $M^{0}_{Q}(\alpha)=M^{0}(\alpha)=\Spec 
\Gamma(\mb M(V), \mo_{\mb M(V)})^{G}$, and write by
$\Pi = \Pi_{Q} \colon M^{\zeta}(\alpha) \to M^{0}(\alpha)$ 
the natural projection.

\subsection{Equivariant framed quiver}
\label{subsec:eq}

We introduce an algebraic torus $\mb T=(\C^{\ast})^{n}$
and use torus $\mb T$-equivariant Chow rings to define 
integrals in general.
For a character $\chi \colon \mb T \to \C^{\ast}$, 
we write by $\C_{\chi}$ the weight space, 
that is, the one dimensional $\mb T$-representation 
with the eigenvalue $\chi$.
For the coordinates $(s_{1}, \ldots, s_{n})$ of 
$\mb T$, we regard a monomial $s_{1}^{m_{1}} \cdots 
s_{n}^{m_{n}}$ a character $\mb T \to \C^{\ast}$ sending
$p \in \mb T$ to $s_{1}(p)^{m_{1}} \cdots 
s_{n}(p)^{m_{n}}$.
The monomial $s_{1}^{m_{1}} \cdots s_{n}^{m_{n}}$ is also 
used as an abbreviation of the weight space
$\C_{s_{1}^{m_{1}} \cdots s_{n}^{m_{n}}}$.
In particular, the $\mb T$-equivariant Chow ring
$A_{\mb T}=A_{\mb T}^{\bullet}(\text{pt})$ over the point pt,
is isomorphic to the polynomial ring
$\Z[\xi_{1}, \ldots, \xi_{n}]$, where 
$\xi_{1}=c_{1}(s_{1}), \ldots, \xi_{n}=c_{1}(s_{n})$ 
are called the {\it equivariant parameters}.

We take a pair $(Q, \chi_{Q})$ where 
\[
\chi_{Q}=((\chi_{h})_{h \in Q_{1}}, 
(\chi_{\gamma})_{\gamma \in Q_{2}})
\]
is a collection of 
characters of $\mb T$.
We replace \eqref{mbM} with
\begin{align}
\label{tmbM}
\mb M_{Q}(V)=\prod_{h \in Q_{1}} \Hom_{\C} 
(V_{\out (h)}, V_{\inn (h)}) \otimes \C_{\chi_{h}}, 
\quad 
\mb L_{Q}(V)=\prod_{\gamma \in Q_{2}} \Hom_{\C} 
(V_{\out (\gamma)}, V_{\inn (\gamma)})
\otimes \C_{\chi_{\gamma}}.
\end{align}

\begin{defn}
\label{eqfrq}
A pair $(Q, \chi_{Q})$ is called a 
{\it $\mb T$-equivariant framed quiver} 
when the {\color{black} map}
$\mu=\mu_{Q} \colon \mb M_{Q}(V) 
\to \mb L_{Q}(V)$ is $\mb T$-equivariant for the
$\mb T$-action induced from \eqref{tmbM}.
\end{defn} 
 
For a $\mb T$-equivariant framed quiver 
$Q=(Q, \chi_{Q})$, we have natural $\mb T$-action on 
moduli spaces $M^{\zeta}(\alpha)$ of $\zeta$-semistable
$Q$-representations for any $\zeta \in \mb R^{I}$.
Furthermore $\Pi_{Q} \colon M^{\zeta}(\alpha) \to 
M^{0}(\alpha)$ are $\mb T$-equivariant.
In the following, we assume that
the fixed points sets $M^{0}_{Q}(\alpha)^{\mb T}$ are 
isomorphic to the point $\pt$ for 
any framed quiver $Q$.
In particular, $Q^{\sharp}=(Q^{\sharp}, 
\chi_{Q^{\sharp}})$ are regarded as 
$\mb T$-equivariant framed quiver with 
$\chi_{Q^{\sharp}}$ naturally induced from 
$\chi_{Q}$.

\subsection{Equivariant Rieman-Roch}
\label{subsec:k-th}

We recall Riemann-Roch theorem for equivariant setting 
from \cite{EG}.
We consider a Deligne-Mumford stack $X$ with a 
$\mb T$-action. 
Let $CH^{\bullet}_{\mb T}(X)$ denote the 
$\mb T$-equivariant 
Chow group of $X$,
$A^{\bullet}_{\mb T}(X)$ the $\mb T$-equivariant Chow ring 
of $X$,
and $K^{\mb T}(X)$ the Grothendieck group of 
$\mb T$-equivariant coherent sheaves on $X$.
When $X$ is smooth, we identify 
$CH^{\bullet}_{\mb T}(X)$ with $A^{\bullet}_{\mb T}(X)$, 
and $K^{\mb T}(X)$ with the Grothendieck ring 
$K_{\mb T}(X)$ of 
$\mb T$-equivariant 
vector bundles on $X$.

When $X$ is a point $\pt$, then $K^{\mb T}(\pt)$ is 
equal to the 
$\mb T$-representation ring 
$R_{\mb T}=\Z[s_{1}^{\pm 1}, \ldots, s_{n}^{\pm 1}]$.
Here $(s_{1}, \ldots, s_{n})$ is the coordinate of 
$\mb T=(\C^{\ast})^{n}$, and 
for $l_{1}, \ldots, l_{n} \in \Z$ the monomial 
$s_{1}^{l_{1}} \cdots s_{n}^{l_{n}}$ 
denotes the weight space $\C_{s_{1}^{l_{1}} \cdots s_{n}^{l_{n}} }$ 
with the eigen-value
$s_{1}^{l_{1}} \cdots s_{n}^{l_{n}}$.
We can also see that $CH^{\bullet}_{\mb T}(\pt)
= A^{\bullet}_{\mb T}(\pt)$ is equal to
the polynomial ring $\Z[\xi_{1}, \ldots, \xi_{n}]$ where 
$\xi_{1} = c_{1}(s_{1}), \ldots, \xi_{n}=c_{1}(s_{n})$.

For a subset 
\[
S_{\mb T}=\lbrace 1- s_{1}^{l_{1}} \cdots s_{n}^{l_{n}} 
\mid (l_{1}, \ldots, l_{n}) \in \Z^{n+1} \setminus 
\lbrace 0 \rbrace \rbrace,
\]
we consider the ideal $I_{\mb T}=\langle S_{\mb T} 
\rangle$ of $R_{\mb T}$ generated by $S_{\mb T}$, and
write by $\wh{R}_{\mb T}$ the completion of $R_{\mb T}$ 
along $I_{\mb T}$.
We set $\wh{K}^{\mb T}(X)=
K^{\mb T}(X)\otimes_{R_{\mb T}} \wh{R}_{\mb T}$, and
$\wh{CH}^{\bullet}_{\mb T}(X)= 
\prod_{i \ge 0} CH^{i}_{\mb T}(X)$.
In particular, we set
\[
A_{\mb T}=\Z[e^{\pm \xi_{1}}, \ldots, e^{\pm \xi_{n}}]
\subset 
\wh{CH}_{\mb T}^{\bullet} (\pt) 
= \Z[[ \xi_{1}, \ldots, \xi_{n}]].
\]
We set $Td=Td_{X}=\frac{\text{Eu}( T X)}{ \euk(T X)}$, 
and define a map $\tau=\tau_{X} \colon K^{\mb T}(X) \to 
\wh{CH}^{\bullet}_{\mb T} (X)$ by 
\[
\tau(\E) = 
\left( 
\ch(\E) \cdot \td_{X} 
\right) \cap [X]
\] 
for $\E \in 
K^{\mb T}(X)$.
The map $\tau$ is compatible 
with push-forwards, and 
factors through the induced isomorphism
\begin{align}
\label{isom}
\wh{K}^{\mb T}(X) \cong 
\wh{CH}^{\bullet}_{\mb T}(X).
\end{align}
We often use identification 
$R_{\mb T} \cong A_{\mb T}$
via the exponential map $\tau_{\pt}$.

By the localization theorem \cite[Theorem 2.1]{T}, 
we have an isomorphism 
$S_{\mb T}^{-1} K^{\mb T}(X) \cong 
S_{\mb T}^{-1} K^{\mb T}(X^{\mb T})$.
Using \eqref{isom}, we also have
$T_{\mb T}^{-1} \wh{CH}^{\bullet}_{\mb T}(X) \cong 
T_{\mb T}^{-1} \wh{CH}^{\bullet}_{\mb T}(X^{\mb T})$, 
where 
\[
T_{\mb T}=\tau(S_{\mb T})=\lbrace 1 - e^{l_{1} \xi_{1}
+ \cdots + l_{n} \xi_{n} } 
\mid (l_{1}, \ldots, l_{n}) \in \Z^{n} \setminus 
\lbrace 0 \rbrace \rbrace.
\]
In our setting, we have the diagram:
\[
\xymatrix{
K^{\mb T}(M^{\zeta}(\alpha)) 
\ar[r]^-{\tau_{M^{\zeta}(\alpha)}} 
\ar[d]_{\Pi_{\ast}}
&
\wh{CH}^{\bullet}_{\mb T} (M^{\zeta}(\alpha))
\ar[d]^{\Pi_{\ast}}
\\
K^{\mb T}(M^{0}(\alpha)) \ar[r]^-{\tau_{M^{0}(\alpha)}} 
\ar[d]_{(\iota_{\ast})^{-1}}
&
\wh{CH}^{\bullet}_{\mb T} (M^{0}(\alpha))
\ar[d]^{(\iota_{\ast})^{-1}}
\\
S_{\mb T}^{-1} K^{\mb T}( \pt ) \ar[r]^-{\tau_{\pt}} 
&
T_{\mb T}^{-1} \wh{CH}^{\bullet}_{\mb T} ( \pt )
}
\]
where $\iota \colon \pt = M^{0}(\alpha)^{\mb T} \to 
M^{0}(\alpha)$ is the embedding.


\subsection{Virtual fundamental cycle}
{\color{black}
We recall constructions of virtual fundamental cycles \cite{BF} in the following setting.
We introduce a triple $\mc K=(Y, P, \sigma)$ 
{\color{black}called {\it Kuranishi structure} originated from 
\cite{FO}.
} 
Here $Y= (Y, \mo_{Y}(1))$ is a smooth scheme with an ample line bundle $\mo_{Y}(1)$, $P$ is a group scheme acting on $(Y, \mo_{Y}(1))$, and $\sigma \colon Y \to E$ is a $P$-invariant section of $P$-equivariant vector bundle $E$ over $Y$.
We take a semistable locus $Y^{ss}$ of $Y$ with respect to $\mo_{Y}(1)$, and put $Z=\sigma^{-1}(0) \cap Y^{ss}$.
We assume that any point in $Y^{ss}$ has a finite stabilizer group.
Hence $\mc Y=[Y^{ss}/P]$ is a Deligne-Mumford stack.

We define a moduli stack $\mc Z=\mc Z_{\mc K}=[Z/P]$ as the quotient 
stack, 
{\color{black}
and say that $\mc K=(Y, P, \sigma)$  
gives a {\it Kuranishi structure} of $\mc Z$.}
In our case, $E$ is always a trivial bundle $Y \times \C^{l}$, and we have $\sigma=\id_{Y} \times \varphi \in \Gamma(Y, E)$
where $\varphi \colon Y \to \C^{l}$ is a $P$-invariant function. 
In such a case, we just write $\mc K=(Y, P, \varphi)$ instead of $\sigma = \id_{Y} \times \varphi$.

To define {\it virtual fundamental cycle} $[\mc Z]^{vir}$, we consider 
a vector bundle $\mc E=[(E|_{Y^{ss}}) / P]$ over $\mc Y$.
We consider the zero section $0 \colon \mc Y \to \mc E$ 
and write by the same letter $\sigma \colon \mc Y \to \mc E$
the section induced by $\sigma \in \Gamma(Y, E)$. 
We also consider a fundamental cycle $[\mc Y]$ of smooth Deligne-Mumford stack $\mc Y$.
For example, we can take a Kuranishi structure $\mc K_{\mc M}=(\wh{\mb M}, \wt{G}, \hat{\mu})$ 
of the enhanced master space $\mc M$ 
as we saw in Definition \ref{defn:ems}.

Using these we define the virtual fundamental cycle of $\mc Z$ by $[\mc Z]^{vir}=0^! [\mc Y]$, 
where $0^!$ is the refined Gysin homomorphism by the zero section $0 \colon \mc Y \to \mc E$
in the following diagram :
\[
\xymatrix{
	\mc Y \ar[r]^{\sigma} & \E \\
	\ar[u] \mc Z \ar[r] & \mc Y \ar[u]^{0}},
\]

We take another Kuranishi structure $\mc K_{1}=(Y_{1}, P_{1}, \sigma_{1})$ of another stack 
$\mc Z_{1}=\mc Z_{\mc K_{1}}$ where $\sigma_{1} \in \Gamma(Y_{1}, E_{1})$ and $E_{1}$
is a vector bundle on $Y_{1}$. 
We assume that there exists a bundle morphism 
$E_{1} \to E$ equivariant via a group homomorphism 
$P_{1} \to P$ such that sections $\sigma_{1}$ and $\sigma$ coincide 
via the isomorphism $E_{1} \to E|_{Y_{1}}$.
This collection of morphisms are called a morphism of Kuranishi structures.
Furthermore when these morphisms induces an \'{e}tale morphism $[Y_{1}^{ss}/P_{1}] \to [Y^{ss}/P]$, 
we call such a collection of morphisms a {\it coordinate change} of Kuranishi structures. 

When we have a morphism $\mc K_{1}=(Y_{1}, P_{1}, \sigma_{1}) \to 
\mc K=(Y, P, \sigma)$ of Kuranishi structures, such that $Y_{1} \to Y$ is flat, then we have 
$[\mc Z_{1}] = \Psi^{\ast} [\mc Z]$.
In general, virtual fundamental cycles can be different depending on choice of Kuranishi chatrts. 


{\color{black}We remark that we can always take Kuranishi structure consisting of 
single {\it Kuranishi neighborhood} globally covering a moduli stack $\mc Z$ in our setting.
In general, we need an open covering 
of $\mc Z$ consisting of 
Kuranishi neighborhoods as in 
\cite{FO}, \cite{FOOO}.
In algebro-geometric setting, we have obstruction theory developed in \cite{BF} to construct virtual fundamental cycle in general case.
}}

\subsection{$K$-theoretic integral}
\label{subsec:k-th}
Given a collection of $\mb T$-equivariant framed quivers
$(Q, \chi_{Q}),
(Q^{(1)}, \chi_{Q^{(1)}}), \ldots, 
(Q^{(j)}, \chi_{Q^{(j)}})$, 
we take generic stability parameter 
$\zeta \in \mb R^{Q \setminus 
\lbrace \infty \rbrace}$, and 
$\zeta^{(k)} \in \mb R^{Q^{(k)} \setminus 
\lbrace \infty \rbrace}$ for $k=1,\ldots, j$.
We consider moduli spaces 
$M^{\zeta}_{Q}(\alpha)$
for $\alpha \in \Z^{Q^{\infty}_{0}}$, and 
$M^{\zeta^{(k)}}_{Q^{(k)}}(\alpha^{(k)})$  
for $\alpha^{(k)} \in \Z^{Q^{(k)}_{0}}$
for $k=1,\ldots,j$. 
We consider the products
\[
\Pi=
\id \times \prod_{k=1}^{j} \Pi_{Q^{(k)}} \colon
M^{\zeta}_{Q}(\alpha)
\times
\prod_{k=1}^{j}
M^{\zeta^{(k)}}_{Q^{(k)}}(\alpha^{(k)})
\to
M^{\zeta}_{Q}(\alpha)
\times \prod_{k=1}^{j}
M^{0}_{Q^{(k)}}(\alpha^{(k)})
\]
The products $\prod_{k=1}^{j}
M^{\zeta^{(k)}}_{Q^{(k)}}(\alpha^{(k)})$ 
and $M^{\zeta}_{Q}(\alpha)
\times
\prod_{k=1}^{j}
M^{\zeta^{(k)}}_{Q^{(k)}}(\alpha^{(k)})$
of moduli 
spaces  carry 
{\it virtual fundamental cycles} 
$\left[ \prod_{k=1}^{j}
M^{\zeta^{(k)}}_{Q^{(k)}}(\alpha^{(k)}) \right]^{vir}$ 
and 
$\left[ M^{\zeta}_{Q}(\alpha)
\times
\prod_{k=1}^{j}
M^{\zeta^{(k)}}_{Q^{(k)}}(\alpha^{(k)})
\right]^{vir}$ defined as in 
\cite[5.1]{O5}.

For a cohomology class $C \in 
A^{\bullet}_{\mb T} \left( M^{\zeta}_{Q}(\alpha)
\times \prod_{k=1}^{j}
M^{\zeta^{(k)}}_{Q^{(k)}}(\alpha^{(k)})
\right)$, we define an
integral of $C$ along the virtual fundamental cycle
$\left[ \prod_{k=1}^{j}
M^{\zeta^{(k)}}_{Q^{(k)}}(\alpha^{(k)}) \right]^{vir}$
as an element 
\[
\int_{\left[ \prod_{k=1}^{j}
M^{\zeta^{(k)}}_{Q^{(k)}}(\alpha^{(k)}) \right]^{vir}} 
C
=
\left( \text{id} \otimes \bigotimes_{k=1}^{j} 
\iota^{(k)}_{\ast} \right)^{-1}
\left ( \Pi_{\ast} \left( C \cap 
\left[ 
M^{\zeta}_{Q}(\alpha) \times \prod_{k=1}^{j}
M^{\zeta^{(k)}}_{Q^{(k)}}(\alpha^{(k)}) 
\right]^{vir} 
\right) \right )
\]
in 
$T_{\mb T}^{-1} CH^{\bullet}_{\mb T} 
(M^{\zeta}_{Q}(\alpha))$.
Here $\iota^{(k)} \colon M^{0}
(\alpha^{(k)})^{\mb T} 
\to M^{0}(\alpha^{(k)})$ is the embedding for $k=1, \ldots,
j$.

\begin{defn}
For a cohomology class $C \in A^{\bullet} (M_{\mbi d})$, 
we define an $K$-theoretic integral of $C$ along the 
virtual fundamental cycle $\left[ \prod_{k=1}^{j}
M^{\zeta^{(k)}}_{Q^{(k)}}(\alpha^{(k)}) \right]^{vir}$
as an element
\[
\intk_{\left[ \prod_{k=1}^{j}
M^{\zeta^{(k)}}_{Q^{(k)}}(\alpha^{(k)}) \right]^{vir}} C
=
\int_{\left[ \prod_{k=1}^{j}
M^{\zeta^{(k)}}_{Q^{(k)}}(\alpha^{(k)}) \right]^{vir}} 
C \cdot 
\boxtimes_{k=1}^{j} \td_{M_{Q^{(k)}}
^{\zeta^{(k)}}(\alpha^{(k)})}
\]
in 
$T_{\mb T}^{-1} CH^{\bullet}_{\mb T} 
(M^{\zeta}_{Q}(\alpha))$.
\end{defn}

For later use, we set $\mb T_{0}=\mb T \times 
\C^{\ast}_{t}$, and consider 
$\mb T_{0}$-action on $M^{\zeta}(\alpha)$ via
the projection $\mb T_{0}=\mb T \times \C^{\ast}_{t} \to 
\mb T$, that is, the component $\C^{\ast}_{t}$ trivially 
acts.  
For a $\mb T$-equivariant vector bundle $\E$ with the 
rank $r$, we set 
\[
\wedge_{-t} \mc E = \sum_{i=0}^{r} (-1)^{i} 
\wedge^{i} \left( \mc E \otimes \C_{t} \right)
= \sum_{i=0}^{r} (-t)^{i} \wedge^{i} \mc E 
\]
For $K$-theory class $\Lambda= [\E ]- [\F]$ 
with vector bundles $\E, \F$, we set 
$\wedge_{-t} \Lambda = \wedge_{-t} E/ \wedge_{-t} F
\in S_{\mb T_{0}}^{-1} R_{\mb T_{0}}$.
We further define the $K$-theoretic Euler class of 
$\Lambda= [\E ]- [\F] \in K_{\mb T}(M^{\zeta}
(\alpha))$ by
\[
\eukt (\Lambda ) = \ch (\wedge_{-t} \Lambda^{\vee}) =
{\ch \left ( \wedge_{-t} \E^{\vee} \right) \over 
\ch \left( \wedge_{-t} \F^{\vee} \right) }.
\in T_{\mb T_{0}}^{-1} \widehat{A}^{\bullet}
_{\mb T_{0}}(M^{\zeta}(\alpha)).
\]
When we can substitute $t=1$, we set 
$\wedge_{-1} \Lambda = \wedge_{-t} \Lambda 
\Big|_{t=1}$ and 
$\euk \left( \Lambda \right) = 
\euk^{t} \left( \Lambda \right) 
\Big |_{t=1}$.

For polynomials $F \left( \left(
X_{i,1}, \ldots, X_{i,m} 
\right)_{i\in I}, \left(
Y_{i,1}, \ldots, Y_{i,m} 
\right)_{i\in I} \right)$
with the coefficients in 
$S_{\mb T_{0}} ^{-1} 
R_{\mb T_{0}}$ 
and 
$I$-graded vector bundle $\mc E=\bigoplus_{i \in I} 
\E_{i}$, set
\[
F \left(\mc E \right):=
F \left( \left(
\wedge^{1} \mc E_{i}, \ldots, \wedge^{m} \mc E_{i} 
\right)_{i \in I}, 
\left(
\wedge^{1} \mc E_{i}^{\vee}, \ldots, \wedge^{m} 
\mc E_{i}^{\vee} \right)_{i \in I} \right).
\]
In more general, we consider $K$-theory classes 
$F(\mc V)$ for tautological bundles 
$\mc V$ on $M^{\zeta}(\alpha)$ and cohomology classes
$C=\ch F(\mc V)$.
Then we have 
$\intk_{M^{\zeta}(\alpha)} \ch 
{\color{black}F(\mc V) }\in 
T_{\mb T_{0}}^{-1} A_{\mb T_{0}} \cong S_{\mb T_{0}}^{-1}
R_{\mb T_{0}}$.
When $F(\mc V)=\wedge_{-t} \Lambda^{\vee}$ for 
$\Lambda=\Lambda(\mc V)$, we have 
$\ch F(\mc V) = \eukt(\Lambda)$.

For simplicity, we assume that 
$M^{\zeta}_{Q}(\alpha), 
M^{\zeta^{\sharp}}_{Q^{\sharp}}(\alpha^{\sharp})$ are 
smooth, 
the fixed points sets
$M^{\zeta}_{Q}(\alpha)^{\mb T}, M^{\zeta^{\sharp}}
_{Q^{\sharp}}(\alpha^{\sharp})$ are finite, and the 
virtual fundamental cycles 
$[M^{\zeta}_{Q}(\alpha)]^{vir}$ coincide with 
usual ones $[M^{\zeta}_{Q}(\alpha)]$.
We note that $[M^{\zeta^{\sharp}}_{Q^{\sharp}}
(\alpha^{\sharp})]^{vir}$ 
can be different from 
$[M^{\zeta^{\sharp}}_{Q^{\sharp}}(\alpha^{\sharp})]$.
We often identify
cohomology classes $C$ on $M^{\zeta}_{Q}(\alpha)$ and 
the Poincare dual $C \cap 
[M^{\zeta}_{Q}(\alpha)]$.


\subsection{Localization for an enhanced master space $\mc M$}
\label{subsec:local1}
For a wall-crossing formula, we consider 
$\C^{\ast}_{\hbar}$-action on 
$\mc M$ defined in \eqref{fiberwise}, and set
$\mb T_{1}=\mb T \times \C^{\ast}_{t} \times 
\C^{\ast}_{\hbar}$.
Using \cite[(1)]{GP}, we consider the following 
commutative diagram :
\[
\xymatrix{
K^{\mb T_{1}} (\mc M ) \ar[d]_{\Pi_{\ast}} \ar[r] & 
S_{\mb T_{1}}^{-1} K^{\mb T_{1}} 
(\mc M)  \ar[d] \ar[r] & 
S_{\mb T_{1}}^{-1} K^{\mb T_{1}} 
( \mc M^{\C^{\ast}_{\hbar}}) \ar[d] \ar[r]^-{\tau} &
T_{\mb T_{1}}^{-1} 
\widehat{CH}_{_{\mb T_{1}}}^{\bullet} 
( \mc M^{\C^{\ast}_{\hbar}}) \ar[d]
\\
K^{\mb T_{1}} (M_{0}) \ar[r] & 
S_{\mb T_{1}}^{-1} K^{\mb T_{1} } (M_{0}) 
\ar@{=}[r] & 
S_{\mb T_{1}}^{-1} K^{\mb T_{1} } ( M_{0} ) 
\ar[r]^-{\tau} &
T_{\mb T_{1}}^{-1} \widehat{CH}_{\mb T_{1}}^{\bullet} 
( M_{0}) 
}
\]
The middle upper horizontal arrow in the above diagram is given by 
\[
{\iota_{+}^{\ast} \over \wedge_{-1} (N_{+})^{\vee}} + 
{\iota_{-}^{\ast} \over \wedge_{-1} (N_{-})^{\vee}} + 
\sum_{\mk I^{\sharp} \in 
\mc D_{\beta_{\ast}}^{\ell}(\mk I)} 
{\iota_{\mk I^{\sharp}}^{\ast} \over 
\wedge_{-1} (N_{\mk I^{\sharp}})^{\vee}}.
\]
Here $\iota_{\pm}$ and $\iota_{\mk I^{\sharp}}$ are 
embeddings of $\mc M_{\pm}$ and $\mc M_{\mk I^{\sharp}}$ 
into $\mc M$.

We consider a $K$-theory class $F(\mc V) \in 
K_{\mb T_{1}}(\mc M)$ as above, and  
consider $C=\ch F(\mc V)$.
By the above diagram, we have
\begin{align}
\label{localization}
\intk_{[\mc M]^{vir}} \ch F(\mc V) =  
\intk_{[\mc M_{+}]^{vir}} {\ch F(\mc V)|_{\mc M_{+}}
\over \euk (N_{+})} + 
\intk_{[\mc M_{-}]^{vir}} {\ch F(\mc V)|_{\mc M_{-}}
\over \euk (N_{-})} + 
\sum_{\mk I^{\sharp} \in 
\mc D_{\beta_{\ast}}^{\ell}(\mk I)} 
\intk_{[\mc M_{\mk I^{\sharp}}]^{vir}} 
{\ch F(\mc V)|_{\mc M_{\mk I^{\sharp}}} \over 
\euk (N_{\mk I^{\sharp}})}.
\end{align}
Since $\intk_{[\mc M]^{vir}} \ch F(\mc V)$ is holomorphic in $\hbar$, we have 
\begin{align}
\label{localization}
 0= 
\oint d\hbar \intk_{[\mc M_{+}]^{vir}} 
{\ch F(\mc V)|_{\mc M_{-}} \over \euk (N_{+})} + 
\oint d\hbar \intk_{[\mc M_{-}]^{vir}} 
{ \ch F(\mc V)|_{\mc M_{-}} \over \euk (N_{-})} + 
\sum_{\mk I^{\sharp} \in 
\mc D_{\beta_{\ast}}^{\ell}(\mk I)}
\oint d\hbar \intk_{[\mc M_{\mk I^{\sharp}}]^{vir}} 
{\ch F(\mc V)|_{\mc M_{\mk I^{\sharp}}} 
\over \euk (N_{\mk I^{\sharp}})}
\end{align}
for contour integrals over arbitrary closed path.

We set $\wt {F}(\mc V)
= F(\mc V) \cdot \wedge_{-t} \Theta
(\mc F_{\bullet})^{\vee}$.
When we consider \eqref{cls1}, we set $t=1$.
Then substituting $\wt C= \ch \wt{F}(\mc V) / 
|\mk I|_{t}!$ 
we compute each term in \eqref{localization}. 
Here we introduce $[m]_{t} =(1-t^{m} ) / (1-t)$ for a non-negative integer $m$,
and we set $|\mk I|_{t} = [|\mk I|]_{t}$ for a finite set $\mk I$.
We also set $[m]_{t} ! = [1]_{t} \cdot [2]_{t} \cdots [m]_{t}$, and
$c=\sum_{i \in \wt{Q}_{0}\setminus \lbrace \infty \rbrace} 
(\theta^{+}_{i} - \theta^{-}_{i}) c_{1}(\mc V_{i})$.
{\color{black}
To compute further, we delete line bundles $L$ and  integers $m^{\sharp}$ in Theorem \ref{thm:decomp} 
(iii) and Lemma \ref{nmki} by taking residues 
$\displaystyle \oint f(\hbar) d\hbar = m^{\sharp} \oint f( m^{\sharp} \hbar+a) d\hbar$ 
(cf. \cite[\S 8.2]{O1}). 
Then we have
}
\begin{align*}
\oint d\hbar 
\intk_{[\mc M_{-}]^{vir}} 
{\wt C|_{\mc M_{-}} \over \euk( N_{-})}
&=
\oint d\hbar \intk_{[\wt{M}^{0}(\alpha, \mk I)]^{vir}} 
{\ch F(\mc V {\color{black}\otimes e^{\hbar}}) 
\over 1 - e^{\color{black} \hbar }}
\cdot { \eukt( \Theta(\mc F_{\bullet} )) 
\over |\mk I|_{t} !}
\\
&= 
\oint d\hbar \intk_{[M^{\zeta^{-}}(\alpha)]^{vir}} 
{ \ch F(\mc V {\color{black} \otimes e^{\hbar}} ) 
\over 1 - e^{\color{black} \hbar }}
\\
&=
-\intk_{[M^{\zeta^{-}}(\alpha)]^{vir}} \ch F(\mc V) 
=
-\intk_{[M^{\zeta^{-}}(\alpha)]^{vir}} C.
\end{align*}
Similarly we have
\begin{align*} 
\oint d\hbar \intk_{[\mc M_{+}]^{vir}} 
{\wt C|_{\mc M_{+}} \over \euk( N_{+})}
&=
\oint d\hbar 
\intk_{[\wt{M}^{\ell}(\alpha, \mk I)]^{vir}} 
{\ch F(\mc V {\color{black} \otimes e^{\hbar}} )
\over 1 - e^{ \color{black} -\hbar}}
\cdot { \eukt ( \Theta(\mc F^{\bullet} )) 
\over |\mk I|_{t}!}
\\
&= 
{1 \over |\mk I|_{t}! }
\intk_{[\wt M^{\ell}(\alpha, \mk I)]^{vir}} 
\ch F(\mc V) \cdot \eukt ( \Theta(\mc F^{\bullet} )) 
=
\intk_{[\wt M^{\ell}(\alpha, \mk I)]^{vir}} \wt C.
\end{align*}
Thus we have
\begin{align}
\label{euler0}
\int_{[\wt M^{\ell}(\alpha, \mk I)]^{vir}} \wt C
- 
\int_{[M^{\zeta^{-}}(\alpha)]^{vir}} C
= - \sum_{\mk I^{\sharp} \in 
\mc D_{\beta_{\ast}}^{\ell}(\mk I)} 
\oint d\hbar \intk_{[\mc M_{\mk I^{\sharp}}]^{vir}} 
{\ch \wt{F}(\mc V))|_{\mc M_{\mk I^{\sharp}}}
\over \euk(N_{\mk I^{\sharp}})}.
\end{align}
This coefficient $m^{\sharp}$ vanishes by comparing integrals over $\mc M_{\mk I_{\sharp}}$ and 
$Fl_{H_{Q}(d^{\sharp} \beta)}(\bar{\mc V}^{\sharp}, \bar{\mk I}^{\sharp})
\times
\wt{M}^{\min(\mk I^{\sharp})-1} (\alpha -d^{\sharp} \beta, \mk I^{\flat})$ via $\mc M'_{\mk I}$ in Theorem \ref{thm:decomp}.
Then together with \eqref{rel}, 
the last summand in \eqref{euler0} is equal to
\begin{align}
\label{euler}
&
\oint d\hbar
{|\mk I^{\flat}|_{t}! \over |\mk I|_{t}!} 
\intk_{[\wt{M}^{\min(\mk I^{\sharp})-1} 
(\alpha -d^{\sharp} \beta, \mk I^{\flat})]^{vir}}
{1 \over |\mk I^{\flat}|_{t}!}
\notag \\
&\cdot
\intk_{[Fl_{H_{Q}(d^{\sharp} \beta)}
(\bar{\mc V}^{\sharp}, \bar{\mk I}^{\sharp})]^{vir}} 
{\ch \left[F(\mc V^{\flat} + \mc V^{\sharp} 
\otimes \C_{e^{\hbar}}) \cdot \wedge_{-t} 
\Theta( \mc F^{\sharp}_{\bullet} 
\otimes \C_{e^{\hbar}} + \mc F^{\flat}_{\bullet} )^{\vee} 
\right]
\over 
\euk \left[ \mk N(\mc V^{\sharp} \otimes \C_{e^{\hbar}}, 
\mc V^{\flat}) 
+ \mk H( \mc F^{\sharp}_{\bullet} \otimes \C_{e^{\hbar}}, 
\mc F^{\flat}_{\bullet} ) \right]
}.
\end{align}
Here we have  
$\Theta( \mc F^{\sharp}_{\bullet} \otimes \C_{e^{\hbar}} + \mc F^{\flat}_{\bullet} )=
\Theta(\mc F_{\bullet}^{\sharp}) +\Theta(\mc F_{\bullet}^{\flat})
+\mk H( \mc F^{\sharp}_{\bullet} \otimes \C_{e^{\hbar}}, \mc F^{\flat}_{\bullet} )$, 
and further decomposition $\Theta(\mc F_{\bullet}^{\sharp})=\Theta(\bar{\mc F}_{\bullet}^{\sharp}) + 
\bar{\mc V}^{\sharp}$ with the relative tangent bundle 
$\Theta(\bar{\mc F}_{\bullet}^{\sharp})$ of
$Fl_{H_{Q}(d^{\sharp} \beta)}(\bar{\mc V}^{\sharp}, \bar{\mk I}^{\sharp})$
over $H_{Q}(d^{\sharp} \beta)$ from \eqref{ob4}.

\subsection{Iterated classes}
\label{subsec:iter1}

In the following, we deduce wall-crossing formula from analysis 
in the previous section.
These are the similar calculations to \cite[\S 6]{NY2}, 
hence we omit detail explanation.

In the following, we fix a $Q_{0}$-graded vector space 
$V=\bigoplus_{v \in Q_{0}} V_{v}$ 
with $\dim V_{\infty}=1$.
We set $\alpha=(\alpha_{v})_{v\in Q_{0}} = 
(\dim V_{v})_{v \in Q_{0}} \in (\Z_{\ge 0})^{Q_{0}}$.
For $\mbi d= (d^{(1)}, \ldots, d^{(j)}) \in 
\Z_{>0}^{j}$, 
we set $|\mbi d| = d^{(1)} + \cdots + d^{(j)}$.
Let $\Dec{}_{\beta_{\ast}, j}^{\alpha_{\ast}}$ 
be the set of collections 
$\mbi {\mk I}=( \mk I^{(1)}, \ldots, \mk I^{(j)})$ 
(decomposition data) such that 
\begin{enumerate}
\item[$\bullet$]
$\mk I^{(1)}, \ldots, \mk I^{(j)}$ are 
disjoint non-empty subsets of $[\alpha_{\ast}]=
\lbrace 1, 2, \ldots, \alpha_{\ast} \rbrace$, 
\item[$\bullet$] $|\mk I^{(i)}| = d^{(i)} 
\beta_{\ast}$ with $d^{(i)} \in \Z_{>0}$
for $i=1, \ldots,j$, and
\item[$\bullet$] $\min(\mk I^{(1)}) > \cdots > 
\min(\mk I^{(j)})$.
\end{enumerate}
We set $\mbi d_{\mbi{\mk I}} = 
\frac{1}{\beta_{\ast}}(|\mk I^{(1)}|, \ldots, 
|\mk I^{(j)}|)
=(d^{(1)}, \ldots, d^{(j)}) \in \Z_{>0}^{j}$
and $\mk I^{\infty} =[\alpha_{\ast}] \setminus 
\bigsqcup_{i=1}^{j} \mk I^{(i)}$.
We note that $\Dec{}_{\beta_{\ast}, 1}^{\alpha_{\ast}} 
= \mc D^{\ell}_{\beta_{\ast}}([\alpha_{\ast}])$
and
\[
\Dec{}_{\beta_{\ast}, j+1}^{\alpha_{\ast}}
=
\lbrace
(\mk I^{(1)}, \ldots, \mk I^{(j)}, \mk I^{(j+1)})
\mid 
\mbi{\mk I}=(\mk I^{(1)}, \ldots, \mk I^{(j)}) \in 
\Dec{}_{\beta_{\ast}, j}^{\alpha_{\ast}}, \ 
\mk I^{(j+1)} \in \mc D^{\min(\mk I^{(j)})-1}
_{\beta_{\ast}}( \mk I^{\infty}) \rbrace.
\]
We define $\sigma \colon 
\Dec{}_{\beta_{\ast}, j+1}^{\alpha_{\ast}} 
\to \Dec{}_{\beta_{\ast}, j}^{\alpha_{\ast}}$ by
$\sigma(\mk I^{(1)}, \ldots, \mk I^{(j)}, 
\mk I^{(j+1)}) 
= (\mk I^{(1)}, \ldots, \mk I^{(j)})$.
When $\beta_{\ast}=1$, we set 
$\Dec^{\alpha_{\ast}}_{j}=\Dec^{\alpha_{\ast}}_{\beta_{\ast}, j}$.

\indent For $\mbi d =(d^{(1)}, \ldots, d^{(j)})\in \Z_{>0}^{j}$ and $\ell=0, 1, \ldots, \alpha_{\ast}$, we
set $H_{\mbi d}=\prod_{i=1}^{j} H_{Q}(d^{(i)} \beta)$, and consider a product 
\[
M_{\mbi d} = M^{\zeta^{-}} (\alpha - |\mbi d| \beta) 
\times H_{\mbi d}
\] 
We write by $\mc V^{(i)}$ the tautological bundle
$\mc V^{\sharp} \oplus \mc V^{\sharp}_{\infty^{\sharp}}$ on the component 
$H_{Q}(d^{(i)} \beta)$.

For $\mbi{ \mk I}=(\mk I^{(1)}, \ldots, \mk I^{(j)}) \in 
\Dec{}_{\beta_{\ast}, j}^{\alpha_{\ast}}$, we set
$Fl_{H_{\mbi d}}(\mbi{\mk I})=\prod_{i=1}^{j} 
Fl_{H_{Q}(d^{(i)} \beta)}(\bar{\mc V}_{\ast}^{(i)}, 
\bar{\mk I}_{i})$
\[
\wt M_{\mbi{\mk I}}^{\ell}  = \wt M^{\ell}(\alpha - 
|\mbi d_{\mbi{\mk I}}| \beta, \mk I^{\infty}) \times
Fl_{H_{\mbi d}} (\mbi{\mk I})
\]
where $\bar{\mc V}_{\ast}^{(k)}=\mc V_{\ast}^{(k)} / 
\mc V^{(k)}_{\infty^{\sharp}}$, 
$\bar{\mk I}_{i}=\mk I^{(i)} \setminus \lbrace 
\min(\mk I^{(i)})\rbrace$,
and write by $\mc F_{\bullet}^{(i)}$ the pull-back to 
$\wt M_{\mbi{\mk I}}^{\ell}$ of 
the universal flag on each component 
$Fl_{H_{Q}(d^{(i)} \beta)}(\bar{\mc V}_{\ast}^{(i)}, 
\bar{\mk I}_{i})$.
We regard $\lbrace \mc V_{\ast (k)} 
\rbrace_{k=1}^{\alpha_{\ast}}$
on $\wt M^{\ell}( \alpha - |\mbi d_{\mbi{\mk I}}| \beta, 
\mk I^{\infty} )$
as a flag $\mc F_{\bullet}$ of $\mc V_{\ast}$, and  
write by $\mc F_{\bullet}^{\infty}$ the pull-back to 
$\wt M_{\mbi{\mk I}}^{\ell}$.
For the tautological bundle $\mc V$ on 
$\wt M^{\ell}(\alpha - |\mbi d_{\mbi{\mk I}}| \beta, 
\mk I^{\infty} )$, 
we write by the same letter 
the pull-backs to the product 
$\wt M_{\mbi{\mk I}}^{\ell}$.
For $\mc V^{(i)}$, we also write by the same letter the 
pull-back by the projection 
$\wt M_{\mbi{\mk I}}^{\ell} \to H_{Q}(d^{(i)} \beta)$.

We set $\mb T_{j}=\mb T \times \C^{\ast}_{t} 
\times \prod_{i=1}^{j} \C^{\ast}_{\hbar_{i}}$.
We set $\mc F_{\bullet}^{>i}=
\mc F_{\bullet}^{\infty} \oplus \bigoplus_{k > i} 
\mc F_{\bullet}^{(k)} \otimes e^{\hbar_{k}}$, 
$\mc V^{>i} = \mc V \oplus \bigoplus_{k > i } 
\mc V^{(k)} \otimes e^{\hbar_{k}}$, and  
$\mbi{\mk I}^{> i} = \mk I^{\infty} \sqcup 
\bigsqcup_{k >i} \mk I^{(k)}$.
Here $e^{\hbar_{i}}$ denotes the weight space 
$\C_{e^{\hbar_{i}}}$ with the 
weight $e^{\hbar_{i}} \in \C^{\ast}_{\hbar_{i}}$.

For $\mbi{\mk I}=(\mk I^{(1)}, \ldots, \mk I^{(j)}) \in 
\text{Dec}_{\beta_{\ast}, j}^{\alpha_{\ast}} $, we set
\begin{align}
\wt C_{\mbi{\mk I}} (\mc V)
=&
{1 \over |\mk I^{\infty}|_{t}!}
\intk_{[F^{H}_{\mbi{\mk I}}]^{vir}} 
{ 
\ch F \left( \mc V \oplus 
\bigoplus_{i=1}^{j} \mc V^{(i)} \otimes e^{\hbar_{i}} 
\right) 
\euk^{t}(\Theta(\mc F^{>0}_{\bullet}))
\over 
\prod_{i=1}^{j}
\euk \left( \mathfrak{N} ( 
\mc V^{(i)} \otimes e^{\hbar_{i}}, 
\mc V^{>i} ) \otimes \chi_{Q}
\right)
\euk(\mk H(\mc F^{(i)}_{\bullet} \otimes e^{\hbar_{i}}, 
\mc F^{>i}_{\bullet}))
}
\label{itcoho}
\end{align}
in $A^{\bullet}_{\mb{T}_{j}}(\wt M^{\ell}(\alpha - 
|\mbi d_{\mbi{\mk I}}| \beta ))$.
Here $\mathfrak{N} ( 
\mc V^{(i)} \otimes e^{\hbar_{i}}, 
\mc V^{>i} ) \otimes \chi_{Q}$ is a modification of
\eqref{normal} by
\begin{align*}
\mathfrak{N}( \mc V^{\sharp}, \mc V^{\flat})
\otimes \chi_{Q}
&=
\sum_{(\heartsuit, \spadesuit)} 
\left(
\sum_{\substack{h \in Q_{1} }} 
\mc H om( \mc V^{\heartsuit}_{\out(h)}, 
\mc V^{\spadesuit}_{\inn(h)}) \otimes \C_{\chi_{h}}
\right.
\\
&\left.
- 
\sum_{\gamma \in Q_{2}} \mc H om( 
\mc V^{\heartsuit}_{\out(\gamma)}, 
\mc V^{\spadesuit}_{\inn(\gamma)})  \otimes 
\C_{\chi_{\gamma}}
-
\sum_{i \in I} 
\mc Hom ( \mc V^{\heartsuit}_{i}, 
\mc V^{\spadesuit}_{i}) \right),
\end{align*}
where the sum is taken for only cross terms 
$(\heartsuit, \spadesuit)= (\sharp, \flat), 
(\flat, \sharp)$.
In the numerator of \eqref{itcoho}, we have further 
decomposition
\[
\euk^{t}(\Theta(\mc F^{>0}_{\bullet}))
=
\euk^{t}(\Theta(\mc F^{\infty}_{\bullet}))
\prod_{i=1}^{j}
\euk^{t}(\Theta(\mc F^{(i)}_{\bullet}))
\euk^{t} (\mk H ( \mc F^{(i)}_{\bullet}, 
\mc F^{>i}_{\bullet})), 
\]
which is useful for recursive computations.

\subsection{Recursions}
\label{subsec:loca}
For $\mbi{\mk I} = (\mk I^{(1)}, \ldots, \mk I^{(j)}) 
\in \Z_{>0}^{j}$, we apply \eqref{euler0} to $\mk I = 
\mk I^{\infty}$, and 
$\ell=\min(\mk I^{(j)})-1$, and take 
an equivariant cohomology class 
$\wt C = \wt C_{\mbi{\mk I }}(\mc V)$ on $\mc M$.
{\color{black} When $j=0$, we
	set $\mbi{\mk I}=()$ and 
	apply \eqref{euler0} to $\mk I = 
	\lbrace 1, \ldots, \alpha_{\ast} 
	\rbrace$ and 
	$\ell=\alpha_{\ast}$.} 
Then \eqref{euler} is equal to
\begin{align}
\label{exc}
\oint d \hbar_{j+1}
{|\mk I^{\flat}|_{t}! \over |\mk I|_{t} !} 
\intk_{[\wt M^{\min(\mk I^{\sharp})-1}
(\alpha - d\beta - d^{\sharp} \beta )]^{vir} } 
\wt C_{ ( \mbi{\mk I}, \mk I^{\sharp} )}(\mc V),
\end{align}
where $d^{\sharp} = |\mk I^{\sharp}|/ \beta_{\ast}$, and $(\mbi{\mk I}, \mk I^{\sharp})
=(\mk I^{(1)}, \ldots, \mk I^{(j)}, \mk I^{\sharp}) \in \Dec^{\alpha_{\ast}}_{\beta_{\ast}, j+1}$.

Using $\wt C_{\mbi{\mk I}} (\mc V)$ defined in \eqref{itcoho},
we deduce recursion formula.
\begin{lem}
For $l \ge 1$, we have
\begin{align}
\notag
&
\intk_{[M^{\zeta^{+}}(\alpha)]^{vir}} 
\ch F(\mc V)  - 
\intk_{[M^{\zeta^{-}}(\alpha )]^{vir}} 
\ch F(\mc V) \\
\notag
&= 
\sum_{j=1}^{l-1} (-1)^{j} 
\oint d \hbar_{1} \cdots \oint d\hbar_{j}
\sum_{\mbi{\mk I} \in \Dec_{\beta_{\ast}, j}^{\alpha_{\ast}}} \frac{ |\mk I^{\infty}|_{t}! }{[ \alpha_{\ast}]_{t}! }
\intk_{[\wt M^{ 0}(\alpha - |\mbi d_{\mbi{\mk I}}|
\beta)]^{vir}} \wt C_{\mbi{\mk I}} (\mc V)\\
\label{formula2>0}
&+
(-1)^{l} \oint d \hbar_{1} \cdots \oint d\hbar_{l}
\sum_{\mbi{\mk I} \in 
\Dec{}_{\beta_{\ast}, l}^{\alpha_{\ast}}} 
{|\mk I^{\infty}|_{t}! \over [\alpha_{\ast}]_{t}! }
\intk_{[\wt M^{\min(\mk I_{l})-1} 
(\alpha - |\mbi d_{\mbi{\mk I}}| \beta )]^{vir}} 
\wt C_{\mbi{\mk I}} (\mc V).
\end{align}
\end{lem}
\proof
We prove by induction on $j$.
For $l=1$, \eqref{formula2>0} is nothing but \eqref{exc} for $j_{0}=0$ and $\ell = n$.
For $l \ge 1$, we assume the formulas \eqref{formula2>0}.
Then again by \eqref{exc}, the last summand for each $\mbi{\mk I} \in \Dec{}_{\beta_{\ast}, l}^{\alpha_{\ast}}$ is equal to 
\begin{align*}
&
{|\mk I^{\infty}|_{t}! \over [\alpha_{\ast}]_{t}! }
\left( \intk_{[\wt M^{0} (\alpha - 
|\mbi d_{\mbi{\mk I}}| \beta )]^{vir}} 
\wt C_{\mbi{\mk I}} (\mc V)  \right. \\
&-
\oint d\hbar_{l+1}
\left. \sum_{\mbi{\mk J} \in \sigma^{-1} (\mbi{\mk I})} 
\frac{| \mk J^{\infty}|_{t}! }{|\mk I^{\infty}|_{t}!} 
\intk_{[\wt M^{\min (\mk J^{(l+1)} )-1 }
(\alpha - |\mbi d_{\mbi{\mk J}}| \beta )]^{vir}} 
\wt C_{\mbi{\mk J}} (\mc V)  \right ),
\end{align*}
where $\mbi{\mk J}=(\mbi{\mk I}, \mk J^{(l+1)})$.
Hence we have \eqref{formula2>0} for general $l \ge 1$.
\endproof

For $l > \alpha_{\ast} / \beta_{\ast}$, the set 
$\Dec{}_{\beta_{\ast}, l}^{\alpha_{\ast}}$ is empty.
Thus we get the following theorem, 
{\color{black} which gives a proof of
Theorem \ref{thm:introwcm}
when $F(\mc V)=\wedge_{-t} \Lambda^{\vee}$ for some $\Lambda=\Lambda(\mc V)$.}

\begin{thm}
\label{thm:wcm}
We have
\begin{align}
&
\intk_{[M^{\zeta^{+}}(\alpha)]^{vir}} 
\ch F(\mc V) - 
\intk_{[M^{\zeta^{-}}(\alpha )]^{vir}} 
\ch F(\mc V)
\notag
\\
&=
\label{main>0}
\sum_{j=1}^{\lfloor \alpha_{\ast}/\beta_{\ast} \rfloor } (-1)^{j}
\sum_{\mbi{\mk I} \in \Dec{}_{\beta_{\ast}, j}
^{\alpha_{\ast}}} 
{ |\mk I^{\infty}|_{t}! \over [\alpha_{\ast}]_{t}! }
\oint d\hbar_{1} \cdots \oint d\hbar_{j}
\intk_{[\wt M^{0} (\alpha - |\mbi d_{\mbi{\mk I}}| 
\beta )]^{vir}} 
\wt C_{\mbi{\mk I}} (\mc V)
\end{align}
where $\wt C_{\mbi{\mk I}} (\mc V)$ is defined 
in \eqref{itcoho}.
\end{thm}

When $t=1$, we set
\begin{align}
C_{\mbi{d}} (\mc V)
=&
\intk_{[H_{\mbi d}]^{vir}} 
{ 
\ch F \left( \mc V \oplus 
\bigoplus_{i=1}^{j} 
\mc V^{(i)} \otimes e^{\hbar_{i}} \right)
\euk \left( 
\mc V_{\ast}^{(k)} / \mo_{M_{\mbi d}} \right)
\over 
\prod_{i=1}^{j} \euk \left( \mathfrak{N} ( 
\mc V^{(i)} \otimes e^{\hbar_{i}}, 
\mc V^{>i} ) \otimes \chi_{Q}
\right)
}
\label{itclasst=1}
\end{align}
for $\mbi d=(d^{(1)}, \ldots, d^{(j)}) \in \Z^{j}_{>0}$. 
Then we can simplify 
\eqref{main>0} using the following lemma
\begin{lem}
For $\mbi d = (d^{(1)}, \ldots, d^{(j)}) \in 
\Z_{>0}^{j}$ with $|\mbi d| \le \alpha_{\ast} / 
\beta_{\ast}$, we have 
\begin{align*}
\left| 
\left\lbrace \left.
\mbi{\mk I}=(\mk I^{(1)}, \ldots, \mk I^{(j)})
\in \Dec{}_{\beta_{\ast}, j}^{\alpha_{\ast}} \ 
\right| \ \mbi d_{\mbi {\mk I}}=\mbi d \right\rbrace 
\right| 
&=
{1 \over \prod_{i=1}^{j} \lvert \mbi d\rvert_{i} } 
{ \alpha_{\ast}! \over |\mk I^{\infty}|! 
\prod_{k=1}^{j} 
(|\mk I^{(k)}|-1)! },
\end{align*}
where we can take any 
$\mbi{\mk I}=(\mk I^{(1)}, \ldots, \mk I^{(j)})$ with  
$\mbi d_{\mbi {\mk I}}=\mbi d$ in the right hand side, 
and 
$\lvert \mbi d\rvert_{i}=d^{(1)} + \cdots + 
d^{(i)}$and 
$\mbi{\mk I}=(\mk I^{(1)}, \ldots, \mk I^{(j)}) \in 
\rho^{-1}(\mbi d)$.
\end{lem}
\proof
This follows from \cite[Lemma 6.8]{NY2} since 
$|\mk I_{k}| = d_{k} \beta_{\ast}$.
\endproof
Then we have
\begin{align}
&
\intk_{[M^{\zeta^{+}}(\alpha)]^{vir}} 
\ch F(\mc V)
- \intk_{[M^{\zeta^{-}}(\alpha )]^{vir}} 
\ch F (\mc V)
\notag
\\
&=
\label{main>1}
\sum_{j=1}^{\lfloor \alpha_{\ast}/\beta_{\ast} \rfloor } (-1)^{j}
\sum_{\mbi d \in \Z^{j}_{>0}} 
{1 \over \prod_{i=1}^{j} \lvert \mbi d\rvert_{i} }
\oint d\hbar_{1} \cdots \oint d\hbar_{j}
\intk_{[M^{\zeta^{-}} 
(\alpha - \lvert \mbi d \rvert \beta )]^{vir}} 
C_{\mbi d} (\mc V).
\end{align}

\subsection{Euler class of the tangent bundle}
\label{subsec:euler2} 
For a $\mb T$-equivariant framed quiver $(Q, \chi_{Q})$
and 
$Q_{0}$-graded vector bundle $\mc V$, we set 
\[
\Lambda_{\adj} = \Lambda_{\adj} (\mc V)=
\sum_{h \in Q_{1}} \mc Hom (\mc V_{\out (h)}, 
\mc V_{\inn (h)}) 
{\color{black}
\otimes \C_{\chi_{h}}
}
- \sum_{\gamma \in Q_{2}} 
\mc Hom (\mc V_{\out (\gamma)}, \mc V_{\inn (\gamma)})
{\color{black}
\otimes \C_{\chi_{\gamma}}
}
- \sum_{i \in I} \mc End(\mc V_{i}).
\]
{\color{black}
To evaluate \eqref{main>0}, we use the following lemma.
\begin{lem}
\label{lem:residue1}
We have
\[
\oint d \hbar \prod_{k=1}^{M} 
{1-e^{\hbar + b_{k}} \over 1 - e^{\hbar + a_{k}}}
\prod_{l=1}^{N} 
{1-e^{-\hbar + d_{l}} \over 1 - e^{-\hbar + c_{l}}}
=
e^{\sum_{k=1}^{M} b_{k} - a_{k}} - 
e^{\sum_{l=1}^{N} d_{l} - c_{l}}.
\]
\end{lem}
\proof
We can easily reduce to the case
where $N=0$.
To prove this case, we use multiplicative variables $u=e^{\hbar}, z_{k}=e^{b_{k}}$ and $w_{k}=e^{a_{k}}$ for
$k=1, \ldots, M$.
Then we have 
\begin{align*}
\oint d \hbar \prod_{k=1}^{M} \frac{1-e^{\hbar + b_{k}}}{1 - e^{\hbar + a_{k}}}
&=
-
\Res_{u=0, \infty} 
{du \over u} 
\prod_{k=1}^{M} 
{1-u z_{k} \over 1-u w_{k} }\\
&=-1 - \Res_{u=\infty} {du \over u}
\prod_{k=1}^{M} {z_{k} \over w_{k}} {1-u^{-1} z_{k}^{-1} \over 1- u^{-1} w_{k}^{-1} } 
 =-1 + {z_{1} \cdots z_{M}
 	\over w_{1} \cdots w_{M}}.
\end{align*}
\endproof
}

In \eqref{itcoho}, we substitute
\[
\ch F \left( \mc V \oplus \bigoplus_{i=1}^{j} \mc V^{(i)} 
\otimes e^{\hbar_{i}} \right) =
\eukt \left(\Lambda_{\adj}(\mc V)
\right) 
\cdot
\prod_{i=1}^{j} 
\eukt(\Lambda_{\adj}
\left(
\mc  V^{(i)} \otimes e^{\hbar_{i}}
\right)
\cdot 
\eukt \left( \mk N \left( \mc V^{(i)} \otimes 
e^{\hbar_{i}}, \mc V^{>i} \right)
\right),  
\]
and divide the computations into three parts:
\begin{lem}
\label{lem:residue2}
For each $i=1,\ldots,j$, we have the following.\\
(1) We have 
\[
\displaystyle
\oint d\hbar_{i} 
{
\eukt \left( \mk N( \mc V^{(i)} \otimes e^{\hbar_{i}}, 
\mc V^{>i} ) \otimes \chi_{Q} \right)
\over
\euk \left( \mk N( \mc V^{(i)} \otimes e^{\hbar_{i}}, 
\mc V^{>i} ) \otimes \chi_{Q} \right)
}
= 
t^{\color{black} 
\chi \left(d^{(i)} \beta, 
\alpha- d^{\le i} \beta \right)} -
t^{\color{black}
\chi \left( \alpha- d^{\le i} \beta, d^{(i)} \beta \right)}, 
\]
where
the Euler form 
{ \color{black} 
$\chi \left( \alpha', \beta'\right)$ is defined
in \eqref{eulerform}}.
We set 
$d^{\le i}= \sum_{l \le i} d^{(l)}$ for  
$\mbi d_{\mbi{\mk I}}=
\left(d^{(1)}, \ldots, d^{(j)} \right)=
{1 \over \beta_{\ast}}
\left(|\mk I^{(1)}|,\ldots, |\mk I^{(i)}| \right) 
$.
\\
(2) 
$
\int_{[Fl_{H_{Q}( d^{(i)} \beta )}
( \bar{\mc V}^{(i)}_{0}, \bar{\mk I}_{i})]^{vir}}
{\color{black}
\eukt 
\left( \Lambda_{\adj} (
(\mc  V^{(i)} \otimes e^{\hbar_{i}}))
\right) 
}
\eukt 
(\Theta( \mc F^{(i)}_{\bullet} )) 
= 
{[d^{(i)} \beta_{\ast} -1]_{t}! \over 1-t}
\gamma_{d^{(i)}} (
{\color{black}t})$, where 
\begin{align}
\label{gamma}
\gamma_{d^{(i)}} (\theta)
&= 
\int_{[H_{Q}( d^{(i)} \beta )]^{vir}} 
\eukt \left( \Lambda_{Q^{\sharp}} \left(\mc V^{(i)} 
\oplus \mc V^{(i)}_{\infty^{\sharp}} \right) \right),
\\
\notag
{\color{black}
\Lambda_{Q^{\sharp}}(\mc V^{(i)} \oplus 
\mc V^{(i)}_{\infty^{\sharp}}) 
}
&
{\color{black}
=
\sum_{h \in Q_{1} \atop \inn(h), \out(h) \neq \infty} \mc Hom (\mc V^{(i)}_{\out (h)}, 
\mc V^{(i)}_{\inn (h)}) 
\otimes \C_{\chi_{h}}
+
\mc Hom( \mc V^{(i)}_{\infty^{\sharp}}, 
\mc V^{(i)}_{\ast})
}
\\
\label{LQsharp}
&
{\color{black}
- \sum_{ \gamma \in Q_{2}} 
\mc Hom (
\mc V^{(i)}_{\out (\gamma)}, 
\mc V^{(i)}_{\inn (\gamma)})
\otimes \C_{\chi_{\gamma}}
- \sum_{i \in I} \mc End(\mc V^{(i)}_{i}).
}
\end{align}
(3) 
$
\displaystyle
\oint d\hbar_{i} 
\eukt (\mk H(\mc F^{(i)}_{\bullet} \otimes e^{\hbar_{i}}, \mc F^{>i}_{\bullet} ))/
\euk (\mk H(\mc F^{(i)}_{\bullet} \otimes e^{\hbar_{i}}, \mc F^{>i}_{\bullet} ))
 = 
t^{s( \mk I^{(i)}, \mbi{\mk I}^{>i})} - t^{s( \mbi{\mk I}^{>i}, \mk I^{(i)})}$,
where 
\begin{align}
\label{sii}
s(\mk I, \mk I') = 
\left| \lbrace (l, l') \in \mk I \times \mk I' \mid 
l < l' \rbrace \right|.
\end{align}
\end{lem}
\proof
(1)
For $ \mc V^{\sharp}=\mc V^{(i)} \otimes e^{\hbar_{i}}, 
\mc V^{\flat}=\mc V^{>i}$, 
we see that 
$\mathfrak{N}( \mc V^{\sharp}, \mc V^{\flat})$ is equal to 
\begin{align*}
&
\sum_{(\heartsuit, \spadesuit)} 
\left(
\sum_{h \in Q_{1} \atop \inn(h), \out(h) \neq \infty} 
\mc H om( \mc V^{\heartsuit}_{\out(h)}, 
\mc V^{\spadesuit}_{\inn(h)}) 
- 
\sum_{\gamma \in Q_{2}} \mc H om( 
\mc V^{\heartsuit}_{\out(\gamma)}, 
\mc V^{\spadesuit}_{\inn(\gamma)}) 
-
\sum_{v \in I} 
\mc Hom ( \mc V^{\heartsuit}_{v}, 
\mc V^{\spadesuit}_{v}) \right)\\
&+
\sum_{h \in Q_{1} \atop \out(h) = \infty} 
\mc H om( \mc V^{\flat}_{\infty}, 
\mc V^{\sharp}_{\inn(h)}) 
+\sum_{h \in Q_{1} \atop \inn(h) = \infty} 
\mc H om( \mc V^{\sharp}_{\out(h)}, 
\mc V^{\flat}_{\infty}) 
\end{align*}
where the sum is taken for only cross terms $(\heartsuit, \spadesuit)= (\sharp, \flat), (\flat, \sharp)$.
By Lemma \ref{lem:residue1}, we do not need $\otimes \chi_{Q}$, and dimension countings give the assertion.
To prove (2) we remark that \eqref{ob4} implies
\begin{align*}
\Theta(\mc F^{(i)}_{\bullet})
&=
\mc Hom( \mc V^{(i)}_{\infty^{\sharp}}, \mc V^{(i)}_{\ast} ) - \mo_{H_{Q} (d^{(i)} \beta)} 
+ \Theta(\bar{\mc F}^{(i)}_{\bullet}),
\end{align*}
and the right hand side does not include $e^{\hbar_{i}}$.
Together with \eqref{LQsharp}, 
the assertion follows from the integral
{\color{black}
$\int_{Fl(\C^{d^{(i)}\beta_{\ast} -1}, \bar{\mk I}_{(i)})} 
\eukt ( T Fl(\C^{d^{(i)} \beta_{\ast} -1}, \bar{\mk I}_{i}))=
[d^{(i)}\beta_{\ast} -1]_{t}!$.
}
Finally (3) follows from direct computations.
\endproof

By Lemma \ref{lem:residue2}, we compute \eqref{main>0}, 
and we have
\begin{align}
&
\notag
\int_{M^{\zeta^{+}}(\alpha)} \euk^{t}(\Lambda_{\adj})
- \int_{M^{\zeta^{-}}(\alpha)} \euk^{t}(\Lambda_{\adj})
\\
\notag
&=
\sum_{j=1}^{\lfloor \alpha_{\ast} / \beta_{\ast} \rfloor} 
\sum_{\mbi{\mk I} \in \Dec{}_{\beta_{\ast}, j}^{\alpha_{\ast}}} 
\frac{ |\mk I^{\infty}|_{t} !  }{ [\alpha_{\ast}]_{t} ! }
\prod_{i=1}^{j} \frac{[d^{(i)} \beta_{\ast} -1]_{t}!}{t-1} \gamma_{d^{(i)}} ({\color{black} t}) 
\\
\label{formula4}
& \cdot  
( t^{s( \mk I^{(i)}, \mbi{\mk I}^{>i})+ {\color{black}
\chi \left(d^{(i)} \beta, \alpha- d^{\le i} \beta \right)
} } - t^{s( \mbi{\mk I}^{>i},  \mk I^{(i)}) + {\color{black} 
\chi \left( 
\alpha- d^{\le i} \beta,
d^{(i)} \beta \right)
}}) 
\int_{ M^{\zeta^{-}}(\alpha - |\mbi d_{\mbi{\mk I}}| 
\beta )} \euk^{t}(\Lambda_{\adj}),
\end{align}
where $d^{(i)} = |\mk I^{(i)}|/\beta_{\ast}$ and 
$\mbi{\mk I}^{>i} = \bigsqcup_{k >i} \mk I_{k} \sqcup 
\mk I^{\infty}$ 
for $\mbi{\mk I}=( \mk I^{(1)}, \ldots, \mk I^{(j)}) \in 
\text{Dec}_{\beta_{\ast}, j}^{\alpha_{\ast}}$,  
and $\gamma_{d^{(i)}} 
({\color{black} t})$ are defined in 
\eqref{gamma}.

We set $\Dec (\alpha_{\ast})= \bigsqcup_{j=1}
^{\lfloor \alpha_{\ast}/\beta_{\ast} \rfloor} 
\Dec{}_{\beta_{\ast}, j}^{\alpha_{\ast}}$ and, 
we get the following theorem.
\begin{thm}
[Theorem \ref{thm:introadj}]
\label{thm:adjoint}
We have
\begin{align*}
&
\int_{[M^{\zeta^{+}}(\alpha)]^{vir}} 
\euk^{t}(\Lambda_{\adj})
- \int_{[M^{\zeta^{-}}(\alpha)]^{vir}} 
\euk^{t}(\Lambda_{\adj})\\
&=
\sum_{k=1}^{\lfloor \alpha_{\ast}/\beta_{\ast} \rfloor} 
\sum_{ \mbi{\mk I} \in \Dec(\alpha_{\ast}) \atop 
|\mbi d_{\mbi{\mk I}}|=k} 
{ |\mk I^{\infty}|_{t} ! \over [\alpha_{\ast}]_{t} ! }
\prod_{i=1}^{k} 
{[d^{(i)} \beta_{\ast} -1]_{t}! \over t-1} 
\gamma_{d^{(i)}} ({\color{black} t}) 
\\
& \cdot  
( t^{s( \mk I^{(i)}, \mbi{\mk I}^{>i})+ {\color{black}
\chi \left(d^{(i)} \beta, \alpha- d^{\le i} \beta \right)
} } - t^{s( \mbi{\mk I}^{>i},  \mk I^{(i)}) + {\color{black} 
\chi \left( 
\alpha- d^{\le i} \beta,
d^{(i)} \beta \right)
}
}) 
\int_{ [M^{\zeta^{-}}(\alpha - k \beta )]^{vir}} 
 \euk^{t}(\Lambda_{\adj}).
\end{align*}
\end{thm}

\section{Affine Laumon space of type $A_{N-1}^{(1)}$}

We consider the framed quiver $Q^{N} \colon$
\begin{center}
\includegraphics[scale=1]{chainsaw}
\end{center}
and the associated framed quiver moduli spaces
$M^{\zeta}(\mbi{r}, \mbi{v})$ called 
the {\it chainsaw quiver varieties}
of type $A^{(1)}_{N-1}$ where
$\mbi{r} = (r_{1}, \ldots, r_{N})$ and 
$\mbi{v}=(v_{1}, \ldots, v_{N})$.
When $N=1$, the framed quiver moduli spaces 
associated to $Q^{1} \colon$
\begin{center}
\includegraphics[scale=1]{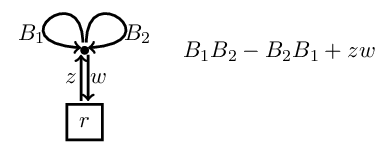}
\end{center}
are called the {\it Jordan quiver varieties}.

We consider the group 
$\Gamma_{N} = \lbrace \zeta \in \C^{\ast} \mid 
\zeta^{N}=1
\rbrace$ of $N$-th roots of the unit.

\subsection{Affine Laumon space as a component of $\Gamma_{N}$-fixed points set}

{\color{black}
In the following, we idenfity 
$I=Q_{0} \setminus \lbrace \infty 
\rbrace$ with 
$\Z/N \Z$, and slightly change the notation
from the previous section.
Here we introduce $W=\bigoplus_{i \in I} W_{i}$ with $W_{i}=\Hom(V_{\infty}, V_{i})^{\oplus r_{i}}$ corresponding to the framings to $i \in I$ in $Q^{N}$, and set
$V= \bigoplus_{i \in I} V_{i}$, 
omitting $V_{\infty}$, 
with $\dim V_{i}=v_{i}$ for $i \in I$.
Their shifts $W[l]$, $V[l]$ have the gradings given by
$W[l]_{i}=W_{i+l}$, $V[l]_{i}=V_{i+l}$.
They can be viewed as $\Gamma_{N}$-representations by regarding
gradings as weight 
space decompositions.
}

{\color{black}
We write by 
$\C_{\kappa/ \zeta}, \C_{q}$ and 
$\C_{\kappa q/ \zeta}$ the weight spaces
of the group $\Gamma_{N} \times \mb T$
with the weights $\kappa/\zeta$, $q$, and 
$\kappa q/ \zeta$ respectively,
where 
$\zeta \in \Gamma$, and
$\mb T = \C^{\ast}_{q} \times \C^{\ast}_{\kappa} \times 
\prod_{\alpha=1}^{r} 
\C^{\ast}_{e_{\alpha}} \times \C^{\ast}_{\mu_{\alpha}}
\times \C^{\ast}_{\nu_{\alpha}}$.
}
We set 
\begin{align}
\mb M = \mb M_{Q^{1}}
&= \mb M(W, V)
=
\End_{\C}(V) \otimes 
\left( \C_{\kappa/\zeta } \oplus \C_{q} \right) \times 
\Hom_{\C} (W, V) \times \Hom_{\C}(V, W) \otimes \C_{\kappa q/\zeta},
\\
\label{jordanml}
\mb L = \mb L_{Q^{1}}
&=
\mb L( V)
=
\Hom_{\C} (V, V) \otimes \C_{ \kappa q/\zeta}.
\end{align}
We regard the component 
$\prod_{\alpha=1}^{r} 
\C^{\ast}_{e_{\alpha}}$
as the diagonal torus of $\GL(W_{1}) \times \cdots \times \GL(W_{N})$
naturally acting on $\mb M$,
and the other component $\C^{\ast}_{\mu_{1}} \times \cdots \times 
\C^{\ast}_{\mu_{r}}
\times \C^{\ast}_{\nu_{1}} \times \cdots \times \C^{\ast}_{\nu_{r}}$
acts on $\mb M$ and $\mb L$ trivially.
Then the {\color{black} map} $\mu=\mu_{Q^{1}} \colon \mb M=\mb M_{Q^{1}} \to
\mb L=\mb L_{Q^{1}}$ sending
$(B_{1}, B_{2}, z,w )$ to $B_{1}B_{2} - B_{2} B_{1} +zw$ is 
$\mb T$-equivariant.
Thus we get the $\mb T$-action on the Jordan quiver varieties
$M^{\pm}(r, v)=[\mu^{-1} (0)^{\pm} / \GL(V)]$ 
where $r=r_{1} + \cdots + r_{N}$, $v=v_{1} + \cdots + v_{N}$, and
$\mu^{-1}(0)^{\pm}$ is the stable loci in $\mb \mu^{-1}(0)$
with respect to the stability parameters $\pm 1$.
Elements in $\mu^{-1}(0)$ are called
{\it ADHM data}.

To define $\Gamma_{N}$-action on $M^{-}(r,v)$, we forget 
the grading of $V$, and regard it as a just vector space $\C^{v}$.
Via the only $\Gamma_{N}$-action on $W$,
we have $\Gamma_{N} \times \GL(V)$-actions on $\mb M$, $\mb L$.
This induces $\Gamma_{N}$-action 
on $M^{-}(r,v)$ compatible with $\mb T$-action.
{\color{black} We take ADHM datum $(B,z,w)$
representing $\Gamma_{N}$-fixed point in 
$M^{-}(r, v)$.
For $\zeta \in \Gamma_{N}$, there exists
$g(\zeta) \in \GL(V)$ such that 
$(\zeta, g(\zeta)) \cdot (B, z,w)=   
(B,z,w)$.
Such an element $g(\zeta)$ is uniquely determined by the stability for $(B,z,w)$. 
This induces a group homomorphism
$\Gamma_{N} \to \GL(V)$
giving $I=\Z/N \Z$-grading on $V$.}
Taking the degree $0$ parts, we have
\begin{align}
\nonumber
\mb M_{Q^{N}}
&= 
\Hom_{\Gamma_{N}} (V, V[1]) \otimes \C_{\kappa} \times 
\End_{\Gamma_{N}}(V) \otimes \C_{q} \times 
\Hom_{\Gamma_{N}} (W, V) \times 
\Hom_{\Gamma_{N}}(V, W[1]) \otimes \C_{\kappa q},
\\
\label{alml}
\mb L_{Q^{N}}
&=
\Hom_{\Gamma_{N}} (V, V[1]) \otimes \C_{\kappa q}.
\end{align}
The {\color{black} map} $\mu \colon \mb M \to \mb L$ keeps 
these degree $0$ parts and we have
$\mu|_{\mb M_{Q^{N}}}=\mu_{Q^{N}}$.
For the stability parameters
$\zeta^{\al}=(- 1)_{i \in I} , \zeta^{\cal}= (1)_{i \in I} 
\in \mb R^{I}$, 
the associated framed quiver moduli spaces
\[
M^{\al}(\mbi{r}, \mbi{v})=\left[\mu_{Q^{N}}^{-1}(0)^{\zeta^{\al}} \Big / 
\prod_{i \in I} \GL(V_{i})\right], \quad
M^{\cal}(\mbi{r}, \mbi{v})=\left[\mu_{Q^{N}}^{-1}(0)^{\zeta^{\cal}} \Big / 
\prod_{i \in I} \GL(V_{i})\right]
\] 
are called 
the {\it affine Laumon space} of type $A^{(1)}_{N-1}$, 
and the {\it dual}.
Here we notice that $\mu^{-1}(0)^{\zeta^{\al}}=
\mu^{-1}(0)^{-} \cap \mb M_{Q^{N}}$ and 
$\mu^{-1}(0)^{\zeta^{\cal}}=
\mu^{-1}(0)^{+} \cap \mb M_{Q^{N}}$.

{\color{black}
We consider the $\Gamma_{N}$-fixed points set 
$M^{-}(r, v)^{\Gamma_{N}}$.
By the construction, we have 
$\mb T$-equivariant embeddings 
$M^{\al}(\mbi{r},\mbi{v}) \to M^{-}(r, v)
^{\Gamma_{N}}$}.
By \cite{N0}, all the connected components of
$M^{-}(r, v)^{\Gamma_{N}}$ are of this form taking all
possible $\mbi{v}$ such that $v_{1} + \cdots + v_{N}=v$.
In particular, the $\mb T$-fixed points sets 
$M^{\al}(\mbi{r}, \mbi{v})^{\mb T}$ are
regarded as a subset of 
{\color{black}
$M^{-}(r, v)^{\mb T \times \Gamma_{N}}$.
In fact, 
we see that 
\[
M^{-}(r, v)^{\mb T}
=
\bigsqcup_{v_{1} + \cdots + v_{N}=v} 
M^{\al}(\mbi{r}, \mbi{v})^{\mb T},
\]
and consequently $M^{-}(r, v)^{\mb T} = 
M^{-}(r,v)^{\mb T \times \Gamma_{N}}$.} 
This follows from combinatorial descriptions of
$M^{-}(r, v)^{\color{black} \mb T}$ 
summarized in the next subsection.

\subsection{Combinatorial description 
of $M^{\al}(\mbi{r}, \mbi{v})^{\color{black} \mb T}$}
We recall combinatorial descriptions of $M^{-}(r, v)^{\mb T}$,
the tautological bundle, and the 
tangent bundles over the fixed points from \cite{NY1}.

We can identify the all fixed points in 
$M^{-}(r, v)^{\mb T}$ with
the collections $\vec{Y}=
(Y^{(1)}, \ldots, Y^{(N)})$ of Young diagrams
$Y^{(1)}, \ldots, Y^{(N)}$ such that 
$|Y^{(1)}| + \cdots + |Y^{(N)}|=v$, where
$|Y|$ denotes the number of boxes in $Y$ 
for a Young diagram $Y$.
By the gauge transformation, we have the $\mb T$-weight
decompositions 
\begin{align}
\label{fixedpt-}
\mc V \Big|_{\vec{Y}}
&= \bigoplus_{\alpha=1}^{r} 
\bigoplus_{(l, m) \in Y_{\alpha}} 
\C_{e_{\alpha} \kappa^{-l+1} q^{-m+1} }, 
&
\mc W \Big|_{\vec{Y}}
&=
\bigoplus_{\alpha =1}^{r} \C_{e_{\alpha}}
\end{align}
of tautological bundles
$\mc V = [\mu^{-1}(0)^{-} \times V / \GL(V) ]$
and $\mc W =
M^{-}(r,v) \times W$.
We define $B_{1}, B_{2} \colon V \to V$, 
$\colon W \to V$, and $w=0 \colon V \to W$ so that 
\begin{align}
\label{b1action-}
B_{1} (\C_{e_{\alpha} \kappa^{-l+1} q^{-m+1} } )
&=\begin{cases}
\C_{e_{\alpha} \kappa^{-l} q^{-m+1} } & (l+1, m) \in Y_{\alpha} \\
0 & (l+1, m) \not \in Y_{\alpha}, \\
\end{cases}
\\
\label{b2action-}
B_{2} (\C_{e_{\alpha} \kappa^{-l+1} q^{-m+1} } )
&=
\begin{cases}
\C_{e_{\alpha} \kappa^{-l+1} q^{-m} } & (l, m+1) \in Y_{\alpha} \\
0 & (l, m+1) \not \in Y_{\alpha}, \\
\end{cases}
\\
\label{zwaction-}
z &= \bigoplus_{\alpha=1}^{r} \id_{\C_{e_{\alpha}}},
& 
w=0.
\end{align}
Then we have $(B_{1}, B_{2}, z, w) \in \mu^{-1}(0)^{-}$ and 
conjugations by the weight decompositions \eqref{fixedpt-} gives
gauge transformations from  the $\mb T$-orbits of 
$(B_{1}, B_{2}, z, w)$ to itself.
Thus the equivalence class $[B_{1}, B_{2}, z,w] \in M^{-}(r,v)$
is the $\mb T$-fixed point corresponding to $\vec{Y}$.
 
We identify a Young diagram $Y$ with the partition
$\lambda=(\lambda_{1}, 
\lambda_{2}, \ldots,)$ where $\lambda_{i}$ is the number
of boxes in $i$-th row in $Y$, and set $Y=Y_{\lambda}$.
For pairs $\lambda, \mu$ of partitions, 
we consider the {\it Nekrasov factor}
\begin{align}
\nonumber
\mathsf{N}_{\lambda, \mu}^{(k|N)}(u)
&=\mathsf{N}_{\lambda, \mu}^{(k|N)}(u|q, \kappa)
=\mathsf{N}_{Y_{\lambda}, Y_{\mu}}^{(k|N)}(u)
=\mathsf{N}_{Y_{\lambda}, Y_{\mu}}^{(k|N)}(u|q, \kappa)
\\
\nonumber
&=\prod_{1 \le l \le l' \atop -l+l' \equiv k \mod N} 
( u \kappa^{-l+l'} q^{-\mu_{l} +\lambda_{l'+1}}; q)
_{\lambda_{l'}-\lambda_{l'+1}}
\prod_{1 \le l \le l'  \atop l-l' -1\equiv k \mod N} 
( u \kappa^{l-l'-1} q^{\lambda_{l} -\mu_{l'}};q )
_{\mu_{l'}-\mu_{l'+1}}
\\
\label{nekfactor}
&=
\prod_{s \in Y_{\lambda} \atop L_{Y_{\lambda}}(s) \equiv k \mod N}
( 1-u \kappa^{L_{Y_{\lambda}}(s)} q^{-A_{Y_{\mu}}(s) -1})
\prod_{t \in Y_{\mu} \atop 
-L_{Y_{\mu}}(t) -1 \equiv k \mod N}
( 1-u \kappa^{-L_{Y_{\mu}}(t) -1} q^{A_{Y_{\lambda}}(t) }).
\end{align}
We often omit $N$ and write 
the Nekrasov factor \eqref{nekfactor} 
by $\mathsf{N}_{\lambda, \mu}^{(k)}(u)$, and so on.
In particular, we have
\begin{align}
\label{nekempty}
\mathsf{N}_{\lambda, \emptyset}^{(k)}(u)
&=
\prod_{(l,m) \in Y_{\lambda} \atop l-1 \equiv k \mod N} 
(1- u \kappa^{l-1} q^{m-1}),
&
\mathsf{N}_{\emptyset, \lambda}^{(k)}(u)
&=
\prod_{(l,m) \in Y_{\lambda} \atop -l \equiv k \mod N} 
(1- u \kappa^{-l} q^{-m}),
\end{align}
where $\emptyset$ denotes the empty partition.
By comparing \eqref{nekempty} 
with $\mc V \Big|_{\vec{Y}}$ in \eqref{fixedpt-},
we have
$\wedge_{-1/\kappa q} \mc V \Big|_{\vec{Y}} = \prod_{\alpha=1}^{r} 
\mathsf{N}_{\emptyset, Y_{\alpha}}^{(0|1)}(e_{\alpha})$
and 
$\wedge_{-1} \mc V^{\vee} \Big|_{\vec{Y}} = \prod_{\alpha=1}^{r} 
\mathsf{N}_{Y_{\alpha}, \emptyset}^{(0|1)}(e_{\alpha}^{-1})$.
By \cite[Theorem 2.11]{NY1}, we have
\begin{align}
\label{tangent-}
\wedge_{-1}T_{\vec{Y}}^{\ast} M^{-}(r,n) 
&=
\prod_{\alpha, \beta=1}^{r} 
\mathsf{N}_{Y_{\alpha}, Y_{\beta}}^{(0|1)} 
(e_{\beta}/ e_{\alpha} | q, \kappa).
\end{align}
 
For $\Gamma_{N}$-action, we rewrite 
$\mbi e=(e_{1}, \ldots, e_{r})$, 
$\mbi \mu=(\mu_{1}, \ldots, \mu_{r})$, and
$\mbi \nu=(\nu_{1}, \ldots, \nu_{r})$ by 
\begin{align}
\nonumber
\mbi e&=(\mbi e_{i})_{i\in I},
&
\mbi \mu&=(\mbi \mu_{i})_{i \in I},
&
\mbi \nu&= (\mbi \nu_{i})_{i \in I} 
&
\end{align} 
with 
$\mbi e_{i}=(e_{(i, \alpha)})_{\alpha=1}^{r_{i}}$
so that
$W_{i} = \bigoplus_{\alpha=1}^{r_{i}} \C_{e_{(i, \alpha)}}$, 
$\mbi \mu_{i}=(\mu_{(i, \alpha)})_{\alpha=1}^{r_{i}}$, and
$\mbi \nu_{i}=(\nu_{(i, \alpha)})_{\alpha=1}^{r_{i}}$.
Set
\begin{align}
\label{degree}
\deg \kappa=-1, 
\quad
\deg q =0,
\quad
\deg e_{(i, \alpha)}=i,
\end{align}
and
$\deg \mu_{(i,\alpha)}= \deg \nu_{(i, \alpha)}=0 \in \Z/ N \Z$ 
for $i \in I$, 
and $1 \le \alpha \le r_{i}$.
We define a direct summand $V_{j}$ of $V$ by taking 
the direct sum of components of $\mc V|_{\vec{Y}}$ in \eqref{fixedpt-} 
with the weights of degree $j \in \Z/N\Z$. 
Explicitly, we have
\begin{align}
\label{altaut-}
V_{j}
&=
\bigoplus_{i=1}^{N} \bigoplus_{\alpha=1}^{r_{i}} \bigoplus
_{\substack{(l, m ) \in Y_{(i,\alpha)} 
\\ i+l-1 \equiv j \mod N}} 
\C_{e_{(i,\alpha)} \kappa^{-l+1} q^{- m +1}}. 
\end{align}

Then we have $(B_{1}, B_{2}, z, w ) \in \mb M_{Q^{N}}$, and
$[B_{1}, B_{2}, z, w] \in M^{\al} (\mbi{r}, \mbi{v})^{\mb T}$.
Thus we get the proposition.
\begin{prop}
\label{prop:comb-}
(a) We have a bijection from the fixed points set 
$M^{\al}(\mbi{r}, \mbi{v})^{\mb T}$ to the set of collections
$\vec{Y}=(Y_{(i, \alpha)})_{i, \alpha}$ of Young diagrams
$Y_{(i, \alpha)}$ indexed by $i \in I$, and 
$1 \le \alpha \le r_{i}$ such that
\[
\sum_{i=1}^{N} \sum_{\alpha=1}^{r_{i}}
\left| \left\lbrace \left.
(l, m) \in Y_{(i,\alpha)} \ \right| \ 
i + l -1 \equiv j \mod N
\right\rbrace\right| 
= v_{j}.
\]
{\color{black}The corresponding $Q^{N}$-representation at each fixed point $\vec{Y}$ is 
described by  
\eqref{fixedpt-}, \eqref{b1action-}, \eqref{b2action-}, and \eqref{zwaction-}. }
\\
(b) For the fixed point $\vec{Y} \in M^{\al}
(\mbi{r}, \mbi{v})^{\mb T}$
via the bijection in (a), we have
\begin{align}
\label{taut-nek}
\wedge_{-\kappa q} \mc V_{i} \Big|_{\vec{Y}}
&=\prod_{j=1}^{N} \prod_{\beta=1}^{r_{j}}
\mathsf{N}_{\emptyset, Y_{(j,\beta)}}^{(j-i-1)}
(e_{(j,\beta)}),
&
\wedge_{-1} \mc V_{j}^{\vee} \Big|_{\vec{Y}} 
&=\prod_{i=1}^{N} \prod_{\alpha=1}^{r_{i}}
\mathsf{N}_{Y_{(i,\alpha)}, \emptyset}^{(j-i)}
(e_{(i,\alpha)}^{-1}).
\end{align}
\\
(c) For the fixed point $\vec{Y} \in 
M^{\al}(\mbi{r}, \mbi{v})^{\mb T}$
via the bijection in (a), we have
\begin{align}
\label{tangent-nek}
\wedge_{-1} T^{\ast}_{\vec{Y}} M^{\al}(\mbi{r}, \mbi{v})
=
\prod_{i,j \in I}
\prod_{\alpha=1}^{r_{i}} \prod_{\beta=1}^{r_{j}} 
\mathsf{N}_{Y_{(i,\alpha)}, Y_{(j,\beta)}}^{(j-i|N)} 
( e_{(j,\beta)}/e_{(i,\alpha)} | q, \kappa).
\end{align}
\end{prop} 
\proof
(a) follows from the preceding arguments.
(b) follows by comparing \eqref{altaut-} and \eqref{nekempty}.
For (c), notice that the tangent space 
$T_{\vec{Y}} M^{\al} (\mbi{r}, \mbi{v})$ 
at the fixed point $\vec{Y}=(Y_{(i, \alpha)})$ 
is isomorphic to the degree $0$
part of $T_{\vec{Y}} M^{-}(r,v)$.
We compute the weight space decomposition of 
$T_{\vec{Y}} M^{-}(r, v)= \mb M - \mb L - \End_{\C}(V)$
in the representation ring $K_{\mb T}(\pt)$, where
$V$ is identified with $\mc V|_{\vec{Y}}$ in \eqref{fixedpt-}.
Then $\Gamma_{N}$-action defined in \eqref{jordanml} is compatible 
with the degree \eqref{degree}.
Hence we have the assertion. 
\endproof

\subsection{Relations between $M^{\al}(\mbi{r}, \mbi{v})$ and $M^{\cal}(\mbi{r}, \mbi{v})$}

For combinatorial descriptions of $M^{\cal}(\mbi{r},\mbi{v})$,
we consider a map $\mb M(W, V) \to \mb M(W^{\vee}, V^{\vee})$
sending $(B_{1}, B_{2}, z, w)$ to 
$(- B_{1}^{\vee}, B_{2}^{\vee}, w^{\vee}, z^{\vee} )$ where
${}^{\vee}$ denotes taking the transpose
$ f^{\vee} \colon U_{2}^{\vee} \to U_{1}^{\vee}$ for a 
linear map $f \colon U_{1} \to U_{2}$. 
By the remark after Definition \ref{def:stab}, this map
interchanges stability conditions 
with respect to the stability parameters
$\pm 1 \in \mb R$.
Hence it induces a map $\mu^{-1}(0)^{+} \to \mu^{-1}(0)^{-}$.
To make this map 
$\GL(V) \times \mb T$-equivariant,
we take a group homomorphism 
\begin{align}
\label{toruseq}
\GL(V) \times \mb T \to \GL(V^{\vee}) \times \mb T, \quad
(g, \mbi e, \mbi \mu, \mbi \nu, \kappa, q)
\mapsto
( \kappa q \  (g^{\vee})^{-1}, (\mbi e^{\vee})^{-1}, 
(\mbi \nu^{\vee})^{ -1}, (\mbi \mu^{\vee})^{ -1}, \kappa ,q), 
\end{align}
where we identify the diagonal 
torus in $\GL(W)$ with the one in $\GL(W^{\vee})$ by taking
the dual base.

For a $I=\Z/N\Z$-graded vector space $U$, 
we define the grading of the dual space $U^{\vee}$
by $(U^{\vee})_{i}=(U_{-i})^{\vee}$.
Then the above map induces an isomorphism 
$\Phi \colon M^{+}(r, v) \to M^{-}(r, v)$ respecting
sub-manifolds $M^{\cal}(\mbi{r}, \mbi{v})$ and 
$M^{\al}(\mbi{r}^{\vee}, \mbi{v}^{\vee}[1])$ such that 
$\Phi^{\ast} \mc V_{j} = \mc V_{-j-1}^{\vee} \otimes \C_{\kappa q}$
for $j \in I$, and $\mb T$-equivariant after taking
a homomorphism 
\begin{align}
\label{duality}
\mb T \to \mb T, \quad (\mbi e, \mbi \mu, \mbi \nu)
\mapsto ((\mbi e^{\vee})^{-1}, (\mbi \nu^{\vee})^{-1}, 
(\mbi \mu^{\vee})^{-1}). 
\end{align}
Here we set $(\mbi{r}^{\vee})_{i} = r_{-i}$, 
$(\mbi{v}^{\vee}[1])_{i}=(\mbi{v}^{\vee})_{i+1}=v_{-i-1}$, and
$((\mbi e^{\vee})^{-1} )_{i} = ((\mbi e_{-i}^{\vee})^{-1} )
\in GL(W_{-i}^{\vee})$ for $i \in I=\Z/N\Z$.
Hence we have 
the combinatorial descriptions
similar to the previous sub-section.

We identify $M^{+}(r, v)^{\mb T}$ with the set of 
the collections $\vec{Y}=
(Y^{(1)}, \ldots, Y^{(N)})$ of Young diagrams
$Y^{(1)}, \ldots, Y^{(N)}$ such that 
$|Y^{(1)}| + \cdots + |Y^{(N)}|=v$.
We have the $\mb T$-weight decompositions 
\begin{align}
\label{fixedpt+}
\mc V \Big|_{\vec{Y}}
&= \bigoplus_{\alpha=1}^{r} 
\bigoplus_{(l, m) \in Y_{\alpha}} 
\C_{e_{\alpha} \kappa^{l} q^{m} }, 
&
\mc W \Big|_{\vec{Y}}
&=
\bigoplus_{\alpha =1}^{r} \C_{e_{\alpha}}
\end{align}
of tautological bundles
$\mc V = [\mu^{-1}(0)^{+} \times V / \GL(V) ]$
and $\mc W =
M^{+}(r,v) \times W$.
We define $B_{1}, B_{2} \colon V \to V$, 
$z \colon W \to V$, and $w \colon V \to W$
so that 
\begin{align}
\label{b1action+}
B_{1} (\C_{e_{\alpha} \kappa^{l} q^{m} } )
&=\begin{cases}
\C_{e_{\alpha} \kappa^{l-1} q^{m} } & (l-1, m) \in Y_{\alpha} \\
0 & (l-1, m) \not \in Y_{\alpha}, \\
\end{cases}
\\
\label{b2action+}
B_{2} (\C_{e_{\alpha} \kappa^{l} q^{m} } )
&=
\begin{cases}
\C_{e_{\alpha} \kappa^{l} q^{m-1} } & (l, m-1) \in Y_{\alpha} \\
0 & (l, m-1) \not \in Y_{\alpha}, \\
\end{cases}
\\
\label{zwaction+}
z &= 0, &
w(e_{\alpha} \kappa^{l} q^{m} )=
\begin{cases} 
e_{\alpha} & l=m=1 \\
0 & \text{otherwise}.
\end{cases}
\end{align}
Then the equivalence class $[B_{1}, B_{2}, z,w] \in M^{+}(r,v)$
is the $\mb T$-fixed point corresponding to $\vec{Y}$.

Comparing \eqref{nekempty} 
with $\mc V \Big|_{\vec{Y}}$ in \eqref{fixedpt+},
we have
$\wedge_{-1/\kappa q} \mc V \Big|_{\vec{Y}} = \prod_{\alpha=1}^{r} 
\mathsf{N}_{\emptyset, Y_{\alpha}}^{(0|1)}(e_{\alpha})$
and 
$\wedge_{-1} \mc V^{\vee} \Big|_{\vec{Y}} = \prod_{\alpha=1}^{r} 
\mathsf{N}_{Y_{\alpha}, \emptyset}^{(0|1)}(e_{\alpha}^{-1})$.
By applying $\Phi^{\ast}$ to \eqref{tangent-} and 
taking the homomorphism \eqref{duality} of $\mb T$, we have
\begin{align}
\label{tangent+}
\wedge_{-1}T_{\vec{Y}}^{\ast} M^{+}(r,n) 
&=
\prod_{\alpha, \beta=1}^{r} 
\mathsf{N}_{Y_{\alpha}, Y_{\beta}}^{(0|1)} 
( e_{\alpha} / e_{\beta} | q, \kappa).
\end{align}
On $M^{\cal}(\mbi{r}, \mbi{v})$, we have
\begin{align}
\label{altaut+}
V_{j}
&=
\bigoplus_{i=1}^{N} \bigoplus_{\alpha=1}^{r_{i}} \bigoplus
_{\substack{(l, m ) \in Y_{(i,\alpha)} 
\\ i-l \equiv j \mod N}} 
\C_{e_{(i,\alpha)} \kappa^{l} q^{m }},
\end{align}
and the similar proposition 
for $M^{\cal}(\mbi{r}, \mbi{v})$ to Proposition \ref{prop:comb-}.

\begin{prop}
\label{prop:comb+}
(a) We have a bijection from the fixed point set 
$M^{\cal}(\mbi{r}, \mbi{v})$ to the set of collections
$\vec{Y}=(Y_{(i, \alpha)})_{i, \alpha}$ of Young diagrams
$Y_{(i, \alpha)}$ indexed by $i=1, \ldots, N$, and 
$1 \le \alpha \le r_{i}$ such that
\[
\sum_{i=1}^{N} \sum_{\alpha=1}^{r_{i}}
\left| \left\lbrace 
(l, m) \in Y_{(i,\alpha)} \mid
i - l \equiv j \mod N
\right\rbrace\right| 
= v_{j}.
\]
{\color{black}The corresponding $Q^{N}$-representation at each fixed point $\vec{Y}$ is 
	described by  
	\eqref{fixedpt+}, \eqref{b1action+}, \eqref{b2action+}, and \eqref{zwaction+}. }
\\
(b) For the fixed point $\vec{Y} \in M^{\cal}(\mbi{r}, \mbi{v})$
via the bijection in (a), 
we have
\begin{align}
\label{taut+nek}
\wedge_{-1} \mc V_{j}^{\vee} 
&=\prod_{i=1}^{N} \prod_{\alpha=1}^{r_{i}}
\mathsf{N}_{\emptyset, Y_{(i,\alpha)}}^{(j-i)}
(e_{(i,\alpha)}^{-1}),
&
\wedge_{-\kappa q} \mc V_{j} 
&=\prod_{i=1}^{N} \prod_{\alpha=1}^{r_{i}}
\mathsf{N}_{Y_{(i,\alpha), \emptyset }}^{(i-j-1)}
(e_{(i,\alpha)}).
\end{align}
\\
(c) For the fixed point $\vec{Y} \in M^{\cal}(\mbi{r}, \mbi{v})$
via the bijection in (a), 
we have
\begin{align}
\label{tangent-nek}
\wedge_{-1} T^{\ast}_{\vec{Y}} M^{\cal}(\mbi{r}, \mbi{v})
=
\prod_{i,j \in I}
\prod_{\alpha=1}^{r_{i}}  \prod_{\beta=1}^{r_{j}} 
\mathsf{N}_{Y_{(i,\alpha)}, Y_{(j,\beta)}}^{(i-j|N)} 
(e_{(i,\alpha)} / e_{(j,\beta) } | q, \kappa).
\end{align}
The tangent space $T_{\vec{Y}} M^{\cal} (\mbi{r}, \mbi{v})$ 
at the fixed point $\vec{Y}=(Y_{(i, \alpha)})$ 
via the bijection in (a) is isomorphic to the degree $0$
part of $T_{\vec{Y}} M^{+}(r,v)$.
\end{prop} 
\proof
(a) follows from \eqref{altaut+}.
By comparing it with \eqref{nekempty}, we have (b).
By taking degree 0 part in \eqref{tangent+}, we have (c).
\endproof

From the combinatorial point of view, 
we introduce stable and co-stable collections
$\mbi Y^{\pm}_{\mbi \lambda}=
\left(
(Y^{\pm}_{i,\alpha})_{\alpha=1}^{r_{i}}
\right)_{i \in I}$ of colored Young 
diagrams associated to
a collection 
$\mbi \lambda=
\left(
( \lambda^{(i,\alpha)})_{\alpha=1}^{r_{i}}
\right)_{i \in I}$
of partitions. 
We define 
stable collections $\mbi Y^{-}_{\mbi \lambda}=
\left(
(Y^{-}_{i,\alpha})_{\alpha=1}^{r_{i}}
\right)_{i \in I}$
of colored Young 
diagrams 
by labeling boxes in $k$-th row a color $i+k-1 \in I$
in $i$-th Young diagrams $(Y_{i,1}, \ldots, Y_{i,r_{i}})$.
On the other hand, we define 
co-stable collections $\mbi Y^{+}_{\mbi \lambda}=
\left(
(Y^{+}_{i,\alpha})_{\alpha=1}^{r_{i}}
\right)_{i \in I}$
of colored Young 
diagrams 
by labeling boxes in $k$-th row a color $i-k \in I$
in $i$-th Young diagrams $(Y_{i,1}, \ldots, Y_{i,r_{i}})$:
\begin{center}
\includegraphics[scale=1]{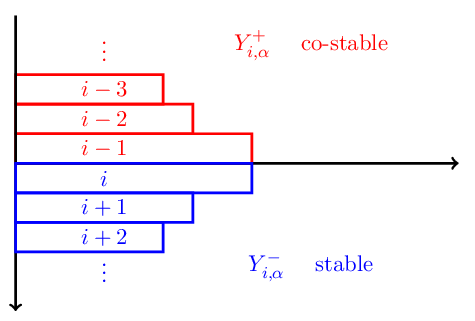}
\end{center}
We set $Y^{\pm}_{\mbi \lambda} =\bigsqcup_{i \in I} 
\bigsqcup_{\alpha=1}^{r_{i}} Y_{i, \alpha}$.
Then $M^{\al}(\mbi w, \mbi v)^{\mb T}$ is parametrized by 
the collections $\mbi \lambda=
\left(
( \lambda^{(i,\alpha)})_{\alpha=1}^{r_{i}}
\right)_{i \in I}$ of partitions such that
\[
\sharp \lbrace s \in Y^{-}_{\mbi \lambda} \mid
s \text{ has a color } i \in I \rbrace
=
v_{i}
\]
for all $i \in I$.
$M^{\cal}(\mbi w, \mbi v)^{\mb T}$ is parametrized by 
the collections $\mbi \lambda=
\left(
( \lambda^{(i,\alpha)})_{\alpha=1}^{r_{i}}
\right)_{i \in I}$ of partitions such that
\[
\sharp \lbrace s \in Y^{+}_{\mbi \lambda} \mid
s \text{ has a color } i \in I \rbrace
=
v_{i}
\]
for all $i \in I$.

\subsection{Non-stationary Ruijsenaars function}

We define the $K$-theoretic class
$\Lambda_{\adj}$ by
\begin{align}
\nonumber
\Lambda_{\adj} 
&=\sum_{j \in I} 
\sum_{\beta=1}^{r_{j}} 
\mc V_{j} \otimes \C_{1/e_{(j, \beta)}} +
\sum_{i \in I} 
\sum_{\alpha=1}^{r_{i+1} } 
\mc V_{i}^{\vee} \otimes \C_{e_{(i+1, \alpha)} \kappa q }
\\
\label{ladj}
&+ \Hom_{\Gamma_{N}}(V, V[1])\otimes (\C_{\kappa} - \C_{\kappa q})
+ \End_{\Gamma_{N}}(V)\otimes (\C_{q} - \C),
\end{align}
using tautological bundles on moduli stacks.
By the construction using \eqref{alml},
we have $\Lambda_{\adj} = T M^{\zeta} (\mbi{r}, \mbi{v})$ 
in $K_{\mb T}(M^{\zeta}(\mbi{r}, \mbi{v}))$.
\begin{prop}
\label{prop:adjnek}
(a)
For $\vec{Y}=((Y_{(i, \alpha)})_{\alpha=1}^{r_{i}} )_{i \in I} 
\in M^{-}(\mbi{r}, \mbi{v})^{\mb T}$, 
we have
\[
\eukt(\Lambda_{\adj})|_{\vec{Y}}
=
\prod_{i,j \in I}
\prod_{\alpha=1}^{r_{i}} \prod_{\beta=1}^{r_{j}}
\mathsf{N}^{(j-i|N)}_{\lambda^{(i, \alpha)},\lambda^{(j, \beta)}} 
(t e_{(j, \beta)}/e_{(i,\alpha)} |q,\kappa).
\]
(b)
For $\vec{Y}=((Y_{(i, \alpha)})_{\alpha=1}^{r_{i}} )_{i \in I} 
\in M^{+}(\mbi{r}, \mbi{v})^{\mb T}$, 
we have
\[
\eukt(\Lambda_{\adj})|_{\vec{Y}}
=
\prod_{i,j \in I}
\prod_{\alpha=1}^{r_{i}} \prod_{\beta=1}^{r_{j}}
\mathsf{N}^{(i-j|N)}_{\lambda^{(i, \alpha)},\lambda^{(j, \beta)}} 
(t e_{(i,\alpha)} /e_{(j, \beta)} |q,\kappa).
\]
\end{prop}

For variables $\mbi p=(p_{i})_{i \in I} =(p_{1}, \ldots, p_{N})$,
set
\[
Z^{\mbi{r}}_{\adj}(\mbi e, t|\mbi p|q, \kappa)
=
\sum_{\mbi{v}\in \Z_{\ge 0}^{I} }
\mbi p^{\mbi{v}} \int_{M^{\al}(\mbi{r}, \mbi{v}) } \Lambda_{\adj},
\quad
\widecheck{Z}^{\mbi{r}}_{\adj}(\mbi e, t|\mbi p|q, \kappa)
=
\sum_{\mbi{v}\in \Z_{\ge 0}^{I} }
\mbi p^{\mbi{v}} \int_{M^{\cal}(\mbi{r}, \mbi{v}) } \Lambda_{\adj}.
\]

\begin{prop}
\label{adj:nek2}
(a) We have
\begin{align}
\nonumber
Z^{\mbi{r}}_{\adj}(\mbi e, t|\mbi p|q, \kappa)
&=
\sum_{\mbi{\lambda} \in \prod_{i \in I} \mathsf{P}^{r_{i}}}
\prod_{i,j \in I}
\prod_{\alpha=1}^{r_{i}} \prod_{\beta=1}^{r_{j}}
{\mathsf{N}^{(j-i|N)}_{\lambda^{(i, \alpha)},\lambda^{(j, \beta)}} 
(t e_{(j, \beta)}/e_{(i,\alpha)} |q,\kappa) \over 
\mathsf{N}^{(j-i|N)}_{\lambda^{(i, \alpha)},\lambda^{(j, \beta)}} 
(e_{(j, \beta)}/e_{(i,\alpha)} |q,\kappa)}
\prod_{i \in I} \prod_{j \ge 1} 
p_{i+j-1}^{\sum_{\alpha=1}^{r_{i}} 
\lambda^{(i,\alpha)}_{j}},
\end{align}
where $\mbi{\lambda}=
\left(
(\lambda^{(i, \alpha)})_{\alpha=1}^{r_{i}}
\right)_{i \in I}$ runs in 
$\prod_{i \in I} \mathsf{P}^{r_{i}}$ for
the set $\mathsf{P}$ of partitions.
\\
(b) We have 
\begin{align}
\nonumber
\widecheck{Z}^{\mbi{r}}_{\adj}(\mbi e, t|\mbi p|q, \kappa)
&=
\sum_{\mbi{\lambda} \in \prod_{i \in I} \mathsf{P}^{r_{i}}}
\prod_{i,j \in I}
\prod_{\alpha=1}^{r_{i}} \prod_{\beta=1}^{r_{j}}
{\mathsf{N}^{(i-j|N)}_{\lambda^{(i, \alpha)},\lambda^{(j, \beta)}} 
(t e_{(i,\alpha)} /e_{(j, \beta)} |q,\kappa) \over 
\mathsf{N}^{(i-j|N)}_{\lambda^{(i, \alpha)},\lambda^{(j, \beta)}} 
(e_{(i,\alpha)} /e_{(j, \beta)}|q,\kappa)}
\prod_{i \in I} \prod_{j \ge 1} 
p_{i-j}^{\sum_{\alpha=1}^{r_{i}} 
\lambda^{(i, \alpha)}_{j}}
\\
&=
\nonumber
Z^{\mbi{r}^{\vee}}_{\adj}
((\mbi e^{\vee})^{-1}, t|(\mbi p[-1])^{\vee}| q, \kappa).
\end{align}
\end{prop}
\proof
These follows from Proposition \ref{prop:comb-} and 
Proposition \ref{prop:comb+}.
We can also check the second equality in (b) from 
\begin{align}
\nonumber
\Lambda_{\adj} 
&=
\Phi^{\ast} \left. T^{\ast} M^{\al} (\mbi{r}^{\vee}, 
\mbi{v}^{\vee}[1])
\right|_{\substack{\mbi e = (\mbi e^{\vee})^{-1}}}.
\end{align}
In fact, the fixed points set $M^{\al}(\mbi r^{\vee}, 
\mbi v^{\vee}[1])^{\mbi T}$
is parametrized by $\mbi \lambda^{\vee}=
(\lambda^{(-i, \alpha)})_{\alpha=1}^{r_{i}})_{i \in I}$
for $\mbi \lambda \in \prod_{i \in I} P^{r_{i}}$, and we have
\begin{align*}
&
\prod_{i,j \in I}
\prod_{\alpha=1}^{r_{i}} \prod_{\beta=1}^{r_{j}}
{\mathsf{N}^{(j-i|N)}_{\lambda^{(-i, \alpha)},\lambda^{(-j, \beta)}} 
(t e_{(-j, \beta)}^{-1} /e_{(-i,\alpha)}^{-1} |q,\kappa) \over 
\mathsf{N}^{(j-i|N)}_{\lambda^{(-i, \alpha)},\lambda^{(-j, \beta)}} 
(e_{(-j, \beta)}^{-1} /e_{(-i,\alpha)}^{-1} |q,\kappa)}
\prod_{i \in I} \prod_{j \ge 1} \prod_{\alpha=1}^{r_{i}}
(\mbi p[-1]^{\vee})^{\lambda^{(-i, \alpha)}_{j}}_{i+j-1}
\\
&=
\prod_{i,j \in I}
\prod_{\alpha=1}^{r_{i}} \prod_{\beta=1}^{r_{j}}
{\mathsf{N}^{(-j+i|N)}_{\lambda^{(i, \alpha)},\lambda^{(j, \beta)}} 
(t e_{(i,\alpha)}/e_{(j, \beta)} |q,\kappa) \over 
\mathsf{N}^{(-j+i|N)}_{\lambda^{(i, \alpha)},\lambda^{(j, \beta)}} 
(e_{(i,\alpha)} /e_{(j, \beta)} |q,\kappa)}
\prod_{i \in I} \prod_{j \ge 1} \prod_{\alpha=1}^{r_{i}}
(\mbi p[-1]^{\vee})^{\lambda^{(i, \alpha)}_{j}}_{-i+j-1}.
\end{align*}
\endproof

When $r_{1} = \cdots = r_{N}=d$ for a positive integer 
$d$, 
following \cite{S}, we introduce
the {\it non-stationary Ruijsenaars function} 
\begin{align}
\nonumber
f^{\hat{\mk{gl}}_{N}}( \mbi x, p |\mbi s,\kappa| q,t )
&=
\sum_{\mbi{\lambda} \in \prod_{i \in I} 
\mathsf{P}^{d}}
\prod_{i,j \in I}
\prod_{\alpha, \beta=1}^{d} 
{
\Nk^{(j-i|N)}_{\lambda^{(i, \alpha)}, 
\lambda^{(j, \beta)}}
(t s_{(j, \beta)} / s_{(i, \alpha)} ) 
\over 
\mathsf{N}_{\lambda^{(i, \alpha)}, 
\lambda^{(j, \beta)}}^{(j-i|N)} 
( s_{(j, \beta)} / s_{(i, \alpha)} )}
\prod_{i \in I} \prod_{j \ge 1} \prod_{\alpha=1}^{d}
(p x_{i + j}/ t^{d} x_{i + j -1})
^{\lambda^{(i, \alpha)}_{j}},
\end{align}
where $\mbi x=(x_{i})_{i \in I}$ are variables, 
$\mbi s=((s_{(i,\alpha)})_{\alpha=1}^{d})_{i \in I}$, 
$p, q, \kappa$, 
and $t$
are parameters.
This is slightly generalized notion, which is more 
suitable here, without assuming
$r_{1}=\cdots =r_{N}=1$.
It is easy to see
\begin{align}
\label{adjruijs}
Z^{\mbi{r}}_{\adj}(\mbi s, t | 
(px_{i+1} /t^{d} x_{i})_{i\in I}| q, \kappa)
=f^{\hat{\mk{gl}}_{N}}(\mbi x, p | \mbi s, \kappa 
| q, t).
\end{align}
For a positive integer $d$, set
\begin{align*}
f_{N, \infty}(X|S|q,t)
&=
\sum_{\theta \in (\widehat{M}_{N})^{d}}
\prod_{\alpha, \beta=1}^{d}
c_{N, \infty}(\theta^{(\alpha)}, \theta^{(\beta)} 
| S_{\alpha}, S_{\beta}  | q, t) 
\prod_{1 \le i < k \le N} 
(X_{k}/X_{i})^{\sum_{\alpha=1}^{d} 
\theta^{(\alpha)}_{ik}}
\end{align*}
with $\widehat{M}_{N}$ the set of infinite, $N$-periodic, strictly upper 
triangular matrices with non-negative integer entries 
which are non-zero only in a finite strip away from the 
diagonal:
\[
(\widehat{M}_{N})^{d}
=
\lbrace \theta=(\theta^{(\alpha)}
)_{\alpha=1}^{d}
\mid \theta^{(\alpha)}_{ik}=\theta^{(\alpha)}_{i+N, k +N} 
\in \Z_{\ge 0} \ 
(i,k \ge 1, \alpha=1,\ldots, d), \ \theta^{(\alpha)}_{ik}
=0 \ 
(k \le i, k \gg i )\rbrace,
\] 
where we set 
$\theta^{(\alpha)}=(\theta^{(\alpha)}_{ik})
_{i,k=1}^{\infty}$ for $\alpha=1, \ldots, d$.

Here $X=(X_{1}, X_{2}, \ldots, )$ and $S=(S_{\alpha})
_{\alpha=1}^{d}$ with  
$S_{\alpha} = (S_{(1, \alpha)}, S_{(2, \alpha)}, 
\ldots, )$ are
$d+1$ types of infinite variables, and
\begin{align*}
c_{N, \infty}(\theta^{(\alpha)}, \theta^{(\beta)} 
|S_{\alpha}, S_{\beta}| q, t )
=
&
\prod_{i=1}^{N} \prod_{i < j \le k < \infty }
{(q^{ \sum_{a>k} \theta^{(\alpha)}_{ia} - 
\theta^{(\beta)}_{ja} } 
t S_{(j,\beta)}/ S_{(i,\alpha)} ;q)_{\theta^{(\alpha)}
_{ik}}
\over 
(q^{ \sum_{a>k} \theta^{(\alpha)}_{ia} - 
\theta^{(\beta)}_{ja} }
S_{(j,\beta)}/ S_{(i,\alpha)} ; q)_{\theta^{(\alpha)}
_{ik}}}
\\
& 
\cdot
\prod_{i=1}^{N} \prod_{i \le j < k < \infty}
{(q^{-\theta^{(\beta)}_{jk} + 
\sum_{a>k} \theta^{(\alpha)}_{ia} - \theta^{(\beta)}_{ja} }
q S_{(j,\beta)}/t S_{(i,\alpha)} ;q)_{\theta^{(\alpha)}
_{ik}}
\over 
(q^{-\theta^{(\beta)}_{jk} + 
\sum_{a>k} \theta^{(\alpha)}_{ia} - 
\theta^{(\beta)}_{ja} } 
S_{(j,\beta)}/S_{(i,\alpha)} ; q)_{\theta^{(\alpha)}
_{ik}}}.
\end{align*}

We substitute $\bar{x} = (\bar{x}_{1}, \bar{x}_{2}, 
\ldots, ), \ 
\bar{s} = 
\left( \bar{s}_{\alpha} \right)_{\alpha=1}^{d}$ with
$\bar{s}_{\alpha}=(\bar{s}_{(1,\alpha)}, 
\bar{s}_{(2,\alpha)}, \ldots, )$ such that
\[
\bar{x}_{1}=x_{1}, \ldots, \bar{x}_{N}=x_{N}, \text{ and } 
\bar{x}_{i+kN}=p^{k} x_{i} \ ( i=1\ldots, N, k=1,2, \ldots,)
\]
\[
\bar{s}_{(1,\alpha)}=s_{(1, \alpha)}, \ldots, 
\bar{s}_{(N, \alpha)}=s_{(N, \alpha)}, \text{ and } 
\bar{s}_{(i+kN, \alpha)}=\kappa^{k} s_{(i, \alpha)} \ 
( i=1\ldots, N, k=1,2, \ldots,).
\]

We set $\delta=(N-1, N-2, \ldots, 1, 0) \in \Z^{N}$.
By \cite[Theorem 3.2]{LNS}, we have the proposition.
\begin{prop}
When $r_{1} = \cdots = r_{N}=d$, we have
\[
f^{\hat{\mk{gl}}_{N}} (p^{\delta/N} \mbi x,p^{1/N}
|\kappa^{\delta/N}\mbi s, \kappa^{1/N}
|q, q/t)
=
f_{N, \infty}(\bar{x}|\bar{s}|q,t).
\]
\end{prop}
\proof 
This is a straight forward generalization of the 
arguments in 
the proof of \cite[Theorem 3.2]{LNS}.
In fact, we set
\begin{align*}
C_{N, \infty}(\mbi \lambda^{(\alpha)}, 
\mbi \lambda^{(\beta)} | S_{\alpha}, S_{\beta} | q,t)
& =
\prod_{i=1}^{N} \prod_{i<j\le k<\infty}
{(q^{\lambda^{(i, \alpha)}_{k-i+1}-
\lambda^{(j, \beta)}_{k-j+1} }
t S_{(j,\beta)}/ S_{(i,\alpha)} ;q)
_{\lambda^{(i, \alpha)}_{k-i} - 
\lambda^{(i, \alpha)}_{k-i+1}}
\over 
(q^{\lambda^{(i, \alpha)}_{k-i+1}-
\lambda^{(j, \beta)}_{k-j+1} }
q S_{(j,\beta)}/ S_{(i,\alpha)} ;q)
_{\lambda^{(i, \alpha)}_{k-i} - 
\lambda^{(i, \alpha)}_{k-i+1}}
}
\\
& \cdot
\prod_{i=1}^{N} \prod_{i \le j < k<\infty}
{(q^{-\lambda^{(j, \beta)}_{k-j} + 
\lambda^{(i, \alpha)}_{k-i+1}}
q S_{(j,\beta)}/ t S_{(i,\alpha)} ;q)
_{\lambda^{(i, \alpha)}_{k-i} - 
\lambda^{(i, \alpha)}_{k-i+1}}
\over 
(q^{-\lambda^{(j, \beta)}_{k-j+1} + 
\lambda^{(i, \alpha)}_{k-i+1}}
S_{(j,\beta)}/ S_{(i,\alpha)} ;q)
_{\lambda^{(i, \alpha)}_{k-i} - 
\lambda^{(i, \alpha)}_{k-i+1}}
}
\end{align*}
where we set $\mbi \lambda^{(\alpha)} =
(\lambda^{(i, \alpha)}_{j})$ and 
$\lambda^{(i+N, \alpha)}_{j} = \lambda^{(i,\alpha)}_{j}$ 
for $\alpha=1,\ldots, d$.
Relations $\theta^{(\alpha)}_{ij}=\lambda^{(i, \alpha)}
_{j-i}$ and $\lambda^{(i,\alpha)}_{k} = 
\sum_{j \ge i+k} \theta^{(\alpha)}_{ij}$ gives a 
bijection of $(\widehat{M})^{d}$ and 
$\prod_{i \in I} \mathsf{P}^{d}$.
Then by \cite[Lemma 3.5]{LNS}, we have
\[
f_{N, \infty}(\bar{x}|\bar{s}|q,t)
=
\sum_{\mbi \lambda \in \prod_{i \in I} \mathsf{P}^{d}}
\prod_{\alpha, \beta=1}^{d}
C_{N, \infty}(\mbi \lambda^{(\alpha)}, \mbi 
\lambda^{(\beta)} 
| \bar{s}_{\alpha}, \bar{s}_{\beta} | q, t) 
\prod_{i=1}^{N} \prod_{ k \ge 1} 
(\bar{x}_{i+k}/ \bar{x}_{i+k-1})^{\sum_{\alpha=1}^{d} 
\lambda^{(i, \alpha)}_{k}}.
\]
On the other hand, we have
\begin{align*}
f^{\hat{\mk{gl}}_{N}}( p^{\delta/N} \mbi x, p^{1/N} 
| \kappa^{\delta/N} \mbi s, \kappa^{1/N} | q, q/t )
&=
\sum_{\mbi{\lambda} \in \prod_{i \in I} 
\mathsf{P}^{d}}
(t/q)^{d \sum_{\alpha=1}^{d} |\mbi \lambda^{(\alpha)}|}
\prod_{\alpha, \beta=1}^{d} 
\wt{C}_{N, \infty}
(\mbi \lambda^{(\alpha)}, \mbi \lambda^{(\beta)} 
| s_{\alpha}, s_{\beta} | q, t, \kappa) 
\\
& \cdot
\prod_{i =1}^{N} \prod_{k \ge 1} 
(\bar{x}_{i + k}/\bar{x}_{i + k -1})
^{\sum_{\alpha=1}^{d} \lambda^{(i, \alpha)}_{k}},
\end{align*}
where $|\mbi \lambda^{(\alpha)}|=
\sum_{i=1}^{N} \sum_{j \ge 1} \lambda^{(i, \alpha)}_{j}$, 
$s_{\alpha}=(s_{(1, \alpha)}, \ldots, s_{(N, \alpha)})$ 
for $\alpha=1,\ldots, d$, and
\[
\wt{C}_{N, \infty}
(\mbi \lambda^{(\alpha)}, \mbi \lambda^{(\beta)} 
| s_{\alpha}, s_{\beta} | q, t, \kappa) 
=\prod_{i,j=1}^{N}
{\Nk^{(j-i|N)}_{\lambda^{(i, \alpha)}, 
\lambda^{(j, \beta)}}
((q/t) (\kappa^{\delta/N} s_{\beta})_{j} 
/ (\kappa^{\delta/N} s_{\alpha})_{i} | q, \kappa^{1/N}) 
\over 
\mathsf{N}_{\lambda^{(i, \alpha)}, 
\lambda^{(j, \beta)}}^{(j-i|N)} 
( (\kappa^{\delta/N} s_{\beta})_{j} 
/ (\kappa^{\delta/N} s_{\alpha})_{i} | q, \kappa^{1/N})}.
\]
By the similar argument as in \cite[Appendix A]{LNS}, 
we have 
\[
(t/q)^{|\mbi \lambda^{\beta}|} \cdot 
\wt{C}_{N, \infty}
(\mbi \lambda^{(\alpha)}, \mbi \lambda^{(\beta)} 
| s_{\alpha}, s_{\beta} | q, t, \kappa) 
=
C_{N, \infty}
(\mbi \lambda^{(\alpha)}, \mbi \lambda^{(\beta)} 
| \bar{s}_{\alpha}, \bar{s}_{\beta} | q, t). 
\]
This completes the proof.
\endproof

\begin{conj}
\label{conj:adj}
(1)
For $\mbi \e_{l} = (\delta_{il})_{i \in I} \in \Z^{I}$, 
we have
\begin{align*}
\widecheck{Z}^{\mbi \e_{l}}_{\adj}(\mbi e, t|\mbi p|q, \kappa)
&=
\prod_{k=1}^{N-1}
{(qp_{l} \cdots p_{l+k-1}; q, t p_{1} \cdots p_{N})_{\infty}
\over
(tp_{l} \cdots p_{l+k-1}; q, t p_{1} \cdots p_{N})_{\infty}}
\\
& \cdot
{(qtp_{1} \cdots p_{N}; q, \kappa^{N}, t p_{1} \cdots p_{N})_{\infty}
(\kappa^{N} tp_1 \cdots p_{N}; q, \kappa^{N}, t p_{1} \cdots p_{N})_{\infty}
\over 
(t^{2}p_1 \cdots p_{N}; q, \kappa^{N}, t p_{1} \cdots p_{N})_{\infty}
(q \kappa^{N}p_1 \cdots p_{N}; q, \kappa^{N}, t p_{1} \cdots p_{N})_{\infty}},
\\
Z^{\mbi \e_{l}}_{\adj}(\mbi e, t|\mbi p|q, \kappa)
&=
\prod_{k=1}^{N-1} 
{(qp_{l-1} \cdots p_{l-k}; q, t p_{1} \cdots p_{N})_{\infty}
\over
(tp_{l-1} \cdots p_{l-k}; q, t p_{1} \cdots p_{N})_{\infty}}
\\
& \cdot
{(qtp_1 \cdots p_{N}; q, \kappa^{N}, t p_{1} \cdots p_{N})_{\infty}
(\kappa^{N} tp_1 \cdots p_{N}; q, \kappa^{N}, t p_{1} \cdots p_{N})_{\infty}
\over 
(t^{2}p_1 \cdots p_{N}; q, \kappa^{N}, t p_{1} \cdots p_{N})_{\infty}
(q \kappa^{N}p_1 \cdots p_{N}; q, \kappa^{N}, t p_{1} \cdots p_{N})_{\infty}}.
\end{align*}
(2)
When $r_{1} = \cdots = r_{N}=d$, we have 
\begin{align*}
Z^{\mbi r}_{\adj}(\mbi e, t|\mbi p|q, \kappa)
=
\widecheck{Z}^{\mbi{r}}_{\adj}(\mbi e, t|\mbi p|q, \kappa).
\end{align*}
\end{conj}

When $N=1$, Conjecture \ref{conj:adj} is already 
checked in 
{\color{black}
\cite{NO} and later in \cite{FOS}}.
See also \cite{CO} for cohomological
version.
\begin{prop}
\label{prop:N=1adj}
When $N=1$, we have
\begin{align*}
\widecheck{Z}^{\mbi \e_{1}}_{{\adj}}
(\mbi e, t|\mbi p|q, \kappa)=
Z^{\mbi \e_{1}}_{\adj}(\mbi e, t|\mbi p|q, \kappa)=
{(qtp_{1}; q, \kappa, t p_{1})_{\infty}
(\kappa tp_{1}; q, \kappa, t p_{1} )_{\infty}
\over 
(t^{2}p_{1}; q, \kappa, t p_{1} )_{\infty}
(q \kappa p_{1}; q, \kappa, t p_{1} )_{\infty}
}.
\end{align*}
\end{prop}
\proof
It follows from \cite[Theorem 3.41]{FOS} after $\kappa 
\leftrightarrow t^{-1}$.
\endproof

\begin{thm}
\label{thm:adj}
We assume Conjecture \ref{conj:adj}.
Then for arbitrary $\mbi r \in \Z^{I}$, we have
\begin{align}
\nonumber
{\widecheck{Z}^{\mbi{r}}_{\adj} (\mbi e, t|\mbi p|q, \kappa)
\over 
Z^{\mbi{r}}_{\adj} (\mbi e, t|\mbi p|q, \kappa)
}
=&
\prod_{k=1}^{N-1}
\prod_{l \in I}
{(qt^{r_{l+1} + \cdots + r_{l+k}}p_{l} \cdots p_{l+k-1}; 
q, t^{\lvert \mbi r \rvert} p_{1} \cdots p_{N}, t)_{\infty}
\over
(t^{r_{l+1} + \cdots + r_{l+k}+1}p_{l} \cdots p_{l+k-1}; 
q, t^{\lvert \mbi r \rvert} p_{1} \cdots p_{N}, t)_{\infty}}
\\
\label{adjrel}
\cdot &
{
(t^{r_{l} + \cdots + r_{l+k-1}+1}p_{l} \cdots p_{l+k-1}; 
q, t^{\lvert \mbi r \rvert} p_{1} \cdots p_{N}, t)_{\infty}
\over
(qt^{r_{l} + \cdots + r_{l+k-1}}p_{l} \cdots p_{l+k-1}; 
q, t^{\lvert \mbi r \rvert} p_{1} \cdots p_{N}, t)_{\infty}
}
\end{align} 
where $\lvert \mbi r \rvert =r_{1} + \cdots + r_{N}$.
\end{thm}

\subsection{Proof of Theorem \ref{thm:adj}}
\label{subsec:check}

For a Young diagram $Y=Y_{\lambda}$ and $k \in I=\Z/N\Z$, set 
\begin{align*}
\lvert Y \rvert_{k} 
&= \lvert \lambda\rvert_{k}
=
\lvert \lbrace s=(i,j) \in Y \mid i \equiv k \mod N 
\rbrace \rvert,
\\
l_{Y}(k)
&=l_{\lambda}(k)
=
\lvert \lbrace s \in Y \mid L_{Y}(s) \equiv k \mod N \rbrace \rvert.
\end{align*}
\begin{lem}
\label{lem:lyk}
For $k=0,1, \ldots, N-1$, we have $\lvert Y \rvert_{k+1}
=l_{Y}(k)$.
\end{lem}
\proof
We consider the one column $Y$ of $L=lN+m$ boxes
with $0 \le m <N$.
We have
\begin{align}
\label{lvert}
\lvert Y \rvert_{k+1}=
\begin{cases} 
l &  m \le k \\
l+1 &  m \ge k+1.
\end{cases}
	\end{align}
On the other hand, we have $l_{Y}(k)=
\lvert \lbrace s=(i,j) \in Y \mid L-i \equiv k \mod N 
\rbrace \rvert$. 
So we label numbers to boxes in the 
one column $Y$
from the bottom beginning with zero.
Then we have the same formula \eqref{lvert} for
$l_{Y}(k)$.
For general $Y$, we can check the assertion by
repeating the similar procedures to each column.
\endproof

For $\mbi \lambda=(\mbi \lambda^{(i)})_{i \in I} \in 
\prod_{i \in I} {\sf P}^{r_{i}}$ with 
$\mbi \lambda^{(i)}=(\lambda^{(i, \alpha)})
_{\alpha=1}^{r_{i}}$ and $k \in I =\Z/N\Z$, set
\begin{align*}
\lvert \mbi \lambda^{(i)} \rvert_{k} = \sum_{\alpha=1}^{r_{i}}
\lvert \lambda^{(i, \alpha)} \rvert_{k}
\end{align*}
for each $i \in I=\Z/N\Z$.
For $\mbi m =(m_{1}, \ldots, m_{N} ) \in \Z^{I}$, set 
$t^{\mbi m} =(t^{m_{1}}, \ldots, t^{m_{N}})$.

\begin{prop}
\label{prop:str}
For $\mbi r=(r_{i})_{i\in I}$ and 
$\mbi \e_{l}=(\delta_{il})
_{i \in I} \in \Z^{I}$ for $l \in I$, we have 
\begin{align*}
\lim_{e_{(l,r_{l}+1)} \to \infty}
Z^{\mbi r+ \mbi \e_{l}}_{\adj}(\tilde{\mbi e}, t|\mbi p|q, \kappa)
&=
Z^{\mbi r}_{\adj}(\mbi e, t|t^{\mbi \e_{l}} \mbi p |q, \kappa)
\cdot
Z^{\mbi \e_{l} }_{\adj}(\mbi e, t|t^{\mbi r[1]} \mbi p|q, \kappa),
\\
\lim_{e_{(l,r_{l}+1)} \to \infty}
\widecheck{Z}^{\mbi r + \mbi \e_{l}}_{\adj}
(\tilde{\mbi e}, t|\mbi p|q, \kappa)
&=
\widecheck{Z}^{\mbi r}_{\adj}
(\mbi e, t|t^{\mbi \e_{l}} \mbi p|q, \kappa)
\cdot
\widecheck{Z}^{\mbi \e_{1}}_{\adj}
(\mbi e, t|t^{\mbi r[1]} \mbi p|q, \kappa),
\end{align*}
where $\mbi e=(\mbi e_{i})_{i \in I}$ and 
$\tilde{\mbi e}=(\tilde{\mbi e}_{i})_{i \in I}$ with
$\tilde{\mbi e}_{l} = (e_{(l,1)}, \ldots, e_{(l, r_{l})},
e_{(l, r_{l}+1)})$ and $\tilde{\mbi e}_{i}=\mbi e_{i}$ 
for $i \neq l$.
\end{prop}
\proof
For each $\mbi \lambda=(\mbi \lambda^{(i)})_{i \in I} \in 
\prod_{i \in I} {\sf P}^{r_{i} + \delta_{il}}$, 
the contribution of the factors
\[
\prod_{i \neq l}
\prod_{\alpha =1 }^{r_{i}}
{{\sf N}^{(l-i)}
_{\lambda_{(i, \alpha)}, \lambda_{(l, r_{l}+1)}}
(te_{(l,r_{l}+1)} /e_{(i, \alpha)} )
\over 
{\sf N}^{(l-i)}
_{\lambda_{(i, \alpha)}, \lambda_{(l, r_{l}+1)}}
(e_{(l,r_{l}+1)} /e_{(i, \alpha)} )}
\cdot
\prod_{\alpha =1 }^{r_{l}}
{{\sf N}^{(0)}_{\lambda_{(l, \alpha)}, \lambda_{(l, r_{l}+1)}}
(te_{(l,r_{l}+1)} /e_{(l, \alpha)} )
\over 
{\sf N}^{(0)}
_{\lambda_{(l, \alpha)}, \lambda_{(l, r_{l}+1)}}
(e_{(l,r_{l}+1)} /e_{(N, \alpha)} )}
\]
is a multiplication of power for $t$ under taking the 
limit
$e_{(l,r_{l}+1)} \to \infty$.
By \eqref{nekfactor} and Lemma \eqref{lem:lyk}, 
{\color{black} the number of power in} $t$ is equal to
\begin{align*}
&
\sum_{i =1}^{N-1} \sum_{\alpha =1 }^{r_{i}} 
\left( 
l_{\lambda^{(i, \alpha)}}(l-i) + l_{\lambda^{(l, r_{l}+1)}}(i-l-1)
\right) + 
\sum_{\alpha =1 }^{r_{l}} 
\left( l_{\lambda^{(l, \alpha)}}(0) + 
l_{\lambda^{(l, r_{l}+1)}}(-1) 
\right)
\\
&=
\sum_{i \neq l}  
\lvert \mbi \lambda^{(i)} \rvert_{l-i+1} + 
\sum_{\alpha =1 }^{r_{l}} 
\lvert \lambda^{(l, \alpha)} \rvert_{1} + 
\sum_{k \neq l} r_{k} \lvert \lambda^{(l, r_{l}+1)} \rvert_{k-l}
+ r_{l}  \lvert \lambda^{(l, r_{l}+1)} \rvert_{0}.  
\end{align*}
On the other hand, the contribution of $\mbi \lambda$
to power of $\mbi p$ in $Z^{\mbi r}_{\adj}$ is equal to
\begin{align*}
\prod_{k=1}^{N} p_{k}^{\sum_{i=1}^{N} |\mbi \lambda^{(i)}|_{k-i+1}}
&=
\prod_{k=1}^{N} p_{k}^{\sum_{i \neq l} |\mbi \lambda^{(i)}|_{k-i+1}
+ \sum_{\alpha=1}^{r_{l} } \lvert \lambda^{(l, \alpha)} \rvert_{k-l+1}}
\\
&\cdot
\prod_{k=1}^{N} p_{k}^{\lvert \lambda^{(l, r_{l}+1)} \rvert_{k-l+1}}.
\end{align*}
Combining both contributions, we get
\[
(tp_{l})^{\sum_{i \neq l} |\mbi \lambda^{(i)}|_{l-i+1}
+ \sum_{\alpha=1}^{r_{l} } \lvert \lambda^{(l, \alpha)} \rvert_{1}}
\cdot
\prod_{k =1}^{N} (t^{r_{k+1}}p_{k})
^{\lvert \lambda^{(l, r_{l}+1)} \rvert_{k-l+1}},
\]
and hence the assertion for $Z^{\mbi r}_{\adj}$.

For $\widecheck{Z}^{\mbi r}_{\adj}$, we consider power in $t$ 
from the contribution of the factors
\[
\prod_{i \neq l}
\prod_{\alpha =1 }^{r_{i}}
{{\sf N}^{(l-i)}_{ \lambda_{(l, r_{l}+1)}, 
\lambda_{(i, \alpha)}}
(te_{(l,r_{l}+1)} /e_{(i, \alpha)} )
\over 
{\sf N}^{(l-i)}_{ \lambda_{(l, r_{l}+1)}, 
\lambda_{(i, \alpha)}}
(e_{(l,r_{l}+1)} /e_{(i, \alpha)} )}
\cdot
\prod_{\alpha =1 }^{r_{l}}
{{\sf N}^{(0)}_{\lambda_{(l, r_{l}+1)}, 
\lambda_{(l, \alpha)}}
(te_{(l,r_{l}+1)} /e_{(l, \alpha)} )
\over 
{\sf N}^{(0)}_{ \lambda_{(l, r_{l}+1)}, 
\lambda_{(l, \alpha)}}
(e_{(l,r_{l}+1)} /e_{(l, \alpha)} )}.
\]
By \eqref{nekfactor} and Lemma \eqref{lem:lyk}, 
the number of power in $t$ is equal to
\begin{align*}
&
\sum_{i\neq l} \sum_{\alpha =1 }^{r_{i}} 
l_{\lambda^{(l, r_{l}+1)}}(l-i) + l_{\lambda^{(i, \alpha)}}(i-l-1) + 
\sum_{\alpha =1 }^{r_{l}} l_{\lambda^{(l, r_{l}+1)}}(0) + 
l_{\lambda^{(l, \alpha)}}(-1) 
\\
&=
\sum_{i \neq l}  
\lvert \mbi \lambda^{(i)} \rvert_{i-l} + 
\sum_{\alpha =1 }^{r_{l}} 
\lvert \lambda^{(l, \alpha)} \rvert_{0} + 
\sum_{k \neq l} r_{k} \lvert \lambda^{(l, r_{l}+1)} \rvert_{l-k+1}
+ r_{l} \lvert \lambda^{(l, r_{l}+1)} \rvert_{1}.  
\end{align*}
On the other hand, the contribution of $\mbi \lambda$
to power of $\mbi p$ in $Z^{\mbi r}_{\adj}$ is equal to
\begin{align*}
\prod_{k=1}^{N} p_{k}^{\sum_{i=1}^{N} |\mbi \lambda^{(i)}|_{i-k}}
&=
\prod_{k=1}^{N} p_{k}^{\sum_{i \neq l} 
|\mbi \lambda^{(i)}|_{i-k}
+ \sum_{\alpha=1}^{r_{l}} \lvert \lambda^{(l, \alpha)} \rvert_{l-k}}
\\
&\cdot
\prod_{k=1}^{N} p_{k}^{\lvert \lambda^{(l, r_{l}+1)} \rvert_{l-k}}.
\end{align*}
Combining both contributions, we get
\[
(tp_{l})^{\sum_{i \neq l} |\mbi \lambda^{(i)}|_{i-l}
+ \sum_{\alpha=1}^{r_{l}} \lvert \lambda^{(l, \alpha)} \rvert_{0}}
\cdot
\prod_{k=1}^{N} (t^{r_{k}}p_{k-1})
^{\lvert \lambda^{(l, r_{l}+1)} \rvert_{l-k+1}},
\]
and hence the assertion for $\widecheck{Z}^{\mbi r}_{\adj}$.
\endproof

Combining Conjecture \ref{conj:adj} and Proposition \ref{prop:str},
we give a proof of Theorem \ref{thm:adj}.
For $\mbi r \in \Z^{I}$, we write by $\Phi^{\mbi r}(\mbi p)$ 
the right hand side of \eqref{adjrel}, that is, set
\begin{align*}
\Phi^{\mbi r}(\mbi p)
&=
\prod_{k=1}^{N-1}
\prod_{l \in I}
{(qt^{r_{l+1} + \cdots + r_{l+k}}p_{l} \cdots p_{l+k-1}; 
q, t^{\lvert \mbi r \rvert} p_{1} \cdots p_{N}, t)_{\infty}
\over
(t^{r_{l+1} + \cdots + r_{l+k}+1}p_{l} \cdots p_{l+k-1};
 q, t^{\lvert \mbi r \rvert} p_{1} \cdots p_{N}, t)_{\infty}}
\\
& \cdot
{
(t^{r_{l} + \cdots + r_{l+k-1}+1}p_{l} \cdots p_{l+k-1}; q, t^{\lvert \mbi r \rvert} p_{1} \cdots p_{N},t)_{\infty}
\over
(qt^{r_{l} + \cdots + r_{l+k-1}}p_{l} \cdots p_{l+k-1}; q, t^{\lvert \mbi r \rvert} p_{1} \cdots p_{N}, t)_{\infty}
}
\end{align*}

We take large enough positive integer $d$ such that 
$d \ge r_{l}$
for any $l \in I$.
By Conjecture \ref{conj:adj} (2), we have 
$\widecheck{Z}^{\sum_{l \in I} d \mbi \e_{l}}_{\adj}=
Z^{\sum_{l \in I} d \mbi \e_{l}}_{\adj}$ which is 
equivalent to 
the equation \eqref{adjrel}.
By the induction on $\sum_{l \in I} (d- r_{l})$, 
we can assume
\eqref{adjrel} holds for $\widecheck{Z}^{\mbi r +
\mbi \e_{l}}
_{\adj} / Z^{\mbi r+ \mbi \e_{l}}_{\adj}$ for some 
$l \in I$.
Applying Proposition \ref{prop:str} for the limit 
$e_{(l, r_{l}+1)} \to \infty$ and 
substituting $p_{l} = t^{-1} p_{l}$, we have
\begin{align*}
{
\widecheck{Z}^{\mbi r}_{\adj}(\mbi p) 
\over 
Z^{\mbi r}_{\adj}(\mbi p)
}
=&
\Phi^{\mbi r+\mbi \e_{l}} (t^{-\mbi \e_{l}} \mbi p)
{Z^{\mbi \e_{l}}_{\adj}(t^{\mbi r[1] - \mbi \e_{l}} \mbi p)
\over 
\widecheck{Z}^{\mbi \e_{l}}_{\adj}(t^{\mbi r[1] - \mbi \e_{l}} \mbi p)} 
\\
=&
\Phi^{\mbi r+\mbi \e_{l}} (t^{-\mbi \e_{l}} \mbi p)
\prod_{k=1}^{N-1} 
{(qt^{r_{l} + \cdots + r_{l-k+1}}p_{l-1} \cdots p_{l-k}; 
q, t^{\lvert \mbi r \rvert} p_{1} \cdots p_{N})_{\infty}
\over
(t^{r_{l} + \cdots + r_{l-k+1}+1}p_{l-1} \cdots p_{l-k};
q, t^{\lvert \mbi r \rvert} p_{1} \cdots p_{N})_{\infty}}
\\
&\cdot
{
(t^{r_{l+1} + \cdots + r_{l+k}}
p_{l} \cdots p_{l+k-1}; q, t^{\lvert \mbi r \rvert} p_{1} \cdots p_{N})_{\infty}
\over
(qt^{r_{l+1} + \cdots + r_{l+k}-1}p_{l} \cdots p_{l+k-1}; 
q, t^{\lvert \mbi r \rvert} p_{1} \cdots p_{N})_{\infty}
}
\\
=&
\Phi^{\mbi r} ( \mbi p).
\end{align*}
This gives a proof of Theorem \ref{thm:adj}.

\subsection{Examples of Conjecture \ref{conj:adj}}
\label{subsec:check}
Here we assume Conjecture \ref{conj:adj}.
For $\ell \in \Z$, we set
\[
F_{\ell}(x) =
{(qx; q, t^{\ell} p_{1} \cdots p_{N})_{\infty}
\over
(tx; q, t^{\ell} p_{1} \cdots p_{N})_{\infty}}.
\]
Then Theorem \ref{thm:adj} implies the following 
expressions.

When $N=2$, we have
\begin{align*}
{
\widecheck{Z}^{(r_{1},r_{2})}_{\adj}(\mbi p)
\over 
Z^{(r_{1}, r_{2})}_{\adj}(\mbi p)
}
=&
\prod_{n=0}^{\infty}
{ 
F_{r_{1} +r_{2}}(t^{r_{2}+n} p_{1} )
 \over 
F_{r_{1} +r_{2}}(t^{r_{1}+n} p_{1} ) 
}
{
F_{r_{1} + r_{2}}(t^{r_{1}+n} p_{2} )
\over  
F_{r_{1} + r_{2}}(t^{r_{2}+n} p_{2} )
}.
\end{align*}
For example, we have
\begin{align*}
{
\widecheck{Z}^{(2,1)}_{\adj}(\mbi p)
\over 
Z^{(2,1)}_{\adj}(\mbi p)
}
=&
{
F_{3}(t p_1)
\over 
F_{3}(t p_2)
}.
\end{align*}

When $N$=3, we have
\begin{align*}
{
\widecheck{Z}^{\mbi r}_{\adj}(\mbi p)
\over 
Z^{\mbi r}_{\adj}(\mbi p)
}
=
\prod_{n=0}^{\infty} &
{F_{\lvert \mbi r\rvert}(t^{r_{2}+n} p_{1} ) 
 \over 
F_{\lvert \mbi r\rvert}(t^{r_{1}+n} p_{1} ) 
} 
{
F_{\lvert \mbi r\rvert}(t^{r_{3}+n} p_{2} )
 \over 
 F_{\lvert \mbi r\rvert}(t^{r_{2}+n} p_{2} )
} 
{
F_{\lvert \mbi r\rvert}(t^{r_{1}+n} p_{3} )
 \over 
F_{\lvert \mbi r\rvert}(t^{r_{3}+n} p_{3} )}
\\
\cdot &
{F_{\lvert \mbi r\rvert}(t^{r_{2}+ r_{3} + n} p_{1} p_{2}) 
 \over 
F_{\lvert \mbi r\rvert}(t^{r_{1} + r_{2} +n} p_{1} p_{2})
}
{
F_{\lvert \mbi r\rvert}(t^{r_{1}+r_{3}+n} p_{2}p_{3} )
 \over 
F_{\lvert \mbi r\rvert}(t^{r_{2}+r_{3}+n} p_{2}p_{3} )
} 
{
F_{\lvert \mbi r\rvert}(t^{r_{1} + r_{2}+n} p_{3}p_{1} ) 
 \over 
F_{\lvert \mbi r\rvert}(t^{r_{1}+r_{3} + n} p_{3} p_{1})
}.
\end{align*}

For example, we have
\begin{align*}
{
\widecheck{Z}^{(2,2,1)}_{\adj}(\mbi p)
\over 
Z^{(2,2,1)}_{\adj}(\mbi p)
}
=&
{
F_{5}(t p_2)
F_{5}(t^3 p_1p_2) 
\over 
F_{5}(t p_3)
F_{5}(t^3 p_1p_3)
},
\\
{
\widecheck{Z}^{(2,1,1)}_{\adj}(\mbi p)
\over 
Z^{(2,1,1)}_{\adj}(\mbi p)
}
=&
{
F_{4}(t p_1)
F_{4}(t^2 p_1p_2)
\over
F_{4}(t p_3) 
F_{4}(t^2 p_2p_3)
},
\\
{
\widecheck{Z}^{(3,1,1)}_{\adj}(\mbi p)
\over 
Z^{(3,1,1)}_{\adj}(\mbi p)
}
=&
{F_{5}(t p_1) 
F_{5}(t^2 p_1)
F_{5}(t^2 p_1p_2 )
F_{5}(t^3 p_1p_2)
\over
F_{5}(t p_3) 
F_{5}(t^2 p_3) 
F_{5} (t^2 p_2p_3)
F_{5}(t^3 p_2p_3 )},
\\
{
\widecheck{Z}^{(3,2,1)}_{\adj}(\mbi p)
\over 
Z^{(3,2,1)}_{\adj}(\mbi p)
}
=&
{
F_{6}(t p_2)
F_{6}(t^2 p_1)
F_{6}(t^3 p_1p_2 ) 
F_{6}(t^4 p_1p_2 ) 
\over 
F_{6}(t p_3)
F_{6}(t^2 p_3)
F_{6}(t^4 p_1p_3 ) 
F_{6}(t^3 p_2p_3 )
}.
\end{align*}

When $N=4$, we have
\begin{align*}
{
\widecheck{Z}^{\mbi r}_{\adj}
\over 
Z^{\mbi r}_{\adj}
}
=
\prod_{n=0}^{\infty} &
{F_{\lvert \mbi r\rvert}(t^{r_{2}+n} p_{1} ) 
 \over 
F_{\lvert \mbi r\rvert}(t^{r_{1}+n} p_{1} ) 
} 
{
F_{\lvert \mbi r\rvert}(t^{r_{3}+n} p_{2} )
 \over 
 F_{\lvert \mbi r\rvert}(t^{r_{2}+n} p_{2} )
} 
{
F_{\lvert \mbi r\rvert}(t^{r_{4}+n} p_{3} )
 \over 
F_{\lvert \mbi r\rvert}(t^{r_{3}+n} p_{3} )}
{
F_{\lvert \mbi r\rvert}(t^{r_{1}+n} p_{4} )
 \over 
F_{\lvert \mbi r\rvert}(t^{r_{4}+n} p_{4} )}
\\
\cdot &
{F_{\lvert \mbi r\rvert}(t^{r_{2}+ r_{3} + n} p_{1} p_{2}) 
 \over 
F_{\lvert \mbi r\rvert}(t^{r_{1} + r_{2} +n} p_{1} p_{2})
}
{
F_{\lvert \mbi r\rvert}(t^{r_{3}+r_{4}+n} p_{2}p_{3} )
 \over 
F_{\lvert \mbi r\rvert}(t^{r_{2}+r_{3}+n} p_{2}p_{3} )
}
{
F_{\lvert \mbi r\rvert}(t^{r_{4} + r_{1}+n} p_{3}p_{4} ) 
 \over 
F_{\lvert \mbi r\rvert}(t^{r_{3}+r_{4} + n} p_{3} p_{4})
}
{
F_{\lvert \mbi r\rvert}(t^{r_{1} + r_{2}+n} p_{4}p_{1} ) 
 \over 
F_{\lvert \mbi r\rvert}(t^{r_{4}+r_{1} + n} p_{4} p_{1})
}
\\
\cdot &
{F_{\lvert \mbi r\rvert}(t^{r_{2}+ r_{3}+r_{4} + n} 
p_{1} p_{2} p_{3}) 
 \over 
F_{\lvert \mbi r\rvert}(t^{r_{1} + r_{2} + r_{3} +n} 
p_{1} p_{2} p_{3})
}
{
F_{\lvert \mbi r\rvert}(t^{r_{3}+r_{4} + r_{1} +n} 
p_{2}p_{3} p_{4})
 \over 
F_{\lvert \mbi r\rvert}(t^{r_{2}+r_{3} + r_{4}+n} 
p_{2}p_{3} p_{4})
}
\\
\cdot &
{
F_{\lvert \mbi r\rvert}(t^{r_{4} + r_{1} + r_{2} +n} 
p_{3}p_{4} p_{1}) 
 \over 
F_{\lvert \mbi r\rvert}(t^{r_{3}+r_{4}+r_{1} + n} 
p_{3} p_{4} p_{1})
}
{
F_{\lvert \mbi r\rvert}(t^{r_{1} + r_{2} + r_{3} +n} 
p_{4}p_{1}p_{2} ) 
 \over 
F_{\lvert \mbi r\rvert}(t^{r_{4}+r_{1} + r_{2} + n} 
p_{4} p_{1} p_{2})
}.
\end{align*}
For example, we have
\begin{align*}
{
\widecheck{Z}^{(2,1,1,1)}_{\adj}
\over 
Z^{(2,1,1,1)}_{\adj}
}
=&
{
F_{5}(t p_1)
F_{5}(t^2 p_1p_2)
F_{5}(t^3 p_1p_2p_3)
\over 
F_{5}(t p_4)
F_{5}(t^2 p_3p_4)
F_{5}(t^3 p_2p_3p_4)
}.
\end{align*}

\subsection{Affine Laumon function}
\label{subsec:affine}
We define $K$-theoretic class $\Lambda_{\fund}$
on $M^{\zeta}(\mbi{r}, \mbi{v})$ 
for a generic stability parameter $\zeta \in \mb R^{I}$  
by
\begin{align}
\label{lfund}
\Lambda_{\fund}
=
\sum_{j \in I} 
\sum_{\beta=1}^{r_{j}} 
\mc V_{j} \otimes \C_{1/\mu_{(j, \beta)}} +
\sum_{i \in I} 
\sum_{\alpha=1}^{r_{i+1} } 
\mc V_{i}^{\vee} \otimes \C_{\nu_{(i+1, \alpha)} \kappa q }.
\end{align}
This class is defined on any moduli stacks with $\mb T$-action
and tautological bundles $\mc V_{i}$ for $i \in I$.

By \eqref{taut-nek} and \eqref{taut+nek},
we have the proposition.
\begin{prop}
\label{prop:fundnek}
(a)
For $\vec{Y}=((Y_{(i, \alpha)})_{\alpha=1}
^{r_{i}} )_{i \in I} 
\in M^{\al}(\mbi{r}, \mbi{v})^{\mb T}$, 
we have
\[
\euk(\Lambda_{\fund})|_{\vec{Y}}
=
\prod_{i, j \in I} \prod_{\alpha=1}^{r_{i}}
\prod_{ \beta=1}^{r_{j}} 
\Nk^{(j-i|N)}_{Y_{(i, \alpha)},\emptyset}
(\mu_{(j, \beta)} / e_{(i, \alpha)})
\cdot
\Nk^{(j-i|N)}_{\emptyset, Y_{(j,\beta)}}
(e_{(j, \beta)} / \nu_{(i,\alpha)}).
\]
(b)
For $\vec{Y}=((Y_{(i, \alpha)})_{\alpha=1}
^{r_{i}} )_{i \in I} 
\in M^{+}(\mbi{r}, \mbi{v})^{\mb T}$, 
we have
\[
\euk(\Lambda_{\fund})|_{\vec{Y}}
=
\prod_{i, j \in I} \prod_{\alpha=1}^{r_{i}}
\prod_{ \beta=1}^{r_{j}} 
\Nk^{(j-i|N)}_{\emptyset, Y_{(i,\alpha)}}
(\mu_{(j,\beta)} / e_{(i, \alpha)})
\cdot
\Nk^{(j-i|N)}_{ Y_{(j,\beta)}, \emptyset}
(e_{(j,\beta)} / \nu_{(i,\alpha)}).
\]
\end{prop}

Set 
\begin{align}
\nonumber
Z^{\mbi{r}}_{\fund} (\mbi e, \mbi \mu, \mbi \nu|\mbi p|q,\kappa)
&=
\sum_{\mbi v \in (\Z_{\ge 0})^{N}}
\mbi p^{\mbi{v}} 
\int_{M^{\al}(\mbi{r}, \mbi{v})} \Lambda_{\fund},
\\
\nonumber
\widecheck{Z}^{\mbi{r}}_{\fund} (\mbi e, \mbi \mu, \mbi \nu
|\mbi p|q,\kappa)
&=
\sum_{\mbi v \in (\Z_{\ge 0})^{N}}
\mbi p^{\mbi{v}} 
\int_{M^{\cal}(\mbi{r}, \mbi{v})} \Lambda_{\fund}.
\end{align}
\begin{prop}
\label{prop:duality}
We have 
\begin{align}
\nonumber
\widecheck{Z}^{\mbi{r}}_{\fund}
(\mbi e, \mbi \mu,\mbi \nu|\mbi p| q, \kappa)
=
Z^{\mbi{r}^{\vee}}_{\fund}
((\mbi e^{\vee})^{-1}, (\mbi \nu^{\vee})^{-1}, 
( \mbi \mu^{\vee})^{-1}| (\mbi p[-1])^{\vee}| q, \kappa).
\end{align}
\end{prop}
\proof
We recall that the isomorphism $\Phi \colon M^{\cal}(\mbi r, \mbi v)
\cong M^{\cal}(\mbi r^{\vee}, \mbi v^{\vee}[1])$ satisfies
$\Phi^{\ast} \mc V_{-i-1} \cong 
\mc V_{i}^{\vee} \otimes \C_{\kappa q}$, and is $\mb T$-equivariant
after taking the homomorphism \eqref{duality} of $\mb T$.
Applying these properties to $\Lambda_{\fund} \in K_{\mb T}
(M^{\al}(\mbi r, \mbi v^{\vee}[1]))$, we have 
\begin{align}
\nonumber
\Phi^{\ast} \Lambda_{\fund} 
&=
\Phi^{\ast} 
\left(
\sum_{j \in I} 
\sum_{\beta=1}^{r_{-j-1}} 
\mc V_{-j-1} \otimes \C_{-1/\mu_{(-j-1, \beta)}} +
\sum_{i \in I} 
\sum_{\alpha=1}^{r_{-i} } 
\mc V_{-i-1}^{\vee} \otimes \C_{-\nu_{(-i, \alpha)}/\kappa q }
\right)
\\
\nonumber
&=
\sum_{j \in I} 
\sum_{\beta=1}^{r^{\vee}_{j+1}} 
\mc V_{j}^{\vee} \otimes \C_{-\nu_{(j+1, \beta)}/\kappa q} +
\sum_{i \in I} 
\sum_{\alpha=1}^{r^{\vee}_{i} } 
\mc V_{i} \otimes \C_{-1/\mu_{(i, \alpha)} }
= \Lambda_{\fund} \in K(M^{\cal}(\mbi{r}, \mbi{v})).
\end{align}
Hence we have
\begin{align}
\nonumber
\widecheck{Z}^{\mbi{r}}_{\fund}
(\mbi e, \mbi \mu, \mbi \nu| 
\mbi p | q, \kappa )
&=
\sum_{\mbi{v} \in \Z_{\ge 0}^{I}} \mbi p^{\mbi{v}}
\left.
\int_{M^{\al}(\mbi{r}^{\vee}, \mbi{v}^{\vee}[1])} 
\Lambda_{\fund}
\right|_{\substack{\mbi e = {}^{t} \mbi e^{-1}
\\ \mbi \mu = {}^{t} \mbi \nu^{-1}
\\ \mbi \nu = {}^{t} \mbi \mu^{-1}}}
\\
\nonumber
&=
\sum_{\mbi{v} \in \Z_{\ge 0}^{I}} 
\mbi p^{(\mbi{v}[-1] )^{\vee}}
\left.
\int_{M^{\al}(\mbi{r}^{\vee}, \mbi{v})} 
\Lambda_{\fund}
\right|_{\substack{\mbi e = {}^{t} \mbi e^{-1}
\\ \mbi \mu = {}^{t} \mbi \nu^{-1}
\\ \mbi \nu = {}^{t} \mbi \mu^{-1}
}},
\end{align}
and the assertion follows.
\endproof

For parameters
$\mbi a=
\left(
(a_{(i,\alpha)})_{\alpha=1}^{r_{i}}
\right)_{i \in I}$, 
$\mbi b=\left(
(b_{(i,\alpha)})_{\alpha=1}^{r_{i}}
\right)_{i \in I}$, 
and $\mbi c=\left(
(c_{(i,\alpha)})_{\alpha=1}^{r_{i}}
\right)_{i \in I}$, and
variables $\mbi p=(p_1,\ldots, p_N)$, set
\begin{align}
\nonumber
f\left( \left.\left.\begin{array}{c}
\mbi a \\\mbi b\\ \mbi c\end{array} \right|
\mbi p \right|q,\kappa \right) 
&=
\sum_{\mbi{\lambda} }
\prod_{i,j \in I}\prod_{\alpha=1}^{r_{i}} 
\prod_{\beta=1}^{r_{j}} 
{
\Nk^{(j-i|N)}_{\emptyset, Y_{(j,\beta)}}
(a_{(i,\alpha)} / b_{(j, \beta)}) 
\cdot
\Nk^{(j-i|N)}_{Y_{(i, \alpha)},\emptyset}
(  b_{(i, \alpha)} / c_{(j, \beta)})
\over 
\mathsf{N}_{Y_{(i,\alpha)}, Y_{(j,\beta)}}^{(j-i|N)} 
( b_{(i,\alpha)}/ b_{(j,\beta)} | q, \kappa)}
\\
\label{gencauchy}
& \cdot \prod_{i,j \in I} \prod_{\alpha=1}^{r_{i}}
p_{j}^{\sum|\lambda^{(i, \alpha)}|_{j-i+1}}, 
\end{align}
where $\mbi{\lambda}=
\left(
(\lambda_{(i, \alpha)})_{\alpha=1}^{r_{i}}
\right)_{i \in I}$ runs in 
$\prod_{i \in I} \mathsf{P}^{r_{i}}$ for
the set $\mathsf{P}$ of partitions, and 
$\emptyset$ denotes the empty partition.
For a partition $\lambda$, we set
\[
|\lambda|_{k} = \sum_{l \in \Z} \lambda_{k+lN}.
\]
By Proposition \ref{prop:fundnek}, we have the proposition.
\begin{prop}
\label{prop:f}
(a) We have 
\begin{align}
\nonumber
Z^{\mbi{r}}_{\fund}
(\mbi e, \mbi \mu, \mbi \nu|\mbi p| q, \kappa)
&=
f\left( \left.\left.\begin{array}{c}
\mbi \nu^{-1} \\\mbi e^{-1} \\ \mbi \mu^{-1} \end{array} \right|
\mbi p \right|q,\kappa \right). 
\end{align}
(b) We have 
\begin{align}
\nonumber
\widecheck{Z}^{\mbi{r}}_{\fund}
(\mbi e, \mbi \mu, \mbi \nu|\mbi p| q, \kappa)
&=
f\left( \left.\left.\begin{array}{c}
\mbi \mu^{\vee} 
\\ 
\mbi e^{\vee} 
\\ 
\mbi \nu^{\vee} \end{array} \right|
(\mbi p [-1])^{\vee} \right|q,\kappa \right).
\end{align}
\end{prop}
\proof
We have to prove only (b). 
It follows from (a) and Proposition \ref{prop:duality}.
\endproof

To propose a conjecture, 
set
\begin{align}
\label{param}
e_{i} &=e_{(i,1)} e_{(i,2)} \cdots e_{(i,r_{i})},
&
\mu_{i}&= \mu_{(i,1)} \mu_{(i,2)} \cdots \mu_{(i,r_{i})},
&
\nu_{i}&= \nu_{(i,1)} \nu_{(i,2)} \cdots \nu_{(i,r_{i})},
\end{align}
and $\vartheta_{i}= p_{i} \partial/ \partial p_{i}$
for $i \in I=\mathbb{Z} /N \mathbb{Z}$, and
\begin{align}
\label{qlaplacian}
\Delta=\vartheta_1^2+\vartheta_2^2+ \cdots +
\vartheta_{N}^2
- \vartheta_1\vartheta_2-\vartheta_2\vartheta_3- \cdots -
\vartheta_{N}\vartheta_1.
\end{align}

\begin{conj}
\label{conj:fund}
We have
\begin{align}
\nonumber
&
\left( {e_1e_2 \cdots e_N  p_1 p_2 \cdots p_N
\over 
\nu_1\nu_2 \cdots \nu_N};\kappa^N \right)_\infty \cdot
\sum_{l_{1} \geq 0} \cdots \sum_{l_{N} \geq 0}
{q^{(l_{1} l_{2} + l_{2} l_{3} + \cdots + l_{N} l_{1})/2}
\over (q;q)_{l_{1}} (q;q)_{l_{2}} \cdots (q;q)_{l_{N}}}
\\
\nonumber
& \cdot
\left(-q^{1/2} e_2 p_1\over \nu_2\right)^{l_{1}} 
\left(-q^{1/2} e_3 p_2\over \nu_3\right)^{l_{2}}
\cdots
\left(-q^{1/2} e_1 p_{N} \over \nu_1\right)^{l_{N}} 
\\
\label{fundrel}
& \cdot
q^{\sum_{i \in I} (l_{i+1}- l_{i-1})\vartheta_i/2}
q^{-{1\over 2}\Delta}
\widecheck{Z}^{\mbi{r}}_{\fund}
(\mbi e, \mbi \mu, \mbi \nu|\mbi p| q, \kappa)\\
\nonumber
=&
\left({\mu_1\mu_2 \cdots \mu_N  p_1 p_2 \cdots p_N 
\over 
e_1 e_2 \cdots e_N};\kappa^N \right)_\infty \cdot 
\sum_{l_{1} \geq 0} \cdots \sum_{l_{N} \geq 0}
{q^{(l_{1} l_{2} + l_{2} l_{3} + \cdots + l_{N} l_{1} )/2}
\over 
(q;q)_{l_{1}} (q;q)_{l_{2}} \cdots (q;q)_{l_{N}}}
\\
\nonumber
&\cdot
\left(-q^{1/2} \mu_1 p_1\over e_1\right)^{l_{1}} 
\left(-q^{1/2} \mu_2 p_2\over e_2\right)^{l_{2}} 
\cdots
\left(-q^{1/2} \mu_N p_N \over e_N\right)^{l_{N}} 
\\
\nonumber
& \cdot
q^{\sum_{i \in I} (-l_{i+1} + l_{i-1}) \vartheta_i/2}
q^{-{1\over 2}\Delta}
Z^{\mbi{r}}_{\fund}
(\mbi e, \mbi \mu, \mbi \nu|\mbi p| q, \kappa)
.
\end{align}
\end{conj}

When $N=1$, we have
\begin{align*}
{Z}^{r_{1}}_{\fund}
(\mbi e_{1}, \mbi \mu_{1}, \mbi \nu_{1} |p_{1}| q, \kappa)
&=f\left( \begin{array}{c}  
\mbi \nu_{1}^{-1}\\
\mbi e_{1}^{-1}\\
\mbi \mu_{1}^{-1}\end{array} 
\Biggl|p_1 \Biggl| q,\kappa\right),
\\
\widecheck{Z}^{r_{1}}_{\fund}
(\mbi e_{1}, \mbi \mu_{1}, \mbi \nu_{1} |p_{1}| q, \kappa)
&=f\left( \begin{array}{c}  
\mbi \mu_{1} \\
\mbi e_{1} \\
\mbi \nu_{1} \end{array} 
\Biggl|p_1\Biggl| q,\kappa\right).
\end{align*}

\begin{prop}
\label{prop:N=1}
When $N=1$, the conjecture \ref{conj:fund} is equivalent to the 
following equation: 
\begin{align}
\notag
&
\widecheck{Z}^{r_{1}}_{\fund}
(\mbi e_{1}, \mbi \mu_{1}, \mbi \nu_{1} |p_{1}| q, \kappa).
\\
\label{N=1}
&=
{
(q\kappa e_{1} p_1/\nu_{1};q,\kappa)_\infty  
(\mu_{1} p_1/e_{1};q,\kappa)_\infty 
\over 
( e_{1} p_1/\nu_{1};q,\kappa)_\infty  
(q \kappa \mu_{1} p_1/e_{1};q,\kappa)_\infty 
}
\cdot Z^{r_{1}}_{\fund}
(\mbi e_{1}, \mbi \mu_{1}, \mbi \nu_{1} 
|p_{1}| q, \kappa)
\end{align}
\end{prop}
\proof 
This follows from 
\[
{(A;q,\kappa)_{\infty} 
\over (q \kappa A; q, \kappa)_{\infty} }
=
(A; \kappa)_{\infty} (q A; q)_{\infty},
\]
and the $q$-binomial theorem.
\endproof
For cohomology theory, we have functional equations \cite[(6)]{O1} 
corresponding to \eqref{N=1}.
Analoguous ones \cite[Theorem 3.3]{O2} for quiver varieties of type $A^{(1)}_{1}$ are aslo obtained.

When $N=2$, we have
\begin{align}
\nonumber
Z^{\mbi{r}}_{\fund}
(\mbi e, \mbi \mu, \mbi \nu|\mbi p| q, \kappa)
&=
f\left( \left.\left.\begin{array}{c}
\mbi \nu_{1}^{-1}, \mbi \nu_{2}^{-1} 
\\ 
\mbi e_{1}^{-1}, \mbi e_{2}^{-1} 
\\ 
\mbi \mu_{1}^{-1}, \mbi \mu_{2}^{-1} \end{array} \right|
p_{1}, p_{2} \right|q,\kappa \right),
\\
\nonumber
\widecheck{Z}^{\mbi{r}}_{\fund}
(\mbi e, \mbi \mu, \mbi \nu|\mbi p| q, \kappa)
&=
f\left( \left.\left.\begin{array}{c}
\mbi \mu_{1}, \mbi \mu_{2} 
\\ 
\mbi e_{1}, \mbi e_{2} \\ 
\mbi \nu_{1}, \mbi \nu_{2} 
\end{array} \right|
p_{2}, p_{1} \right|q,\kappa \right).
\end{align}
\begin{prop}
\label{prop:N=2}
When $N=2$, the conjecture \eqref{conj:fund} 
is equivalent
to the following equation:
\begin{align*}
q^{-{1\over 2}\Delta}
\widecheck{Z}^{\mbi{r}}_{\fund}
(\mbi e, \mbi \mu, \mbi \nu|\mbi p| q, \kappa)
&=
{(-q^{1/2} e_{2} p_1/\nu_{2};q)_\infty
\over 
(-q^{1/2} \mu_{1} p_1/e_{1};q)_\infty
}
{(-q^{1/2} e_{1} p_2/\nu_{1};q)_\infty
\over 
(-q^{1/2} \mu_{2} p_2/e_{2};q)_\infty
}\\
&\cdot
{(q\kappa^2 e_{1}e_{2} p_1 p_2/\nu_{1}\nu_{2};q,\kappa^2)_\infty  
(\mu_{1}\mu_{2} p_1 p_2/e_{1}e_{2};q,\kappa^2)_\infty 
\over 
( e_{1}e_{2} p_1 p_2/\nu_{1}\nu_{2};q,\kappa^2)_\infty  
(q \kappa^2 \mu_{1} \mu_{2} p_1 p_2/e_{1} e_{2};q,\kappa^2)_\infty 
}
q^{-{1\over 2}\Delta}
Z^{\mbi{r}}_{\fund}
(\mbi e, \mbi \mu, \mbi \nu|\mbi p| q, \kappa).
\end{align*}
\end{prop}
\proof
This follows from 
\[
{(AB;q)_{\infty}
\over 
(A;q)_{\infty} (B;q)_{\infty}}
=
\sum_{m,n=0}^{\infty}
{q^{mn} \over (q;q)_{m} (q;q)_{n}
}
A^{m} B^{n},
\]
which is obtained by the two times 
usages of $q$-binomial theorem
and $(B;q)_{m} / (B;q)_{\infty}=1/(q^{m} B;q)_{\infty}$.
\endproof

When $N=3$, we have
\begin{align}
\nonumber
Z^{\mbi{r}}_{\fund}
(\mbi e, \mbi \mu, \mbi \nu|\mbi p| q, \kappa)
&=
f\left( \left.\left.\begin{array}{c}
\mbi \nu_{1}^{-1}, \mbi \nu_{2}^{-1}, \mbi \nu_{3}^{-1} 
\\ 
\mbi e_{1}^{-1}, \mbi e_{2}^{-1}, \mbi e_{3}^{-1} 
\\ 
\mbi \mu_{1}^{-1}, \mbi \mu_{2}^{-1}, \mbi \mu_{3}^{-1} 
\end{array} \right|
p_{1}, p_{2}, p_{3} \right|q,\kappa \right),
\\
\nonumber
\widecheck{Z}^{\mbi{r}}_{\fund}
(\mbi e, \mbi \mu, \mbi \nu|\mbi x| q, \kappa)
&=
f\left( \left.\left.
\begin{array}{c}
\mbi \mu_{2}, \mbi \mu_{1}, \mbi \mu_{3} 
\\ 
\mbi e_{2}, \mbi e_{1}, \mbi e_{3} 
\\ 
\mbi \nu_{2}, \mbi \nu_{1}, \mbi \nu_{3} 
\end{array} \right|
p_{1}, p_{3}, p_{2} \right|q,\kappa \right).
\end{align}
When $N \ge 3$, the factor $q^{\sum_{i \in I} (l_{i+1}- l_{i-1})
\vartheta_i/2}$ in Conjecture \ref{conj:fund} 
is non-trivial, and prevents us 
from factorize $q^{-{1\over 2}\Delta}
\widecheck{Z}^{\mbi{r}}_{\fund}
/q^{-{1\over 2}\Delta}
Z^{\mbi{r}}_{\fund}
$
unlike as in Proposition \ref{prop:N=1} and Proposition 
\ref{prop:N=2}.

\subsection{Hilbert scheme of
points on the plane ($A^{(1)}_{0}$ and $N=r=1$)}
Here we study the case where 
$N=1$ and $r_{1}=1$.
Then it is known that $M^{\al}
(\mbi r, \mbi v)=
M^{-}(1,n)$ is isomorphic to the moduli space of rank $1$ 
framed torsion free
sheaves on $\PP^{2}$, that is, the Hilbert scheme of $n$ points
on $\C^{2}$.
We set $e=e_{1}=e_{(1,1)}, \mu=
\mu_{1}=\mu_{(1,1)}$, $\nu =
\nu_{1}=\nu_{(1,1)}$, and $p=p_{1}$.

By Proposition \ref{prop:N=1adj}, we have 
\begin{align*}
Z^{\mbi{r}}_{\adj}(e, t|p|q, \kappa)
&=
\sum_{\lambda \in \mathsf{P}}
{\mathsf{N}^{(0|1)}_{\lambda,\lambda} 
(t |q,\kappa) \over 
\mathsf{N}^{(0|1)}_{\lambda,\lambda} 
(1 |q,\kappa)}
p^{|\lambda|}
=
\widecheck{Z}^{\mbi{r}}_{\adj}
(e, t|p|q, \kappa).
\end{align*}
{\color{black}
For the fundamental matter class, we have the following.
\begin{prop}
\label{prop:cauchy}
When $N=1$ and $r_{1}=1$, we have
\begin{align}
\nonumber
f\left( \left.\left.\begin{array}{c}
a \\ b \\ c \end{array} \right|
p \right|q,\kappa \right) \nonumber 
&=
\exp\left(
\sum_{n=1}^{\infty}
{
(1- b^{n}/c^{n}) 
(1- a^{n}/q^{n} \kappa^{n} b^{n})
\over
(1-q^{-n}) (1-\kappa^{-n})} \cdot
{p^{n} \over n}
\right).
\end{align}
\end{prop}
\proof
From the definition \eqref{gencauchy}, 
we have
\begin{align*}
f\left( \left.\left.\begin{array}{c}
a \\ b \\ c \end{array} \right|
p \right|q,\kappa \right) \nonumber 
&=
\sum_{\lambda \in \mathsf{P}}
{
\Nk^{(0|1)}_{\emptyset, \lambda}
(a / b) 
\cdot
\Nk^{(0|1)}_{\lambda,\emptyset}
(  b / c)
\over 
\mathsf{N}_{\lambda, \lambda}^{(0|1)} 
( 1 | q, \kappa)}
p^{|\lambda|}
\\
&=
\sum_{\lambda \in \mathsf{P}}
\prod_{s=(l,m) \in Y_{\lambda}}
{
(1- a \kappa^{-l} q^{-m}/b)
(1- b \kappa^{l-1}q^{m-1} b /c)
\over
(1- \kappa^{L_{Y_{\lambda}}(s)}
q^{-A_{Y_{\lambda}}(s)-1})
(1- \kappa^{-L_{Y_{\lambda}}(s)-1}
q^{A_{Y_{\lambda}}(s)})
}.
\end{align*}
Then we get the assertion by substituting $(q,t)=
(q^{-1}, \kappa^{-1})$, 
$X=(1-b/c), Y=(1- a/q \kappa b)$ into \cite[(1.7)]{N2} 
under the principal specialization \cite[(6.17)]{mac}
of the Macdonald polynomial.
\endproof
By Proposition \ref{prop:f} and Proposition \ref{prop:cauchy}, 
we have
\begin{align}
\nonumber
Z^{\mbi{r}}_{\fund}
(e, \mu, \nu|p| q, \kappa)
&=
f\left( \left.\left.\begin{array}{c}
\nu^{-1} \\e^{-1} \\ \mu^{-1} \end{array} \right|
p \right|q,\kappa \right)
=
\exp\left(
\sum_{n=1}^{\infty}
{
(1-e^{n}/q^{n} \kappa^{n} \nu^{n}) (1-\mu^{n}/e^{n})
\over
(1-q^{-n}) (1-\kappa^{-n})} \cdot
{p^{n} \over n}
\right),
\\
\nonumber
\widecheck{Z}^{\mbi{r}}_{\fund}
(e, \mu, \nu|p| q, \kappa)
&=
f\left( \left.\left.\begin{array}{c}
\mu \\ e \\ \nu \end{array} \right|
p \right|q,\kappa \right)
=
\exp\left(
\sum_{n=1}^{\infty}
{
(1- e^{n}/\nu^{n}) 
(1- \mu^{n}/q^{n} \kappa^{n} e^{n})
\over
(1-q^{-n}) (1-\kappa^{-n})} \cdot
{p^{n} \over n}
\right).
\end{align}
Then Conjecture \ref{conj:fund} 
follows from 
Proposition \ref{prop:N=1}, since
we have
\begin{align*}
{
\widecheck{Z}^{\mbi{r}}_{\fund}
(e, \mu, \nu|p| q, \kappa)
\over
Z^{\mbi{r}}_{\fund}
(e, \mu, \nu|p| q, \kappa)
}
&=
\exp
\left(
-
\sum_{n=1}^{\infty}
{
\kappa^{n} q^{n} e^{n} /\nu^{n} - e^{n}/\nu^{n}
-\kappa^{n} q^{n} \mu^{n}/e^{n} - \mu^{n}/ e^{n}
\over
(1-q^{n}) (1-\kappa^{n})} \cdot
{p^{n} \over n}
\right)
\\
&=
{
(q \kappa e p/\nu;q,\kappa)_{\infty}
( \mu p /e;q,\kappa)_{\infty}
\over
( e p /\nu;q,\kappa)_{\infty}
(q \kappa \mu p/e;q,\kappa)_{\infty}
}.
\end{align*}
}

\subsection{$A_{N-1}$ Limit}
When $v_{i}=0$ for some $i \in I$, 
the framed quiver moduli space 
$M^{\zeta}(\mbi{r}, \mbi{v})$ is called the
{\it handsaw quiver variety} of type $A_{N-1}$.
In particular, $M^{\al}(\mbi{r}, \mbi{v})$ with $v_{i}=0$
is called the {\it Laumon space}.
Generating series of integrals of $\Lambda_{\adj}$ and 
$\Lambda_{\fund}$ over Laumon spaces with the fixed framing
$\mbi{r}$ are obtained by substituting $p_{i}=0$ into
$Z^{\mbi{r}}_{\adj}(\mbi e, t|\mbi p|q, \kappa)$ and 
$Z^{\mbi{r}}_{\fund}(\mbi e, \mbi \mu, \mbi \nu
|\mbi p|q, \kappa)$.
Here we consider the case where $v_{N}=0$, and hence
$p_{N}=0$ for simplicity.

For variables $\mbi x = (x_{i})_{i \in I}$ and 
$\mbi s=(\mbi s_{\alpha})_{\alpha=1}^{d}$ with 
$\mbi s_{\alpha} = (s_{(i, \alpha)})_{i \in I}$ for 
$\alpha=1, \ldots, r_{i}$, set
\begin{align*}
f^{\mk{gl}_N}(\mbi x|\mbi s|q,t)
&=
\sum_{\theta \in (M_{N})^{d}}
\prod_{\alpha, \beta=1}^{d}
c_{N}(\theta^{(\alpha)}, \theta^{(\beta)} |
\mbi s_{\alpha}, \mbi s_{\beta}| q, t) 
\prod_{1 \le i < k \le N} (x_{k}/x_{i})^{
\sum_{\alpha=1}^{d} \theta^{(\alpha)}_{ik}},
\end{align*}
where $M_{N}$ is the set of $N \times N$ strictly upper 
triangular matrices with non-negative integer entries, 
and for $\theta=(\theta^{(\alpha)})_{\alpha=1}^{d} 
\in (M_{N})^{d}$ with $\theta^{(\alpha)}=
(\theta^{(\alpha)}_{i,j} )_{i,j \in I} \in M_{N}$, 
we set
\begin{align*}
c_{N}(\theta^{(\alpha)}, \theta^{(\beta)} |
\mbi s_{\alpha}, \mbi s_{\beta}| q, t )
&=
\prod_{1\le i < j \le k \le N}
{(q^{ \sum_{a>k} \theta^{(\alpha)}_{ia} - 
\theta^{(\beta)}_{ja} } 
t s_{(j,\beta)}/ 
s_{(i,\alpha)} ;q)_{\theta^{(\alpha)}_{ik}}
\over 
(q^{ \sum_{a>k} \theta^{(\alpha)}_{ia} - 
\theta^{(\beta)}_{ja} } 
q s_{(j, \beta)}/
s_{(i, \alpha)} ; q)_{\theta^{(\alpha)}_{ik}}}
\\
& 
\cdot
\prod_{1\le i \le j < k \le N}
\frac{(q^{-\theta^{(\beta)}_{jk} + \sum_{a>k} 
\theta^{(\alpha)}_{ia} - \theta^{(\beta)}_{ja} }
q s_{(j,\beta)}/t s_{(i, \alpha)} ;q)_{\theta^{(\alpha)}_{ik}}}
{
(q^{-\theta^{(\beta)}_{jk} + \sum_{a>k} \theta^{(\alpha)}_{ia} - 
\theta^{(\beta)}_{ja} } 
s_{(j,\beta)}/s_{(i,\alpha)} ; q)_{\theta^{(\alpha)}_{ik}}}.
\end{align*}

\begin{prop}
\label{prop:macdlimit}
(a) 
When $r_{1}=\cdots = r_{N}=d$, we have 
\[
Z^{\mbi r}_{\adj}(\kappa^{\delta}\mbi s, q/t|
x_{2} /tx_{1}, \ldots, x_{N}/tx_{N-1},0  
|q, \kappa)
=
f^{\mk{gl}_{N}}
\left( \left. \left.
\mbi x
\right| 
\mbi s  
\right|
q, t
\right).
\]
(b) 
We have 
\[
\left.
Z^{\mbi r}_{\fund}(\mbi e, \mbi \mu, \mbi \nu
|\mbi p|q, \kappa)
\right|_{x_{N}=0}
=
f
\left( \left. \left.
\begin{array}{ccccc}
\mbi 0 & \mbi \nu_{2}^{-1} & \ldots & 
\mbi \nu_{N-1}^{-1}& \mbi \nu_{N}^{-1}\\
\mbi e_{1}^{-1}& \mbi e_{2}^{-1}& \ldots & 
\mbi e_{N-1}^{-1} & \mbi e_{N}^{-1}\\
\mbi \mu_{1}^{-1}& \mbi \mu_{2}^{-1}& \ldots & 
\mbi \mu_{N-1}^{-1} & \mbi 0
\end{array}
\right| 
p_{1}, \ldots, p_{N-1} , 0 
\right|
q, \kappa
\right).
\]
\end{prop}
\proof
(a) By \cite{BFS} and \cite[Corollary 3.4]{LNS}, we have 
\[
\lim_{p \to 0}
f^{\hat{\mk{gl}}_{N}} (p^{\delta/N} \mbi x,p^{1/N}
|\kappa^{\delta/N}\mbi s, \kappa^{1/N}
|q, q/t)
=f^{\mk{gl}_{N}}
\left( \left. \left.
\mbi x
\right| 
\mbi s  
\right|
q, t
\right).
\]
Hence the assertion follows from \eqref{adjruijs}
substituting $p^{1/N}=p$ and $\kappa^{1/N}= \kappa$
since $p$ and $\kappa$ vanish in the right hand side.
Here we give a direct proof for the convenience.
Taking $p_{N}=0$, the partitions producing non-vanishing
contribution to the summation satisfy 
$\ell(\lambda^{(i, \alpha)}) \le N-i$ for 
$i =1, \ldots, N$ and $1 \le \alpha \le d$.
Hence we can parametrize them by using the set 
$(M_{N})^{d}$ as
$\lambda^{(i, \alpha)}_{k}= 
\sum_{j=i+k}^{N} \theta^{(\alpha)}_{i,j}$.
Namely $\theta^{(\alpha)}_{i,j}= 
\lambda^{(i, \alpha)}_{j-i} - 
\lambda^{(i, \alpha)}_{j-i+1}$.
Then for $1 \le i \le j \le N$, we have 
\begin{align}
\nonumber
\mathsf{N}^{(j-i|N)}_{\lambda^{(i,\alpha)}, 
\lambda^{(j,\beta)}}(u|q, \kappa)
&=
\prod_{1\le l \le l' \le N-i \atop l'=l+j-i}
(uq^{-\lambda^{(j,\beta)}_{l}+\lambda^{(i,\alpha)}
_{l'+1}} 
\kappa^{-l+l'} ;q)_{\lambda^{(i, \alpha)}_{l'} - 
\lambda^{(i, \alpha)}_{l'+1}}
\\
\nonumber
&=
\prod_{k=j+1}^{N} 
(uq^{-\sum_{a=k}^{N} \theta^{(\beta)}_{j,a}+
\sum_{a=k+1}^{N} \theta^{(\alpha)}_{i,a}} \kappa^{j-i} ;q)
_{\theta^{(\alpha)}_{i,k}}
\end{align}
setting $k=l+j$.
For $1 \le j < i \le N$, we have
\begin{align}
\nonumber
\mathsf{N}^{(j-i|N)}_{\lambda^{(i,\alpha)}, 
\lambda^{(j,\beta)}}
(u|q, \kappa)
&=
\prod_{1\le l \le l' \le N-j \atop l'=l+i-j-1}
(uq^{\lambda^{(i,\alpha)}_{l} - \lambda^{(j,\beta)}_{l'}} 
\kappa^{l-l'-1} ;q)_{\lambda^{(j, \beta)}_{l'} - 
\lambda^{(j, \beta)}_{l'+1}}
\\
\nonumber
&=
\prod_{k=i}^{N} 
(uq^{\sum_{a=k+1}^{N} \theta^{(\alpha)}_{i,a} - 
\sum_{a=k}^{N} \theta^{(\beta)}_{j,a}} 
\kappa^{j-i} ;q)_{\theta^{(\beta)}_{j,k}}
\\
\nonumber
&=
\prod_{k=i}^{N} 
A^{\theta^{(\beta)}_{jk}} 
q^{\theta^{(\alpha)}_{i,k} 
(\theta^{(\alpha)}_{i,k}-1)/2 }
(uq^{\sum_{a=k+1}^{N} \theta^{(\beta)}_{j,a}-
\theta^{(\alpha)}_{i,a}} q \kappa^{i-j} ;q)
_{\theta^{(\alpha)}_{i,k}}
\end{align}
setting $k=i+l-1$ and 
$A=uq^{\sum_{a=k+1}^{N} \theta^{(\alpha)}_{i,a} - 
\sum_{a=k}^{N} \theta^{(\beta)}_{j,a}} 
\kappa^{j-i}$.
Thus we obtain
\begin{align}
\nonumber
&
\prod_{\alpha, \beta=1}^{d} 
\prod_{i,j \in I}
{
\Nk^{(j-i|N)}_{\lambda^{(i, \alpha)}, \lambda^{(j, \beta)}}
(\kappa^{i-j} q s_{(j, \beta)} / t s_{(i, \alpha)} | q, \kappa) 
\over 
\mathsf{N}_{\lambda^{(i, \alpha)}, \lambda^{(j, \beta)}}^{(j-i|N)} 
( \kappa^{i-j} s_{(j, \beta)} / s_{(i, \alpha)} | q, \kappa)}
\prod_{\alpha=1}^{d} \prod_{i =1}^{N} \prod_{k = 1}^{N-i}
( x_{i + k}/ t x_{i + k -1})
^{\lambda^{(i, \alpha)}_{k}}
\\
\nonumber
&
=
\prod_{\alpha, \beta=1}^{d}
c_{N}(\theta^{(\alpha)}, \theta^{(\beta)} | s_{\alpha}, 
s_{\beta} | q, t)
\prod_{1\le i < j \le N}
(x_{j} / x_{i})^{\sum_{\alpha=1}^{d} 
\theta^{(\alpha)}_{i,j}}.
\end{align}

(b) This directly follows from the definition 
\eqref{lfund}.
\endproof


When $N=2$, by Proposition \ref{prop:macdlimit} (b), 
we have 
\begin{align}
\label{a1limit1}
Z^{\mbi r}_{\fund}
(\mbi e, \mbi \mu, \mbi \nu |p_{1}, 0| q, \kappa)
&=
f\left( \left. \left.
\begin{array}{cc}
0 & \mbi \nu_{2}^{-1} \\
\mbi e_{1}^{-1} & \mbi e_{2}^{-1} \\
\mbi \mu_{1}^{-1} & 0 
\end{array}
\right|
p_{1}, 0
\right|
q, 1
\right),
\\
\label{a1limit2}
\widecheck{Z}^{\mbi{r}}_{\fund}
(\mbi e, \mbi \mu, \mbi \nu |p_{1}, 0| q, \kappa)
&=
f\left( \left. \left.
\begin{array}{cc}
\mbi \mu_{1}^{-1} & 0 \\
\mbi e_{1}^{-1} & \mbi e_{2}^{-1} \\
0 & \mbi \nu_{2}^{-1} 
\end{array}
\right|
p_{1}, 0
\right|
q, 1
\right),
\\
\label{a1limit3}
Z^{\mbi r}_{\fund}
(\mbi e, \mbi \mu, \mbi \nu |0, p_{2}| q, \kappa)
&=
f\left( \left. \left.
\begin{array}{cc}
\mbi \nu_{1}^{-1} & 0 \\
\mbi e_{1}^{-1} & \mbi e_{2}^{-1} \\
0 & \mbi \mu_{2}^{-1}
\end{array}
\right|
0, p_{2}
\right|
q, 1
\right),
\\
\label{a1limit4}
\widecheck{Z}^{\mbi{r}}_{\fund}
(\mbi e, \mbi \mu, \mbi \nu |0, p_{2}| q, \kappa)
&=
f\left( \left. \left.
\begin{array}{cc}
0 & \mbi \mu_{1}^{-1} \\
\mbi e_{1}^{-1} & \mbi e_{2}^{-1} \\
\mbi \nu_{1}^{-1} & 0 
\end{array}
\right|
0, p_{2}
\right|
q, 1
\right).
\end{align}

\section{Handsaw quiver variety of type $A_{1}$}
In the rest of the paper, 
we consider the case $N=2$ and $v_{1}=0$.
We use the indices $0, 1$ in $I=\Z/2 \Z$ instead of $1,2$ 
as in the previous section. 
We set $M^{-} (\mbi{r}, n)=M^{\al} (\mbi r, \mbi v)$, 
$M^{+} (\mbi{r}, n)=M^{\cal} (\mbi r, \mbi v)$ for
$\mbi r=(r_{0}, r_{1})$ and $\mbi v=(v_{0}, v_{1})=
(n,0)$.
They are called the {\it handsaw quiver varieties} 
of type $A_{1}$.

In other words, these spaces are regarded as framed 
quiver moduli for the framed quiver $Q \colon$
\begin{center}
\includegraphics[scale=1]{handsaw}
\end{center}
where we have only one vertex $\ast$ other than $\infty$, 
and $r_{0}$ framings and $r_{1}$ co-framings.
Hence a dimension vector $\alpha=(1,n) \in \Z^{Q_{0}}$ 
are prescribed by a non-negative integer $n$. 
We choose stability parameters $\zeta^{-} < 0 <
\zeta^{+}$ so that we have $M^{\pm} (\mbi{r}, n)= 
M_{Q}^{\zeta^{\pm}}(\alpha)$. 

Furthermore we use
$J_{0} = \lbrace 1, \ldots, r_{0} \rbrace$, $J_{1}=
\lbrace r_{0}+1, \ldots, r=r_{0} + r_{1} \rbrace$, 
and set 
\begin{align}
\label{cvar1}
e_{\alpha} = 
\begin{cases}
e_{(0, \alpha)} & \alpha \in J_{0}\\ 
e_{(1, \alpha- r_{0})} & \alpha \in J_{1}
\end{cases}
\end{align}
so that we have
$W_{0} = \bigoplus_{\alpha \in J_{0}} \C_{e_{\alpha}}, 
W_{1} = \bigoplus_{\alpha \in J_{1}} \C_{e_{\alpha}}$.

We set $V=V_{0}=\C^{n}$ and $\kappa=1$.
Then in \eqref{alml}, we have $\mb L=0$, 
\begin{align*}
\mb M &= \mb M(W, V) = \End_{\C}(V) \otimes \C_{q} 
\times \Hom_{\C}(W_{0}, V) 
\times \Hom_{\C}(V, W_{1}) \otimes \C_{q}, 
\end{align*}
and get the explicit descriptions
$M^{\pm}(\mbi{r}, n)=
[ \mb M^{\pm} / \GL(V) ]$ and 
$\mc V=
[ \mb M^{\pm} \times V / \GL(V) ]$, 
where 
$\mb M^{\pm}$ is the $\zeta^{\pm}$-stable locus in 
$\mb M$.
Set 
\begin{align}
\label{cvar2}
\mu_{\alpha} =
\begin{cases}
q/\mu_{(0,\alpha)} & \alpha \in J_{0}\\ 
1/\nu_{(1, \alpha -r_{0})} & \alpha \in J_{1}.
\end{cases}
\end{align}
Then in \eqref{ladj} and \eqref{lfund}, we have
\begin{align}
\label{cls1}
\Lambda_{\adj} &= 
\mc Hom (\mc W_{0}, \mc V ) + \mc Hom (\mc V, \mc W_{1}) \otimes \C_{q}+
\mc End(\mc V) \otimes \C_{q} - \mc End( \mc V),
\\
\label{cls2}
\Lambda_{\fund} &= 
\bigoplus_{\alpha=1}^{r_{0}} \mc V 
\otimes \mb{C}_{\mu_{\alpha} / q}
+
\bigoplus_{\alpha=r_{0}+1}^{r} \mc V^{\vee} 
\otimes \mb{C}_{q/\mu_{\alpha}},
\end{align}
where $\mc W_{0} = \mo_{M^{\pm}(\mbi{r}, n)} \otimes W_{0}$ and 
$\mc W_{1}=\mo_{M^{\pm}(\mbi{r}, n)} \otimes W_{1}$.
We set
\[
Z_{\pm \adj}^{\mbi r}(q,t, \mbi e; p)= 
\sum_{n=0}^{\infty} 
p^{n} \intk_{M^{\pm}(\mbi{r}, n)} 
\eukt ( \Lambda_{\adj}), \quad 
Z_{\pm \fund}^{\mbi r}(q, \mbi e, \mbi \mu; p)= 
\sum_{n=0}^{\infty} 
p^{n} \intk_{M^{\pm}(\mbi{r}, n)} 
\euk ( \Lambda_{\fund}).
\] 
We have 
$Z_{+ \adj}^{\mbi r}(q,t, \mbi e; p)= 
\widecheck{Z}_{\adj}^{\mbi r}(\mbi e, t|p,0|q, 1)$, 
$Z_{- \adj}^{\mbi r}(q,t, \mbi e; p)= 
Z_{\adj}^{\mbi r}(\mbi e, t|p,0|q, 1)$ by Proposition 
\ref{adj:nek2}, and 
\begin{align*}
Z_{+ \fund}^{\mbi r}(q, \mbi e, \mbi \mu; p)
&= 
\widecheck{Z}_{\fund}^{\mbi r}( \mbi e, \mbi \mu, \mbi \nu| 
p, 0 |q, 1),
&
Z_{- \fund}^{\mbi r}(q, \mbi e, \mbi \mu; p)
&= 
Z_{\fund}^{\mbi r}( \mbi e, \mbi \mu, \mbi \nu| 
p, 0 |q, 1)
\end{align*}
by \eqref{a1limit3}, \eqref{a1limit4}.
Here we substitute 
\[
\mbi e=((e_{\alpha})_{\alpha \in 
J_{0}}, (e_{\alpha})_{\alpha 
\in J_{1}}), \mbi \mu=((q/\mu_{\alpha})_{\alpha \in J_{0}}, 
(0)_{\alpha \in J_{1}}), \mbi \nu=
((0)_{\alpha \in J_{0}}, 
(1/\mu_{\alpha})_{\alpha \in J_{1}}).
\]

\subsection{Combinatorial description}
\label{subsec:comb}
The fixed points set $M^{-}(\mbi r, n)^{\mb T}$ 
can be identified with 
$
\lbrace
\mbi k \in (\Z_{\ge 0})^{J_{0}} \mid |\mbi k|=n  
\rbrace
$,
where we regard $(\Z_{\ge 0})^{J_{0}}$ as a subset of $(\Z_{\ge 0})^{r}$ and set 
$|\mbi k|=\sum_{\alpha =1}^{r} k_{\alpha}$.
The $\mb T$-fixed point in $M^{-}(\mbi r, n)$ corresponding to $\mbi k$ is described by
\begin{align}
\label{taut-}
V=\bigoplus_{\alpha \in J_{0}} \bigoplus_{i=1}^{k_{\alpha}} \C_{e_{\alpha} q^{-i+1}}, 
\quad
W_{0}=\bigoplus_{\alpha \in J_{0}} \C_{e_{\alpha}}, 
\quad
W_{1}=\bigoplus_{\alpha \in J_{1}} \C_{e_{\alpha}},
\end{align}
\[
B(\C_{e_{\alpha} q^{-i+1}}) =
\begin{cases}
\C_{e_{\alpha} q^{-i}} & \text{ if } 1 \le i < k_{\alpha} \\
0 & \text{ if } i = k_{\alpha} \\
\end{cases},
\quad
z(\C_{e_{\alpha}} ) =
\C_{e_{\alpha}},
\]
and $w=0$.
Hence the tangent space $T_{\mbi k}M^{-}(\mbi r, n)$ 
at the point corresponding to $\mbi k$ is
\begin{align}
\notag
&
\Hom_{\C}( V, V) \otimes \C_{q} + \Hom_{\C}(W_{0}, V) + \Hom_{\C}(V, W_{1}) \otimes \C_{q} - \Hom_{\C}(V, V)\\
\label{tan-}
&=
\sum_{\alpha, \beta \in J_{0}} e_{\alpha} e_{\beta}^{-1} 
\left( \sum_{ i =1}^{k_{\alpha}} \sum_{j=1}^{k_{\beta}} ( q^{j -i + 1} - q^{j-i} ) 
+ \sum_{i =1}^{k_{\alpha}}  q^{-i+1} \right)
+ \sum_{\alpha \in J_{1}} \sum_{\beta \in J_{0}} e_{\alpha} e_{\beta}^{-1} \sum_{i=1}^{k_{\beta}} q^{i}.
\end{align}
Here we have calculated \eqref{taut-} and \eqref{tan-} in 
$K_{\mb T}(\pt)=\Z[q^{\pm 1}, t^{\pm 1}, \mbi e^{\pm 1}, \mbi \mu^{\pm 1}]$ where
$q, t, \mbi e=(e_{1}, \ldots, e_{r})$ and 
$\mbi \mu =(\mu_{1}, \ldots, \mu_{r})$ are equivariant parameters.

Similarly, we can identify $M^{+}(\mbi r, n)^{\mb T}$ with 
$
\lbrace
\mbi k \in (\Z_{\ge 0})^{J_{1}} \mid | \mbi k|=n
\rbrace.
$
The $\mb T$-fixed point in $M^{+}(\mbi r, n)$ corresponding to $\mbi k$ is described by
\begin{align}
\label{taut+}
V=\bigoplus_{\alpha \in J_{1}} \bigoplus_{i=1}^{k_{\alpha}} \C_{e_{\alpha} q^{i}}, 
\quad
W_{0}=\bigoplus_{\alpha \in J_{0}} \C_{e_{\alpha}}, 
\quad
W_{1}=\bigoplus_{\alpha \in J_{1}} \C_{e_{\alpha}},
\end{align}
\[
B(\C_{e_{\alpha} q^{i}}) =
\begin{cases}
\C_{e_{\alpha} q^{i-1}} & \text{ if } i > 1 \\
0 & \text{ if } i = 1 \\
\end{cases},
\quad
w(\C_{e_{\alpha} q^{i}}) =
\begin{cases}
0 & \text{ if } i > 1 \\
\C_{e_{\alpha}} & \text{ if } i = 1 \\
\end{cases},
\]
and $z=0$.
The tangent space $T_{\mbi k} M^{+}(\mbi r, n)$ 
at the point corresponding to $\mbi k$ is
\begin{align}
\notag
&
\Hom_{\C}(  V, V) \otimes \C_{q} + \Hom_{\C}(W_{0}, V) + \Hom_{\C}( V, W_{1}) \otimes \C_{q} - \Hom_{\C}(V, V)\\
\label{tan+}
&=
\sum_{\alpha, \beta \in J_{1}} e_{\alpha} e_{\beta}^{-1} 
\left( \sum_{ i =1}^{k_{\alpha}} \sum_{j=1}^{k_{\beta}} ( q^{i-j +1} - q^{i-j} ) 
+ \sum_{i =1}^{k_{\beta}}  q^{-i+1} \right)
+ \sum_{\alpha \in J_{1}} \sum_{\beta \in J_{0}} e_{\alpha} e_{\beta}^{-1} \sum_{i=1}^{k_{\alpha}} q^{i}.
\end{align}
Hence we have 
\begin{align}
\label{stvsco}
T_{\mbi k} M^{+}(\mbi{r}, n) = T_{\mbi k} M^{-}(\mbi{r}, n))|_{e_{1}=e_{1}^{-1}, \ldots, e_{r}=e_{r}^{-1}}.
\end{align}

\subsection{Adjoint matter class $\Lambda_{\adj}$}
\label{subsec:adj}

We take $K$-theory classes \eqref{cls1}
called {\it adjoint matter classes} after physics theory.
Set $H_{n}^{\pm} = \intk_{M^{\pm}(\alpha)} \euk(\Lambda_{\adj})$,
and $Z_{\pm \adj}^{\mbi{r}}(q, \mbi e, \mu, p)
=\sum_{n=0}^{\infty} H_{n}^{\pm} p^{n}$.

\begin{prop}
\label{prop:tangent1}
We have
\begin{align}
T_{\mbi k}M^{+}(\mbi r, n) 
&= 
\sum_{\alpha, \beta \in J_{1}} e_{\alpha}^{-1} e_{\beta}  
\sum_{\ell=k_{\beta} - k_{\alpha}+1}^{k_{\beta}}  q^{\ell}
+
\sum_{\alpha \in J_{0}} \sum_{\beta \in J_{1}} e_{\alpha}^{-1} e_{\beta} 
\sum_{i=1}^{k_{\beta}} q^{i},
\label{tm+}
\\
T_{\mbi k}M^{-}(\mbi r, n)  
&= 
\sum_{\alpha, \beta \in J_{0}} e_{\alpha} e_{\beta}^{-1} \sum_{\ell=k_{\beta} - k_{\alpha}+1}^{k_{\beta}}  q^{\ell}
+
\sum_{\alpha \in J_{1}} \sum_{\beta \in J_{0}} e_{\alpha} e_{\beta}^{-1} 
\sum_{i=1}^{k_{\beta}} q^{i}. 
\label{tm-}
\end{align}
\end{prop}
\proof 
It follows from \eqref{tan-} and \eqref{tan+}.
See \cite[Proposition 2.2]{OY}.
\endproof

\begin{prop}
\label{tangent}
For $\mbi k \in M^{\pm}(\mbi{r}, n)^{\mb T}$, we have
\begin{align}
\label{tangent1}
\frac{
\eukt  (T_{\mbi k} M^{+}(\mbi{r}, n))
}
{
\euk ( T_{\mbi k} M^{+}(\mbi r, n) ) 
}
&=
t^{n r_{0}}
\prod_{ \substack{\alpha \in J_{1}\\ \beta \in J_{1}}}
\frac{(tq^{-k_{\beta}} e_{\alpha} / e_{\beta} ;q)_{k_{\alpha}}}
{(q^{-k_{\beta}} e_{\alpha} / e_{\beta} ;q)_{k_{\alpha}}} 
\prod_{\substack{ \beta \in J_{1} \\ \alpha \in J_{0}}} 
\frac{( qe_{\beta} / t e_{\alpha};q)_{k_{\beta}}}
{( qe_{\beta}/e_{\alpha} ;q)_{k_{\beta}}},
\\
\frac{
\eukt  (T_{\mbi k} M^{-}(\mbi{r}, n))
}
{
\euk ( T_{\mbi k} M^{-}(\mbi r, n) ) 
}
&=
t^{nr_{1}}
\prod_{ \substack{\alpha \in J_{0}\\ \beta \in J_{0}}}
\frac{(tq^{-k_{\beta}} e_{\beta} / e_{\alpha} ;q)_{k_{\alpha}}}
{(q^{-k_{\beta}} e_{\beta} / e_{\alpha} ;q)_{k_{\alpha}}} 
\prod_{\substack{ \beta \in J_{0} \\ \alpha \in J_{1}}} 
\frac{( qe_{\alpha}/t e_{\beta} ;q)_{k_{\beta}}}
{( qe_{\alpha} / e_{\beta};q)_{k_{\beta}}
},
\label{tangent2}
\end{align}
where $(x;q)_{k}=(1- x) ( 1-qx ) \cdots (1-q^{k-1}x )$ 
is the $q$-Pochhammer symbol.
\end{prop}

We rewrite Proposition \ref{tangent} in the form of generating series.
\begin{thm}
\label{explicit1}
We have
\begin{align}
Z_{+\adj}^{\mbi{r}}(q,t, \mbi e; p)
&=
\sum_{ \substack{\mbi k \in \mathbb{Z}_{\ge 0} ^{r_{1}} } }
\left( t^{r_{0}} p  \right)^{\lvert \mbi k \rvert}
\prod_{ \substack{\alpha \in J_{1}\\ \beta \in J_{1}}}
\frac{(tq^{-k_{\beta}} e_{\alpha} / e_{\beta} ;q)_{k_{\alpha}}}
{(q^{-k_{\beta}} e_{\alpha} / e_{\beta} ;q)_{k_{\alpha}}} 
\prod_{\substack{ \beta \in J_{1} \\ \alpha \in J_{0}}} 
\frac{( qe_{\beta}/t e_{\alpha}  ;q)_{k_{\beta}}}
{( qe_{\beta}/e_{\alpha};q)_{k_{\beta}}},
\notag
\\
Z_{-\adj}^{\mbi{r}}(q,t, \mbi e; p)
&=
\sum_{ \substack{\mbi k \in \mathbb{Z}_{\ge 0}^{r_{0}} } }
\left( t^{r_{1}} p \right)^{\lvert \mbi k \rvert}
\prod_{ \substack{\alpha \in J_{0}\\ \beta \in J_{0}}}
\frac{(tq^{-k_{\beta}} e_{\beta} / e_{\alpha} ;q)_{k_{\alpha}}}
{(q^{-k_{\beta}} e_{\beta} / e_{\alpha} ;q)_{k_{\alpha}}} 
\prod_{\substack{ \beta \in J_{0} \\ \alpha \in J_{1}}} 
\frac{(q e_{\alpha}/t e_{\beta} ;q)_{k_{\beta}}}
{( qe_{\alpha}/ e_{\beta} ;q)_{k_{\beta}}
}.
\notag
\end{align}
\end{thm}

By applying Theorem \ref{thm:adjoint}, we give a proof of the 
following theorem in the next section.
\begin{thm}
\label{thm:main1}
We have 
\begin{align}
\label{main1}
Z^{\mbi{r}}_{+ \adj} (q, t, \mbi e; p)
&=
\frac{
(q t^{r_{0}} p ;q,t )_{\infty} ( t^{r_{1}+1} p ;q,t )_{\infty} 
}{
(t^{r_{0}+1 } p ;q, t)_{\infty}(qt^{r_{1} } p ;q, t)_{\infty}
}
\cdot Z_{- \adj}^{\mbi{r}} (q, t, \mbi e; p).
\end{align}
\end{thm}
When $(r_{0}, r_{1})=(1,0)$ we have 
$\sum_{n=0}^{\infty} H_{n}^{+} p^{n} =1$ and 
\begin{align*}
\sum_{n=0}^{\infty} H_{n}^{-} p^{n} 
&=
\sum_{n=0}^{\infty} \frac{(tq^{-n})_{n}}{(q^{-n})_{n}} p^{n}  
=
\sum_{n=0}^{\infty} \frac{(q/t)_{n}}{(q)_{n}} (tp)^{n}  
=
\frac{(qp)_{\infty}}{(tp)_{\infty}}.
\end{align*}
When $(r_{0}, r_{1})=(0,1)$ we have 
$\sum_{n=0}^{\infty} H_{n}^{-} p^{n} =1$ and 
\begin{align*}
\sum_{n=0}^{\infty} H_{n}^{-} p^{n} 
&=
\sum_{n=0}^{\infty} \frac{(tq^{-n})_{n}}{(q^{-n})_{n}} p^{n}  
=
\sum_{n=0}^{\infty} \frac{(q/t)_{n}}{(q)_{n}} (tp)^{n}  
=
\frac{(qp)_{\infty}}{(tp)_{\infty}}.
\end{align*}
In these cases, we can easily check \eqref{main1}.
\begin{lem}
\label{lem:indep}
When $r_{0}=1$, the set 
$\left \lbrace H_{n} \right \rbrace_{n=0}^{\infty}
\subset \C(q,t)( x_{1}, \ldots, x_{r_{1}} )$ is 
linearly independent over $\C(q,t)$, where
$x_{i} = e_{i+1} /e_{1}$ for $i=1, \ldots, r_{1}$.   
\end{lem}
\proof Up to scalar, $H_{n}$ is equal to $f_{n}(x;a,b) = \prod_{i=1}^{r_{1}} \frac{(ax_{i})_{n}}{(b x_{i})_{n}}$
with $a=q/e$ and $b=q$ by Proposition \ref{prop:tangent1}.
It is enough to show $\lbrace H_{n} \rbrace_{n=0}^{\infty}$ is linearly independent over $\C(q,t)( x_{2}, \ldots, x_{r_{1}} )$.
Hence we may assume $r_{1}=1$.  
Then we can show that $\sum_{k=0}^{n} c_{k} f_{k}(x;a,b)=0$ implies $c_{k}=0$ for $k=0,1,\ldots, n$ as follows.

For $q$-difference operator $T_{q, x}$, we have $(T_{q,x} -1) f_{k}(x;a,b)=\frac{x(a-b)(1-q^{k} )}{(bx;q)_{2}} f_{k-1}(x;aq,bq^{2})$.
Hence if we set $\mc D^{(\ell)}=\mc D_{q^{\ell-1}a, q^{2(\ell-1)} b} \cdots \mc D_{q a, q^{2} b} \mc D_{a,b}$ with 
\[
\mc D_{a,b} = \frac{(bx;q)_{2}}{ x (a-b)} \left( T_{q,x} -1 \right), 
\] 
then we have $\mc D^{(\ell)} f_{k}(x;a,b) \Big |_{x=1/aq^{\ell}} = (1;q)_{k} \cdot \delta_{k \ell}$.
This gives the assertion.
\endproof

\subsection{Fundamental matter class $\Lambda_{\fund}$}
\label{subsec:fund}

We take $K$-theory classes \eqref{cls2}
called {\it fundamental matter classes} after physics theory.
Set $I_{n}^{\pm} = \intk_{M^{\pm}(\alpha)} 
\euk(\Lambda_{\fund})$,
and $Z_{\pm \fund}^{\mbi{r}}(q, \mbi e, \mbi \mu; p)
=\sum_{n=0}^{\infty} I_{n}^{\pm} p^{n}$.
By \eqref{taut-}, \eqref{taut+} and Proposition \ref{tangent}, 
we have
\begin{align}
\label{in+}
I^{+}_{n}
&=
\sum_{\mbi k \in (\Z_{\ge 0})^{J_{1}}}
\prod_{ \substack{\alpha \in J_{0} \\ \beta \in J_{1}}} 
{
(q^{- k_{\beta} + 1}/e_{\beta} \mu_{\alpha} ; q)_{k_{\beta}}
\over 
(q^{-k_{\beta}} e_{\alpha} / e_{\beta} ;q)_{k_{\beta}}}
\prod_{ \substack{ \alpha \in J_{1} \\ \beta \in J_{1}} }
{
(e_{\beta} \mu_{\alpha} )_{k_{\beta}}
\over 
(q^{-k_{\beta}} e_{\alpha}/e_{\beta} ; q)_{k_{\alpha}}
}
\\
\label{in-}
I^{-}_{n}
&=
\sum_{\mbi k \in (\Z_{\ge 0})^{J_{0}}}
\prod_{\substack{ \alpha \in J_{1} \\ \beta \in J_{0} }} 
{
( q^{-k_{\beta}} e_{\beta} \mu_{\alpha} )_{k_{\beta}}
\over
(q^{-k_{\beta}} e_{\beta}/e_{\alpha} ; q)_{k_{\beta}}
}
\prod_{\substack{ \alpha \in J_{0} \\ \beta \in J_{0} }} 
{
(q /e_{\beta} \mu_{\alpha} ; q)_{k_{\beta}}
\over 
(q^{-k_{\beta}} e_{\beta} / e_{\alpha} ;q)_{k_{\alpha}}}
\end{align}
\begin{lem}
\label{lem:elem}
We have 
\[
\prod_{\alpha, \beta=1}^{m}
\frac{1}
{(q^{-k_{\beta}}x_{\alpha}/x_{\beta};q)_{k_{\alpha}}}
=
(-1)^{|\mbi k|}q^{\frac{|\mbi k|(|\mbi k| +1)}{2}}
\prod_{1\le \alpha < \beta \le m}
\frac{q^{k_{\alpha}} x_{\alpha} - q^{k_{\beta}} x_{\beta}}
{x_{\alpha} - x_{\beta}}
\prod_{\alpha, \beta=1}^{m}
\frac{1}
{(qx_{\alpha}/x_{\beta};q)_{k_{\alpha}}}
\]
\end{lem}
\proof
For $1 \le \alpha = \beta \le m$, 
the corresponding factor in the left hand side is equal to
\[
(q^{-k_{\beta}})_{k_{\beta}} = (-1)^{k_{\beta}}
q^{\frac{k_{\beta}(k_{\beta}+1)}{2}}
(q)_{k_{\beta}}.
\]
For $\alpha \neq \beta$, we have
\begin{align*}
(q^{-k_{\beta}} e_{\alpha}/e_{\beta}  )_{k_{\alpha}}  \cdot
( q^{-k_{\alpha}} e_{\beta}/e_{\alpha} )_{k_{beta}}
&=
q^{- \frac{ k_{\alpha} (2 k_{\beta} - k_{\alpha}+1)}{2}}
(- e_{\alpha}/e_{\beta} )^{ k_{\alpha} }  
( q^{k_{\beta} - k_{\alpha}+1} e_{\beta}/e_{\alpha} )_{k_{\alpha}}
( q^{-k_{\alpha}} e_{\beta}/e_{\alpha} )_{k_{\beta}}\\
&=
q^{- \frac{ k_{\alpha} (2 k_{\beta} - k_{\alpha}+1)}{2}}
(- e_{\alpha}/e_{\beta} )^{ k_{\alpha} }  
\frac{
(q e_{\beta}/e_{\alpha} )_{k_{\beta} }
(1-e_{\beta}/e_{\alpha})
(q^{-k_{\alpha}} e_{\beta}/e_{\alpha} )_{ k_{\alpha}}
}
{1- q^{k_{\beta}-k_{\alpha}} e_{\beta}/e_{\alpha}}\\
&=
q^{- k_{\alpha} - k_{\beta}  k_{\alpha}}
\frac{
1-e_{\beta}/e_{\alpha} 
}
{1- q^{k_{\beta}-k_{\alpha}} e_{\beta}/e_{\alpha}}
(q e_{\beta}/e_{\alpha} )_{k_{\beta}} (q e_{\alpha}/e_{\beta} )_{k_{\alpha}}
\\
&=
q^{ - k_{\beta}  k_{\alpha}}
\frac{
e_{\alpha} -e_{\beta} 
}
{
q^{k_{\alpha}} e_{\alpha}- q^{k_{\beta}} e_{\beta}
}
(q e_{\beta}/e_{\alpha} )_{k_{\beta}} 
(q e_{\alpha}/e_{\beta} )_{k_{\alpha}}.
\end{align*}
Hence we get the assertion.
\endproof

Using this lemma, we rewrite \eqref{in+} and \eqref{in-}. 
\begin{prop}
\label{taut}
For $\mbi k \in M^{\pm}(\mbi{r},n)^{\mb T}$, we have
\begin{align}
\frac{
\wedge_{-1} ( \Lambda_{\fund})^{\vee}_{\mbi k}
}
{
\wedge_{-1} T_{\mbi k}^{\ast} M^{+}(\mbi r, n) 
}
&=
\left( -q^{ \frac{n+1}{2}} 
\prod_{\alpha \in J_{0}}\frac{q}{ e_{\alpha} \mu_{\alpha} }\right)^{n}
\frac{\Delta_{J_{1}}(q^{\mbi k} \mbi e )}{\Delta_{J_{1}}(\mbi e)}
\prod_{\substack{ \alpha \in J_{0} \\ \beta \in J_{1} }} 
\frac{( e_{\beta} \mu_{\alpha} ;q)_{k_{\beta}}}
{(qe_{\beta}/e_{\alpha} ;q)_{k_{\beta}}}
\prod_{ \substack{\alpha \in J_{1}\\ \beta \in J_{1}}}
\frac{(e_{\beta} \mu_{\alpha} ;q)_{k_{\beta}}}
{(qe_{\beta}/e_{\alpha} ;q)_{k_{\beta}}} 
\label{taut+k}
\\
\frac{
\wedge_{-1} ( \Lambda_{\fund} )^{\vee}_{\mbi k}
}
{
\wedge_{-1} T_{\mbi k}^{\ast} M^{-}(\mbi{r}, n) 
}
&=
\left( - q^{ \frac{n+1}{2}} 
\prod_{\alpha \in J_{1} } e_{\alpha} \mu_{\alpha} \right)^{n}
\frac{\Delta_{J_{0}}(q^{\mbi k} \mbi e^{-1})}{\Delta_{J_{0}}(\mbi e^{-1})}
\prod_{ \substack{ \alpha \in J_{1} \\ \beta \in J_{0} }} 
\frac{
(q/ e_{\beta} \mu_{\alpha} )_{k_{\beta}}
}
{
(qe_{\alpha}/e_{\beta})_{k_{\beta}}
}
\prod_{ \substack{\alpha \in J_{0} \\ \beta \in J_{0}}}
\frac{
(q/ e_{\beta} \mu_{\alpha} )_{k_{\beta}}
}
{
(qe_{\alpha}/e_{\beta})_{k_{\beta}} 
}
\label{taut-k}
\end{align}
where $\Delta_{K}(x)=\prod_{\substack{ \alpha, \beta \in K \\
\alpha < \beta}} (x_{\alpha}-  x_{\beta})$ 
for any subset $K \subset [r]$. 
\end{prop}
\proof

We compare \eqref{in+} and \eqref{taut+k}.
In the first product,  we get the coefficient  
$(q /e_{\alpha} \mu_{\alpha})^{k_{\beta}}$ 
corresponding to $\alpha \in J_{0},  \beta \in J_{1}$.
Applying Lemma \ref{lem:elem} to the second product, 
we get \eqref{taut+k}.
Similarly we compare \eqref{in-} and \eqref{taut-k}, and
use Lemma \ref{lem:elem} to get the assertion.
\endproof

We summarize a result of generating series
$Z_{\pm \fund}^{\mbi{r}}(q, \mbi e, \mu, p)
=\sum_{n=0}^{\infty} I_{n}^{\pm} p^{n}$.
\begin{thm}
\label{subsec:explicit2}
We have
\begin{align}
Z_{+\fund}^{\mbi{r}}(q, \mbi e, \mbi \mu, p)
&= 
\sum_{ \substack{\mbi k \in (\Z_{\ge 0})^{J_{1}} \\ |\mbi k|=n} }  
\left( 
-p q^{\frac{n+1}{2}}
\prod_{\alpha \in J_{0} }\frac{q}{e_{\alpha} \mu_{\alpha}}
\right)^{n}
\frac{\Delta_{J_{1}}(q^{\mbi k} \mbi e )}{\Delta_{J_{1}}(\mbi e)}
\prod_{ \substack{\alpha \in J_{1}\\ \beta \in J_{1}}}
\frac{(e_{\beta} \mu_{\alpha} )_{k_{\beta}}}{(qe_{\beta}
/e_{\alpha})_{k_{\beta}}} 
\prod_{\substack{ \beta \in J_{1} \\ \alpha \in J_{0}}} 
\frac{( e_{\beta} \mu_{\alpha} )_{k_{\beta}}}
{(qe_{\beta}/e_{\alpha})_{k_{\beta}}},
\notag
\\
Z_{- \fund}^{\mbi{r}}(q, \mbi e, \mbi \mu, p)
&= 
\sum_{ \substack{\mbi k \in (\Z_{\ge 0})^{J_{0}} \\ |\mbi k|=n} }  
\left( 
- p  q^{\frac{n+1}{2}}
\prod_{\alpha \in J_{1}} e_{\alpha} \mu_{\alpha}  
\right)^{n}
\frac{\Delta_{J_{0}}(q^{\mbi k} \mbi e^{-1})}{\Delta_{J_{0}}(\mbi e^{-1})}
\prod_{ \substack{\alpha \in J_{0} \\ \beta \in J_{0}}}
\frac{
(q/ e_{\beta} \mu_{\alpha} )_{k_{\beta}}
}
{
(qe_{\alpha}/e_{\beta})_{k_{\beta}} 
}
\prod_{ \substack{\beta \in J_{0} \\ \alpha \in J_{1}}} 
\frac{
(q/ e_{\beta} \mu_{\alpha} )_{k_{\beta}}
}
{
(qe_{\alpha}/e_{\beta})_{k_{\beta}}
}.
\notag
\end{align}
\end{thm}

\subsection{Borel transformation}
\label{subsec:borel}

We set $\hat{I}^{\pm}_{n}=q^{- n(n+1)/2} I^{\pm}_{n}$ 
where 
\[
I^{+}_{n}
=
\intk_{M^{+}(\mbi{r}, n)} C(\mc V),
\quad
I^{-}_{n}
=
\intk_{M^{-}(\mbi{r}, n)} C(\mc V).
\]
These are coefficients of 
$\mc B^{-1} 
Z^{\mbi r}_{\pm \fund}( q, \mbi e, \mbi \mu; p) 
=
\mc B^{-1} \left( 
\sum_{n=0}^{\infty} I^{\pm}_{n} p^{n} \right)$.
Here $\mc B$ is the $q$-Borel transformation acting on 
$p$-series defined by
$\mc B \cdot p^{n} = q^{ n(n+1)/2} p^{n}$ for 
$n \in \Z_{\ge 0}$.

By Proposition \ref{taut}, we have 
\begin{align}
\hat{I}^{+}_{n}
&= 
\sum_{ \substack{\mbi k \in (\Z_{\ge 0})^{J_{1}} \\ |\mbi k|=n} }  
\left( \frac{-1}{ e_{1} \mu_{1} \cdots e_{r_{0}} \mu_{r_{0}} /q^{r_{0}} }\right)^{n}
\frac{\Delta_{J_{1}}(q^{\mbi k} \mbi e )}{\Delta_{J_{1}}(\mbi e)}
\prod_{ \substack{\alpha \in J_{1}\\ \beta \in J_{1}}}
\frac{(e_{\beta} \mu_{\alpha} )_{k_{\beta}}}{(qe_{\beta}/e_{\alpha})_{k_{\beta}}} 
\prod_{\substack{ \beta \in J_{1} \\ \alpha \in J_{0}}} 
\frac{( e_{\beta} \mu_{\alpha} )_{k_{\beta}}}
{(qe_{\beta}/e_{\alpha})_{k_{\beta}}}
\label{hati+}
\\
\hat{I}^{-}_{n}
&= 
\sum_{ \substack{\mbi k \in (\Z_{\ge 0})^{J_{0}} \\ |\mbi k|=n} }  
\left( - e_{r_{0}+1} \mu_{r_{0}+1} \cdots e_{r} \mu_{r} \right)^{n}
\frac{\Delta_{J_{0}}(q^{\mbi k} \mbi e^{-1})}{\Delta_{J_{0}}(\mbi e^{-1})}
\prod_{ \substack{\alpha \in J_{0} \\ \beta \in J_{0}}}
\frac{
(q/ e_{\beta} \mu_{\alpha} )_{k_{\beta}}
}
{
(qe_{\alpha}/e_{\beta})_{k_{\beta}} 
}
\prod_{ \substack{\beta \in J_{0} \\ \alpha \in J_{1}}} 
\frac{
(q/ e_{\beta} \mu_{\alpha} )_{k_{\beta}}
}
{
(qe_{\alpha}/e_{\beta})_{k_{\beta}}
}.
\label{hati-}
\end{align}

\begin{thm}
\label{thm:main2}
We have 
\begin{align}
\label{main2}
\sum_{n=0}^{\infty} \hat{I}^{+}_{n} p^{n}
=
{(- pe_{r_{0}+1} \mu_{r_{0}+1} \cdots e_{r} \mu_{r}
; q)_{\infty} 
\over 
(- q^{r_{0}} p/e_{1}\mu_{1} \cdots e_{r_{0}}\mu_{r_{0}} 
; q)_{\infty}
}
\sum_{n=0}^{\infty} \hat{I}^{-}_{n} p^{n}.
\end{align}
\end{thm}
This is equivalent to the Kajihara transformation formula 
\eqref{kajihara} given in \cite[Theorem 1.1]{K}
as we remarked after Theorem \ref{thm:intromain2}.
In the next section,
we give a proof using wall-crossing formula 
\eqref{main>0}.

Comparing $\mc B$ with $q^{-{\Delta \over 2}}$, 
we substitute $p = q^{1/2}p$ in \eqref{main2} to get 
Proposition \ref{prop:N=2} under $p_{1}= 0, p_{2}=p$ 
after the changes 
\eqref{cvar1} and \eqref{cvar2} of variables.

\section{Wall-crossing formula for handsaw quiver variety}
We consider one of $K$-theory classes \eqref{cls1} and 
\eqref{cls2}.
Using Lemma \ref{lem:residue1}, we deduce explicit 
transformation formula. 

We have $H_{Q}(d^{\sharp}) \cong M((1,0), 
d^{\sharp})$.
To apply results in \S \ref{sec:framed} and \S 
\ref{sec:wall} to the handsaw quiver variety, 
we need lemmas.
\begin{lem}
\label{destab}
\label{lem:hilb}
We have an isomorphism $M^{-}((1, 0), d^{\sharp} ) 
\cong \mb A^{d^{\sharp}}
=\Spec \C[x_{1}, \ldots, x_{d^{\sharp}}]$ 
such that the tautological bundle $\mc V^{\sharp}$ 
corresponds to 
$\C[x_{1}, \ldots, x_{d^{\sharp}}]$-module 
$\C[x_{1}, \ldots, x_{d^{\sharp}}, y] / 
(y^{d^{\sharp}} + x_{1} y^{d^{\sharp}-1} + \cdots + 
x_{d^{\sharp}})$ 
via this isomorphism.
Furthermore $q \in T$ acts by $q^{-1} y$ and 
$q^{i} x_{i}$ for 
$i=1,\ldots, d^{\sharp}$.
\end{lem}

\subsection{Adjoint matter theory}
\label{subsec:euler1}

By Lemma \ref{lem:hilb}, we have 
$\gamma_{d} (t)=t^{d} {(q/t)_{d} \over (q)_{d}}$.
We set 
\[
a_{\mbi{\mk I}}(\mbi{r})=
\prod_{i=1}^{j} 
{ t^{d^{(i)}} [d^{(i)} - 1]_{t}! \over t-1} 
{(q/t)_{d^{(i)}} \over (q)_{d^{(i)}}}
\left( t^{s( \mk I^{(i)}, \mbi{\mk I}^{>i}) + 
r_{1} d^{(i)} - r_{0} d^{(i)} } - 
t^{s( \mbi{\mk I}^{>i}, \mk I^{(i)}) } \right).
\]
For $H^{\pm}_{n}=\int_{M^{\pm}(\mbi r, n)} 
\eukt (\Lambda_{\adj})$, we have
\begin{align}
H^{+}_{n} - H^{-}_{n}  
&=
\sum_{j=1}^{n}
\sum_{\mbi{\mk I} \in \Dec_{j}^{n}} 
{
t^{r_{0} |\mbi d_{\mbi{\mk I}}|}
\cdot} 
{ [ n- \lvert \mbi d_{\mbi{\mk I}} \rvert]_{t} ! \over 
[n]_{t}! } 
a_{\mbi{\mk I}}(\mbi{r}) 
H^{-}_{n- \lvert \mbi d_{\mbi{\mk I}} \rvert}
\label{adj}
\end{align} 
by Theorem \ref{thm:adjoint}.

\begin{lem}
\label{lem:vanish}
When $r_{0} = r_{1}$, we have
\[
\sum_{j=1}^{n}
\sum_{\substack{\mbi{\mk I} \in \Dec_{j}^{n} \\ 
\lvert \mbi d_{\mbi{\mk I}} \rvert=k}} 
{ [ n- \lvert \mbi d_{\mbi{\mk I}} \rvert]_{t} ! \over 
[n]_{t}! }
a_{\mbi{\mk I}}(\mbi{r}) 
=0
\]
\end{lem}
\proof
When $r_{0} = r_{1}$, the left hand side does not depend 
on $(r_{0}, r_{1})$. 
Hence, we can reduce computations to the case where 
$r_{0}=r_{1}=1$.
In this case, we see that $H_{n}^{+}-H_{n}^{-}=0$ in 
\eqref{adj}. 
Since $H^{-}_{0}, \ldots, H^{-}_{n-1}$ is linearly 
independent by Lemma \ref{lem:indep}.
This gives the assertion.
\endproof

By this lemma, we have $Z_{+ \adj}^{(r_{0}, r_{1})}
(q, t, \mbi e,  p ) = Z_{- \adj}^{(r_{0}, r_{1})}
(q, t, \mbi e,  p)$ when $r_{0} =r_{1}$.
From this, we deduce \eqref{main1} by 
Theorem \ref{thm:noumi}.
In fact, substituting $\tau=q/t$, 
$x_{\alpha} = e_{\alpha}$ for $\alpha \in J_{1}$, and 
$y_{\alpha} = q/te_{\alpha}$ for $\alpha \in J_{0}$
in \eqref{noumi2}, we get \eqref{introadj+} and 
\eqref{introadj-}.
When $r_{0} \le r_{1}$, we have 
$\prod_{s=1}^{r_{1}-r_{0}} 
{(qt^{r_{0}+s-1}p ; q)_{\infty} \over 
(t^{r_{0}+s}p;q)_{\infty}}
={
(q t^{r_{0}} p ;q,t )_{\infty} (t^{r_{1}+1} p ;q,t )
_{\infty} 
\over 
(t^{r_{0}+1 } p ;q)_{\infty}(qt^{r_{1} } p ;q, t)_{\infty}
}$,
and the transformation formula 
\eqref{noumi2} in Theorem \ref{thm:noumi}
is equivalent to
\begin{align}
\label{adjfactor}
{Z^{\mbi{r}}_{+ \adj} \over Z^{\mbi{r}}_{- \adj}}= 
\prod_{s=1}^{r_{1}-r_{0}} 
\frac{(qt^{r_{0}+s-1}p ; q)_{\infty}}
{(t^{r_{0}+s}p;q)_{\infty}} 
\end{align}

When $r_{0} \ge r_{1}$, we can replace the role of $x$ 
and $y$ to get the similar formula in 
Theorem \ref{thm:noumi}.
Then we have
$Z^{\mbi{r}}_{- \adj} = \prod_{s=1}^{r_{0}-r_{1}} 
\frac{(qt^{r_{1}+s-1}p ; q)_{\infty}}
{(t^{r_{1}+s}p;q)_{\infty}} Z^{\mbi{r}}_{+ \adj}$.
This also implies \eqref{adjfactor}.
This completes the proof of Theorem \ref{thm:main1}.

\subsection{Fundamental matter theory}
\label{subsec:euler1}

We take $\Lambda=\Lambda_{\fund}$ in \eqref{cls2}.
We set
\[
F(\mc V)= \prod_{\alpha \in J_{0}}
\wedge_{-q/\mu_{\alpha}} \mc V^{\vee}
\cdot
\prod_{\alpha \in J_{1}}
\wedge_{-\mu_{\alpha}/q} \mc V.
\]
Then we have $\ch F(\mc V)=\eukt(\Lambda_{fund})$.
For $\mbi d=(d^{(1)}, \ldots, d^{(j)}) \in 
\Z_{\ge 0}^{j}$, we compute $C_{\mbi d}$ in 
\eqref{itclasst=1} part by part.
%
%
By Lemma \ref{nmki}, we have 
\begin{align}
&
\oint d\hbar_{i}
\frac{ \euk \left( \sum_{\alpha =1}^{r_{0}} \mc V^{(i)} \otimes e^{\hbar_{i}} \mu_{\alpha}/q +
\sum_{\alpha =r_{0}+1}^{r} \mc V^{(i) \vee} \otimes e^{-\hbar_{i}} q/\mu_{\alpha}\right)  
} 
{
\euk \left( 
\mk N ( \mc V^{(i)} \otimes e^{\hbar_{i}}, \mc V \oplus \bigoplus_{l=i+1}^{j} \mc V^{(l)} \otimes e^{\hbar_{l}})
\right).
}
\notag
\\
&=
\oint d\hbar_{i}
{
\ch \wedge_{-1} \left(
\sum_{\alpha =1}^{r_{0}}  \mc V^{(i) \vee} \otimes e^{-\hbar_{i}} q/\mu_{\alpha} 
+\sum_{\alpha =r_{0}+1}^{r}  \mc V^{(i)} \otimes e^{\hbar_{i}} \mu_{\alpha}/q 
\right)
\over 
\ch \wedge_{-1}
\left(
\mc V^{(i) \vee} \otimes  \mc W_{0} e^{-\hbar_{i}}
+\mc V^{(i)} \otimes  \mc W_{1}^{\vee} e^{\hbar_{i}}/q
\right)
} 
\notag
\\
& \cdot
{ 
\ch \wedge_{-1} \left(
\mc V^{(i) \vee} \otimes \mc V^{>i} e^{ - \hbar_{i}} 
+ 
\mc V^{>i \vee} \otimes \mc V^{(i)} e^{\hbar_{i}}
\right)
\over 
\ch \wedge_{-1} \left( \mc V^{(i) \vee} \otimes 
\mc V^{>i} e^{ - \hbar_{i}} / q +
\mc V^{>i \vee} \otimes \mc V^{(i)} e^{\hbar_{i}} / q
\right)
}
\notag\\
&=
{q^{r_{0} d^{(i)}+ ( n- |\mbi d|_{i})d^{(i)} }
\over (e_{1} \mu_{1} \cdots e_{r_{0}} \mu_{r_{0}} )
^{d^{(i)}}}
\oint d\hbar_{i}
\prod_{\alpha =1}^{r_{0}} 
{
(e^{\hbar_{i}} q^{-d^{(i)}} \mu_{\alpha} )_{d^{(i)}} 
\over 
(e^{\hbar_{i}} q^{-d^{(i)}+1}/ e_{\alpha})_{d^{(i)}} 
}
\prod_{\alpha =r_{0}+1}^{r} 
{
(e^{\hbar_{i}} q^{-d^{(i)}} \mu_{\alpha} )_{d^{(i)}} 
\over 
(e^{\hbar_{i}} q^{-d^{(i)}}/ e_{\alpha})_{d^{(i)}} 
}
\notag\\
&=
{q^{( n - |\mbi d|_{i}) d^{(i)}}
\over 
(e_{1} \mu_{1} \cdots e_{r_{0}} \mu_{r_{0}}/q^{r_{0}})
^{d^{(i)}}}
((e_{1} \mu_{1} \cdots e_{r} \mu_{r} / q^{r_{0}})
^{d^{(i)}} - 1)
\label{factor1}
\end{align}
where we set $\mc V^{>i}=\mc V \oplus 
\bigoplus_{k=i+1}^{j} \mc V^{(k)} \otimes 
e^{\hbar_{k}}$ and $|\mbi d|_{i} = d^{(1)} + \cdots + 
d^{(i)}$ for $i=0, \ldots, j$.
Furthermore we have $H_{Q}(d^{(i)})=M^{-}((1,0), 
d^{(i)})$, and 
\begin{align}
\intk_{H_{Q}(d^{(i)})}
\ch \wedge_{-1}
\left(\mc V^{(i)} /\mo_{H_{Q}(d^{(i)})} \right)^{\vee}
&=
{(q)_{d^{(i)}-1} \over (q^{-d^{(i)}})_{d^{(i)}}}
=
{(-1)^{d^{(i)}} q^{d^{(i)}(d^{(i)}+1)/2} 
\over 1 - q^{d^{(i)}}}.
\label{factor2}
\end{align}

For $I_{n}= \int_{M^{\pm}(\mbi r, n)} 
\euk(\Lambda_{\fund})$, we have 
\begin{align}
&
I^{+}_{n} - I^{-}_{n} 
\notag
\\
&=
\sum_{j=1}^{n} 
(-1)^{j}
\sum_{{ \mbi d \in \Z_{>0}^{j} \atop | \mbi d | \le n}} 
\prod_{i=1}^{j} 
{(-1)^{d^{(i)}} \over |\mbi d|_{i}}
{q^{ ( n - |\mbi d|_{i}) d^{(i)} +d^{(i)}(d^{(i)}+1)/2}
\over 1 - q^{d^{(i)}}}
{(e_{1} \mu_{1} \cdots e_{r} \mu_{r} / 
q^{r_{0}})^{d^{(i)}} - 1 
\over 
(e_{1} \mu_{1} \cdots e_{r_{0}} \mu_{r_{0}} 
/q^{r_{0}} )^{d^{(i)}}}
I^{-}_{n-|\mbi d|}
\label{formula3'}
\end{align} 
by \eqref{main>1}.
For $\hat{I}^{\pm}_{n}=q^{ - n(n+1)/2} \cdot 
I^{\pm}_{n}$, 
we can rewrite \eqref{formula3'} as follows :
\begin{align}
\hat{I}^{+}_{n} - \hat{I}^{-}_{n}
&=
q^{-n(n+1)/2}
\sum_{l=1}^{n}
q^{(n-l)(n-l+1)/2}
\hat{I}^{-}_{n- l}
\notag\\
&\cdot
(-1)^{l}
\sum_{\substack{0=c_{0}< c_{1} < \cdots < c_{j} = l }} 
\prod_{i=1}^{j}
\frac{1}{ c_{i}} 
\cdot
\frac{q^{(n - c_{i})(c_{i} - c_{i-1} )+(c_{i}-c_{i-1})(c_{i} -c_{i-1}+1)/2}}{1-q^{c_{i} - c_{i-1} }}
\left( B^{c_{i} - c_{i-1} } - A^{c_{i}-c_{i-1}} \right)
\notag
\\
&=
\sum_{l=1}^{n}
\hat{I}^{-}_{n- l}
\cdot (-B)^{l}
\sum_{\substack{0=c_{0}< c_{1} < \cdots < c_{j} = l }} 
\prod_{i=1}^{j}
\frac{1}{ c_{i}} 
\cdot
\frac{1 - (A/B)^{c_{i}-c_{i-1}}}{1-q^{c_{i} - c_{i-1} }},
\label{check2}
\end{align}
where $B= q^{r_{0}} / e_{1} \mu_{1} \cdots e_{r_{0}} 
\mu_{r_{0}} $, 
and $A= e_{r_{0}+1} \mu_{r_{0}+1} \cdots e_{r} \mu_{r}$.

In \eqref{check2}, the coefficients on $\hat{I}_{n-l}$ 
does not depend on $n$ nor $r$, 
hence it is reduced to the case where $r_{0}=1$ and $r_{1}=0$. 
In this case, we see that $\sum_{n=0} \hat{I}^{+}_{n} p^{n}=1$ and 
$\sum_{n=0} \hat{I}^{-}_{n} p^{n}=(-qp/e_{1}\mu_{1})_{\infty}/(-p)_{\infty}$
by \eqref{hati-}.
Hence we have
\[
\sum_{\substack{0=c_{0}< c_{1} < \cdots < c_{j} = l }} 
\prod_{i=1}^{j}
\frac{1}{ c_{i}} 
\cdot
\frac{1 - (1/B)^{c_{i}-c_{i-1}}}{1-q^{c_{i} - c_{i-1} }}
= 
\frac{(1/B)_{l}}{(q)_{l}}.
\]
Substituting $B=B/A$, we have
\[
\sum_{n=0}^{\infty} \hat{I}^{+}_{n} p^{n} 
=
{(- pA ;q)_{\infty} \over (- pB ;q)_{\infty}}
\sum_{n=0}^{\infty} \hat{I}^{-}_{n} p^{n}.
\]
This gives a proof of \eqref{main2}.
\vspace{.5cm}\\
\begin{ack}
We would like to appreciate Takuro Mochizuki and 
Masatoshi Noumi for useful discussions. 
{\color{black}
The author would like to appreciate referees for careful 
suggestions and corrections.
They made the manuscript greatly improved.
}
The work of R.O. was partly supported by Osaka Central Advanced Mathematical
Institute: MEXT Joint Usage/Research Center on Mathematics and
Theoretical Physics JPMXP0619217849, and the Research Institute for Mathematical Sciences,
an International Joint Usage/Research Center located in Kyoto University.
\end{ack}

\noindent {\bf Funding} 
Partial financial support was received from
Grant-in-Aid for Scientific Research (Kakenhi);
Research 21K03180 (R.O.), 19K03512 (J.S.) 
and 19K03530 (J.S.).
\vspace{.5cm}\\
\noindent {\bf Data Availability}
Data sharing not applicable to this article as no 
datasets were generated or analysed during the current 
study.
\vspace{.5cm}\\
\noindent {\bf Declaration}
\vspace{1cm}\\
\noindent {\bf Conflict of interest}
The authors have no conflict of interest directly 
relevant to the content of this article.




\end{document}